\documentclass[11pt]{amsart}
\usepackage[margin=1in]{geometry}
\usepackage{amsmath,amssymb,amsthm,mathtools}
\usepackage{graphicx}
\usepackage{enumitem}
\usepackage{hyperref}
\usepackage{xcolor}

\hypersetup{
	colorlinks=true,
	linkcolor=blue!60!black,
	citecolor=blue!60!black,
	urlcolor=blue!60!black,
	bookmarksdepth=2
}

\newcommand{\C}{\mathbb{C}}
\newcommand{\R}{\mathbb{R}}
\newcommand{\E}{\mathbb{E}}
\newcommand{\id}{\mathrm{id}}
\newcommand{\scal}{\mathrm{sc}}
\DeclareMathOperator{\supp}{supp}
\DeclareMathOperator{\Var}{Var}
\DeclareMathOperator{\dist}{dist}
\DeclareMathOperator{\per}{per}
\DeclareMathOperator{\Her}{He}

\newtheorem{theorem}{Theorem}[section]
\newtheorem{proposition}[theorem]{Proposition}
\newtheorem{corollary}[theorem]{Corollary}
\newtheorem{lemma}[theorem]{Lemma}
\newtheorem{remark}[theorem]{Remark}

\title{Wigner Matrices and the \texorpdfstring{$\alpha$}{alpha}-Determinant}
\author{Frederic Koehler}
\author{Pui Kuen Leung}
\date{}

\begin{document}
	\begin{abstract}
		For a Hermitian Wigner matrix $M_n$ and complex $\alpha$ with
		$|\alpha|\leq1$, we study the $\alpha$-characteristic polynomial
		$\det_\alpha(zI-M_n)$, which interpolates between the characteristic
		polynomial ($\alpha=-1$) and the permanental polynomial ($\alpha=1$).
		Assuming finite fourth moments off the diagonal and second moments on
		it, we show that this polynomial, normalized by its expectation,
		converges in distribution off the segment $\sqrt{-\alpha}\,[-2,2]$ to
		the Wick exponential of an explicit Gaussian analytic function, with
		the same Gaussian modes for every $\alpha$. As consequences, all
		zeros lie near this segment with probability tending to one, the
		empirical zero measure converges to a rotated semicircle law, and
		analytic linear statistics of the zeros satisfy a central
		limit theorem. We also show that $\det_\alpha(zI-M_n)$ is
		approximated, with vanishing relative error, in terms of the
		characteristic polynomial of $M_n$ at a transformed point. The proof
		writes the $\alpha$-characteristic polynomial as a polymer gas, which
		becomes asymptotically ideal after an exact resummation over matchings.
	\end{abstract}
	\maketitle
	\setcounter{tocdepth}{1}
	
	\tableofcontents
	
	\section{Introduction}\label{sec:intro}
	
	For an $n\times n$ matrix $A \in \C^{n\times n}$ and $\alpha\in\C$, we can define the $\alpha$-determinant (see \cite{VereJones1988,ShiraiTakahashi2003}) by 
	\[
	\det_\alpha(A)
	:=
	\sum_{\sigma\in\mathfrak S_n}\alpha^{n-c(\sigma)}\prod_{i=1}^n A_{i,\sigma(i)},
	\]
	where $c(\sigma)$ is the number of cycles of $\sigma$. Thus $\det_{-1}(A)=\det(A)$ and $\det_1(A)=\per(A)$, where $\per(A)$ is known as the \emph{permanent} of the matrix $A$.	
	
	We will study the $\alpha$-determinant in the context of Wigner random matrices.
	Let $M_n$ be Hermitian, with strictly upper-triangular entries
	$\xi_{ij}/\sqrt n$, where the $\xi_{ij}$, $i<j$, are i.i.d.\ with
	\[
	\E\xi_{12}=0,\qquad \E|\xi_{12}|^2=1,\qquad
	\E\xi_{12}^2=\theta,\qquad \E|\xi_{12}|^4<\infty.
	\]
	The holomorphic second moment $\theta$ may be complex and necessarily
	satisfies $|\theta|\leq1$.
	Its diagonal entries are $\eta_i/\sqrt n$, where the $\eta_i$ are
	i.i.d.\ real random variables, independent of the off-diagonal
	entries, with
	\[
	\E\eta_1=0,\qquad \E\eta_1^2=v<\infty.
	\]
	The $\alpha$-characteristic polynomial of this matrix is
	\[
	p_{n,\alpha}(z):=\det_\alpha(zI-M_n).
	\]
	When $\alpha = -1$, this is the usual characteristic polynomial and the roots are the eigenvalues of $M_n$. The eigenvalues converge (in empirical distribution) to the famous \emph{Wigner semicircle law} and, more generally, have been deeply studied in random matrix theory (see, e.g., \cite{AndersonGuionnetZeitouni2010}).
	When $\alpha = 1$, this is known as the \emph{permanental polynomial} --- see, e.g., \cite{merris1981permanental} for general discussion and background.

	Fyodorov initiated the study of the permanental polynomial in the context of random matrices, and conjectured for Gaussian entries (GOE and GUE) that the zeros converge to a semicircle law along the imaginary axis \cite{Fyodorov2006}. This conjecture was recently proved by Feng and Yao \cite{FengYao2026}, who asked about universality. In this work, we develop a new method and prove a functional limit theorem for the $\alpha$-characteristic polynomial: for general i.i.d.\ entries with finite fourth moments off the diagonal, and jointly in $|\alpha|\leq1$, the polynomial normalized by its expectation converges to the Wick exponential of a single Gaussian field, that is, the exponential of the field normalized to have mean one.

	\subsection{Summary of main result and consequences.} 
	
	Write $\Sigma_\alpha:=\sqrt{-\alpha}\,[-2,2]$, using the principal square root, and
	let $g_\alpha$ be the Cauchy transform of the semicircle law rotated
	onto $\Sigma_\alpha$. Let $\kappa:=\Var(|\xi_{12}|^2)$ and
	suppose for simplicity that $\theta$ is real. Then, we show that jointly in
	$|\alpha|\leq1$ and locally uniformly in $z\notin\Sigma_\alpha$,
	\begin{equation}\label{eq:informal-main}
		\frac{p_{n,\alpha}(z)}{\E p_{n,\alpha}(z)}
		\Longrightarrow
		\exp\Bigl(\sum_{\ell\geq1}(-1)^\ell\alpha^{\ell-1}g_\alpha(z)^\ell G_\ell
		-\frac12\sum_{\ell\geq1}\alpha^{2\ell-2}g_\alpha(z)^{2\ell}\sigma_\ell^2\Bigr),
	\end{equation}
	where the $G_\ell$ are independent $N(0,\sigma_\ell^2)$ variables with
	$\sigma_1^2=v$, $\sigma_2^2=\kappa/2$, and
	$\sigma_\ell^2=(1+\theta^\ell)/\ell$ for $\ell\geq3$.
	The $G_\ell$ do not depend on $\alpha$: for polynomials formed from the
	same matrix, one Gaussian field governs all the limits.
	Theorem~\ref{gaussian} gives the precise statement.
	
	The functional limit has many useful consequences:
	\begin{itemize}
		\item \emph{Zeros.} With probability tending to one, every zero lies
		in any prescribed neighborhood of $\Sigma_\alpha$, simultaneously
		for all $|\alpha|\leq1$ (Theorem~\ref{zerofree}). The empirical
		zero measures converge almost surely to the corresponding rotated semicircle
		laws, uniformly over the same disk (Theorem~\ref{thm:almost-sure-semicircle}).
		\item \emph{Fluctuations.} The empirical Cauchy transform has
		Gaussian fluctuations of order $n^{-1}$ jointly in $(\alpha,z)$
		(Theorem~\ref{thm:cauchy-transform-fluctuations}), and analytic
		linear statistics of the zeros satisfy a central limit theorem
		diagonalized by Dickson polynomials (Theorem~\ref{clt} and
		Corollary~\ref{cor:chebyshev-diagonalization}).
		\item \emph{Approximation.} Although permanents are
		$\#\mathrm{P}$-hard\footnote{More generally, $\det_\alpha$ is $\#\mathrm{P}$-hard to compute for every fixed rational $\alpha\notin\{0,-1\}$, as a specialization of the digraph cover polynomial \cite{BlaserDellFouz2012}.} to compute in the worst case \cite{Valiant1979},
		the polynomial $p_{n,\alpha}(z)$ can be approximated with vanishing relative error for any $z \notin \Sigma_{\alpha}$. The approximation is given by an explicit power of an ordinary determinant of the
		same matrix at a transformed spectral parameter
		(Theorem~\ref{thm:determinantal-approximation}).
		Consequently, the normalized polynomial is
		asymptotically determined by the eigenvalues of $M_n$ (Proposition~\ref{prop:spectral-determination}), 
		even though the $\alpha$-determinant depends on the choice of basis.
	\end{itemize}

	\subsection{Related work.}\label{sec:related-work}
	The general $\alpha$-determinant was introduced by Vere-Jones
	\cite{VereJones1988,VereJones1997}; our definition follows the convention of Shirai and
	Takahashi \cite{ShiraiTakahashi2003} (see also \cite{MatsumotoWakayama2006}).

	As discussed above, Fyodorov \cite{Fyodorov2006} and Feng and Yao \cite{FengYao2026}
	previously studied the semicircle law for the permanental case $\alpha = 1$. Their proof techniques are specific to the Gaussian case.

	In the classical random matrix theory setting of $\alpha = -1$, functional limits of normalized characteristic polynomials were proved
	by Bordenave, Chafa\"i and Garc\'ia-Zelada
	\cite{BordenaveChafaiGarciaZelada2022} for matrices with i.i.d.\
	entries, giving spectral-radius convergence from the nonvanishing
	limit, and by Fran\c{c}ois and Garc\'ia-Zelada
	\cite{FrancoisGarciaZelada2026}, who obtain an exponential Gaussian
	analytic limit for the elliptic Ginibre ensemble, including its GUE
	endpoint, via Hermite polynomials and Chebyshev traces. 
	Recent work of 
	Mohanty and Rajaraman
	\cite{MohantyRajaraman2025} also developed general techniques to obtain spectral norm bounds by analyzing the second moment of the characteristic polynomial.
	
	For invariant Hermitian ensembles, Johansson \cite{Johansson1998}
	established central limit theorems for linear eigenvalue statistics,
	with independent Chebyshev modes in the Gaussian case.
	Central limit theorems for linear eigenvalue statistics of Wigner matrices with four moments
	were obtained by Bai and Yao \cite{BaiYao2005}.
	At $\alpha=-1$, Corollary~\ref{cor:chebyshev-diagonalization} recovers
	the Chebyshev modes of Wigner linear statistics
	\cite{SchenkerSchulzBaldes2007,Kargin2015}; Schenker and
	Schulz-Baldes also discuss the Eulerian weights for complex $\theta$
	\cite[equation~(1.4) and Section~6]{SchenkerSchulzBaldes2007}.
	
	Our method builds heavily on the connection between the expected characteristic polynomial and graph matchings.
	The Godsil--Gutman identity \cite{GodsilGutman1981} identifies the
	expected characteristic polynomial of a random signing of a graph with
	its matching polynomial, which is a Hermite polynomial for the complete
	graph. Forrester and Gamburd \cite[Proposition~11]{ForresterGamburd2006}
	give the Hermite identity for Wigner matrices, and Marcus, Spielman and
	Srivastava \cite{MarcusSpielmanSrivastava2015} built on these ideas in their
	works on interlacing polynomials.
	
	Ravichandran and Leake~\cite{RavichandranLeake2020} proved
	that the $\alpha$-characteristic polynomial is real-rooted for every
	Hermitian matrix when $\alpha=-1/r$, with $r$ a positive integer. Their study was motivated by connections to interlacing polynomials and matrix paving.
	While exact real-rootedness fails for other negative values $\alpha \in (-1,0)$ (Remark 3.10 of \cite{RavichandranLeake2020}), our results show that in the random matrix setting the roots asymptotically concentrate on the real axis.

	The method in this paper is motivated by cluster expansions of polymer gas models (see Appendix~\ref{app:cluster-expansion}). In a related previous work, the authors
	\cite{KoehlerLeung2026} studied zero-free regions for polynomials of the form
	$\per(J+zW)$, with $J$ the all-ones matrix and $W$ random, and used cluster expansion ideas to improve on previous zero-free results of Eldar and Mehraban \cite{eldar2018approximating}. Another relevant work, by McCullagh \cite{McCullagh2012}, used related ideas to derive an approximation for 
	permanents of doubly stochastic matrices by a
	determinant; there too, the leading cycle contributions are encoded by
	a determinant. Zero-free regions also give approximation algorithms
	through Barvinok's interpolation method \cite{Barvinok2017}; here,
	Theorem~\ref{zerofree} would give a quasi-polynomial algorithm, whereas
	the determinantal approximation (Theorem~\ref{thm:determinantal-approximation})
	costs $O(n^3)$.

	\subsection{Organization.}
	Readers interested mainly in the results can read
	Sections~\ref{sec:preliminaries}--\ref{sec:consequences},
	\ref{sec:further-consequences}, and~\ref{sec:determinantal-approximation}
	first. Section~\ref{sec:strategy}
	explains the proof developed in Sections~\ref{sec:expansion}--\ref{sec:fixed-order}.
	Appendix~\ref{app:almost-sure-semicircle} proves the uniform almost-sure
	zero law of Section~\ref{sec:consequences}.
	Appendix~\ref{app:analytic-statistics} gives closed formulas for the
	covariance of the limiting field and for the correction factor of the
	determinantal approximation.
	Appendix~\ref{app:exterior-regime} collects the geometry and
	consequences of the exterior limit for $|\alpha|>1$.
	Appendix~\ref{app:cluster-expansion} is not needed for the proofs.
	It gives cluster-expansion and trace-method context for the argument
	and explains why the
	raw cycle expansion must first be resummed over matchings. This
	resummation is an exact form of the tree summation in the
	trace-method proof of the semicircle law, so readers familiar with
	moment methods or with cluster expansions may find the appendix useful.
	
	\subsection{Notation.}\label{sec:notation}
	We write $\widehat\C:=\C\cup\{\infty\}$ for the Riemann sphere.
	For an open set $U\subset\C^d$, write
	$\mathcal H(U)$ for the space of holomorphic functions on $U$,
	equipped with the topology of uniform convergence on compact subsets;
	it is a Polish space.
	A \emph{random holomorphic function} on $U$ is a random element of
	$\mathcal H(U)$.
	For a compact set $K\Subset U$, put 
	\[ \|f\|_K:=\sup_{z\in K}|f(z)|. \]
	The symbol $\Rightarrow$ denotes weak convergence, including convergence
	in distribution for random elements. For a map $f$ and a measure $\nu$,
	$f_\#\nu:=\nu\circ f^{-1}$ denotes the pushforward, and
	$(n)_m:=n(n-1)\cdots(n-m+1)$, with $(n)_0:=1$, the falling factorial.
	For random elements $X,Y$ of a separable Banach space $(E,\|\cdot\|)$
	with $\E\|X\|^2,\E\|Y\|^2<\infty$, the \emph{quadratic Wasserstein
		distance} between their laws is
	\[
	W_2(X,Y):=\inf\bigl(\E\|X'-Y'\|^2\bigr)^{1/2},
	\]
	the infimum over all pairs $(X',Y')$ on a common probability space with
	$X'\sim X$ and $Y'\sim Y$. It is a metric on laws, and
	$W_2(X,Y)\leq(\E\|X-Y\|^2)^{1/2}$ when $X$ and $Y$ are defined jointly.
	Moreover, $W_2(X_n,X)\to0$ if and only if $X_n\Rightarrow X$ and
	$\E\|X_n\|^2\to\E\|X\|^2$ (equivalently, $X_n\Rightarrow X$ and the
	$\|X_n\|^2$ are uniformly integrable) \cite[Definition~6.8 and Theorem~6.9]{Villani2009}.
	For random continuous functions $X,Y$ on a set containing a compact set
	$K$, we write $W_{2,K}(X,Y):=W_2(X|_K,Y|_K)$ in $C(K)$ with the supremum
	norm.
	We use the notation 
	\[ g_\nu(z):=\int\frac{1}{z-x}\,d\nu(x) \]
	for the Cauchy transform of a compactly supported probability measure $\nu$,
	defined outside its support. For measures supported on $\R$, this is
	also called the \emph{Stieltjes transform}, sometimes with a minus sign.
	Compactly supported probability measures on $\R$ are uniquely
	determined by their Stieltjes transforms \cite[Section~2.4.1]{AndersonGuionnetZeitouni2010}.
	Let $\mu_\scal(dt)=(2\pi)^{-1}\sqrt{4-t^2}\,\mathbf 1_{[-2,2]}(t)\,dt$ be
	the semicircle law. Its Cauchy transform is
	\begin{equation}\label{eq:semicircle-transform}
		g_{\mu_\scal}(z)=\frac{z-\sqrt{z^2-4}}{2},
		\qquad z\in\C\setminus[-2,2],
	\end{equation}
	with the square root analytic off $[-2,2]$ and asymptotic to $z$ at
	infinity; it satisfies $g_{\mu_\scal}^2-zg_{\mu_\scal}+1=0$ and
	$g_{\mu_\scal}(z)\sim z^{-1}$ at infinity.
	
	\section{The expected polynomial}\label{sec:preliminaries}
	The expected $\alpha$-characteristic polynomial is a matching
	polynomial (Section~\ref{sec:expected-polynomial}), and the Cauchy
	transform of its zeros obeys a Riccati equation
	(Section~\ref{sec:hermite-riccati}). The limiting equation is solved by
	the Cauchy transform $g_\alpha$ of a rotated semicircle law, which is
	the limiting zero law (Section~\ref{sec:rotated-semicircle}).
	
	\subsection{Cycle decomposition and the expected polynomial.}\label{sec:expected-polynomial}
	Write $P_{n,\alpha}(z):=\E p_{n,\alpha}(z)$ for the expected
	$\alpha$-characteristic polynomial.
	We first decompose $\det_\alpha(zI-M)$ by cycle type; taking
	expectations shows that $P_{n,\alpha}$ is a matching polynomial, and
	hence a Hermite polynomial.
	
	For a Hermitian $n\times n$ matrix $M$ and an oriented cycle
	$\vec c=(i_1\,\cdots\,i_\ell)$, write
	$M_{\vec c}:=M_{i_1i_2}\cdots M_{i_\ell i_1}$.
	For $\sigma\in\mathfrak S_n$, let $\operatorname{Fix}(\sigma)$ be its
	fixed points, $\mathfrak m(\sigma)$ its transpositions, viewed as
	unordered edges, and $C_{\geq3}(\sigma)$ its cycles of length at least
	three.
	
	\begin{lemma}[Cycle decomposition]\label{lem:cycle-decomposition}
		For every Hermitian $n\times n$ matrix $M$ and every $\alpha,z\in\C$,
		\begin{equation}\label{eq:proof-cycle-start}
			\det_\alpha(zI-M)
			=\sum_{\sigma\in\mathfrak S_n}
			\prod_{i\in\operatorname{Fix}(\sigma)}(z-M_{ii})
			\prod_{\{i,j\}\in\mathfrak m(\sigma)}\alpha|M_{ij}|^2
			\prod_{\vec c\in C_{\geq3}(\sigma)}
			(-1)^{|\vec c|}\alpha^{|\vec c|-1}M_{\vec c}.
		\end{equation}
	\end{lemma}
	\begin{proof}
		Since $\alpha^{n-c(\sigma)}=\prod_c\alpha^{|c|-1}$ over the cycles $c$
		of $\sigma$, and the off-diagonal entries of $zI-M$ are $-M_{ij}$, a
		fixed point contributes $z-M_{ii}$, a transposition contributes
		$\alpha M_{ij}M_{ji}=\alpha|M_{ij}|^2$, and a cycle of length
		$\ell\geq3$ contributes $(-1)^\ell\alpha^{\ell-1}M_{\vec c}$.
		Multiplying these factors over the disjoint cycles of each
		permutation gives the identity.
	\end{proof}
	
	For $r\geq0$ and $x\in\C$, the \emph{matching polynomial} of the
	complete graph $K_r$ with edge weight $x$ is
	$\sum_{\mathfrak m}x^{|\mathfrak m|}z^{r-2|\mathfrak m|}$, where
	$\mathfrak m$ ranges over the matchings of $K_r$, that is, the sets of
	pairwise disjoint edges \cite{HeilmannLieb1972,GodsilGutman1981}; for
	$x=-1$ this is the classical matching polynomial. A $k$-edge matching is
	a choice of $2k$ vertices together with one of the $(2k-1)!!$ ways to
	pair them, so collecting terms by $k$ gives
	\begin{equation}\label{eq:hermite-definition}
		\sum_{\mathfrak m}x^{|\mathfrak m|}z^{r-2|\mathfrak m|}
		=\sum_{k=0}^{\lfloor r/2\rfloor}
		\frac{r!}{2^k k!(r-2k)!}x^kz^{r-2k}
		=:\Her_r(z;-x).
	\end{equation}
	That is, the \emph{monic Hermite polynomial with variance parameter
		$s$}, $\Her_r(z;s)$, is the matching polynomial of $K_r$ with edge
	weight $-s$; for $s=1$ these are the probabilists' Hermite
	polynomials.
	
	\begin{lemma}[Expected polynomial]\label{lem:annealed-polynomial}
		For every $n\geq1$ and $\alpha\in\C$, $P_{n,\alpha}$ is the matching
		polynomial of $K_n$ with edge weight $\alpha/n$:
		\begin{equation}\label{eq:annealed-matching}
			P_{n,\alpha}(z)
			=\sum_{\mathfrak m}\left(\frac{\alpha}{n}\right)^{|\mathfrak m|}z^{n-2|\mathfrak m|}
			=\Her_n\Bigl(z;-\frac\alpha n\Bigr).
		\end{equation}
	\end{lemma}
	\begin{proof}
		Take expectations in \eqref{eq:proof-cycle-start} with $M=M_n$.
		Since the cycles of $\sigma$ are vertex-disjoint, distinct cycles use
		disjoint sets of independent matrix entries, so the expectation
		factors over cycles. Every cycle of length at least three has mean
		zero by independence and centering, since each unordered matrix edge
		appears once. Fixed points and transpositions average to $z$ and
		$\alpha/n$, respectively, giving the matching sum, which is
		$\Her_n(z;-\alpha/n)$ by definition.
	\end{proof}
	For the determinant ($\alpha=-1$), this Hermite representation is
	classical~\cite[Proposition~11]{ForresterGamburd2006}.
	
	\subsection{Hermite identities and the Riccati equation.}\label{sec:hermite-riccati}
	The following classical identities give the structure of $P_{n,\alpha}$
	used below; matchings make them transparent.
	
	\begin{lemma}[Hermite identities]\label{lem:hermite-identities}
		For $s,z\in\C$, with $\Her_r(z;0)=z^r$:
		\begin{enumerate}[label=(\alph*)]
			\item for $c\in\C$ and $r\geq0$,
			\begin{equation}\label{eq:hermite-scaling}
				\Her_r(cz;c^2s)=c^r\Her_r(z;s);
			\end{equation}
			\item for $r\geq1$,
			\begin{equation}\label{eq:hermite-derivative}
				\partial_z\Her_r(z;s)=r\Her_{r-1}(z;s);
			\end{equation}
			\item for $r\geq2$,
			\begin{equation}\label{eq:hermite-recurrence}
				\Her_r(z;s)=z\,\Her_{r-1}(z;s)-(r-1)s\,\Her_{r-2}(z;s).
			\end{equation}
		\end{enumerate}
	\end{lemma}
	
	\begin{proof}
		(a) Scaling $z$ by $c$ and the edge weight by $c^2$ multiplies each
		matching term by $c^r$. (b) Differentiating
		$\sum_{\mathfrak m}(-s)^{|\mathfrak m|}z^{r-2|\mathfrak m|}$ brings down
		the number of unmatched vertices, that is, it sums over the choice of
		one unmatched vertex $v$. Removing $v$ leaves a matching of
		$K_r-v\cong K_{r-1}$, so $\partial_z\Her_r=\sum_v\Her_{r-1}=r\Her_{r-1}$.
		(c) Condition on whether a fixed vertex of $K_r$ is unmatched or
		matched to one of the other $r-1$ vertices.
	\end{proof}
	
	Let $\mu_{n,\alpha}$ be the empirical measure of the zeros of
	$P_{n,\alpha}$, counted with multiplicity, and write
	\[
	h_{n,\alpha}(z):=g_{\mu_{n,\alpha}}(z)
	=\frac1n\frac{P_{n,\alpha}'(z)}{P_{n,\alpha}(z)}
	=\frac{\Her_{n-1}(z;-\alpha/n)}{\Her_n(z;-\alpha/n)}
	\]
	for its Cauchy transform, defined away from the zeros of $P_{n,\alpha}$,
	where the last equality uses \eqref{eq:hermite-derivative}.
	
	\begin{lemma}[Riccati equation]\label{lem:hermite-riccati}
		Away from the zeros of $P_{n,\alpha}$,
		\begin{equation}\label{eq:hermite-riccati}
			\frac\alpha n\,h_{n,\alpha}'
			+\alpha h_{n,\alpha}^2+zh_{n,\alpha}-1=0.
		\end{equation}
	\end{lemma}
	
	\begin{proof}
		The recurrence \eqref{eq:hermite-recurrence} at $r=n$, rewritten
		with \eqref{eq:hermite-derivative}, reads
		\[ P_{n,\alpha}=\frac zn P_{n,\alpha}'+\frac{\alpha}{n^2}P_{n,\alpha}''. \]
		Differentiating $h_{n,\alpha}=P_{n,\alpha}'/(nP_{n,\alpha})$ gives
		$h_{n,\alpha}'=P_{n,\alpha}''/(nP_{n,\alpha})-nh_{n,\alpha}^2$, that is,
		$P_{n,\alpha}''/(n^2P_{n,\alpha})=h_{n,\alpha}'/n+h_{n,\alpha}^2$.
		Dividing the recurrence by $P_{n,\alpha}$ and substituting gives
		$1=zh_{n,\alpha}+\alpha\bigl(h_{n,\alpha}'/n+h_{n,\alpha}^2\bigr)$.
	\end{proof}
	
	Dropping the term $\frac\alpha nh_{n,\alpha}'$ gives the quadratic
	equation \eqref{eq:g-quadratic} below, which identifies the zero law
	(Lemma~\ref{lem:hermite-zero-law}).
	
	\subsection{The rotated semicircle law.}\label{sec:rotated-semicircle}
	For $\alpha\in\C$, set
	\[
	\beta_\alpha:=\sqrt{-\alpha},\qquad
	\Sigma_\alpha:=\beta_\alpha[-2,2],\qquad
	\Omega_\alpha:=\C\setminus\Sigma_\alpha,\qquad
	\mu_\alpha:=(\beta_\alpha\id)_\#\mu_\scal,
	\]
	with the principal square root in $\beta_\alpha$. By symmetry of
	$[-2,2]$ and $\mu_\scal$, $\Sigma_\alpha$ and $\mu_\alpha$ do not depend
	on this choice; for real $\alpha$, $\Sigma_\alpha$ is a real interval when
	$\alpha<0$, an imaginary interval when $\alpha>0$, and $\{0\}$ at
	$\alpha=0$, where $\mu_0=\delta_0$. Put
	\[
	\mathcal U:=\{(\alpha,z)\in\C^2:z\notin\Sigma_\alpha\},
	\qquad
	g_\alpha:=g_{\mu_\alpha}\quad\text{on }\Omega_\alpha.
	\]
	
	\begin{lemma}[Cauchy transform of the rotated semicircle law]\label{lem:semicircle-equation}
		For every $\alpha\in\C$, $g_\alpha$ is holomorphic on $\Omega_\alpha$,
		$g_\alpha(z)\sim z^{-1}$ at infinity, and
		\begin{equation}\label{eq:g-quadratic}
			\alpha g_\alpha(z)^2+zg_\alpha(z)-1=0,\qquad z\in\Omega_\alpha.
		\end{equation}
		Moreover, $\mathcal U$ is open and $g_\alpha(z)$ is jointly
		holomorphic in $(\alpha,z)\in\mathcal U$.
	\end{lemma}
	
	\begin{proof}
		At $\alpha=0$, $g_0(z)=1/z$. For $\alpha\neq0$, a change of variables
		gives $g_\alpha(z)=\beta_\alpha^{-1}g_{\mu_\scal}(z/\beta_\alpha)$, and
		the properties of $g_{\mu_\scal}$ after \eqref{eq:semicircle-transform},
		together with $\alpha/\beta_\alpha^2=-1$, give the first three claims.
		For the last one, average over $x$ and $-x$, which is allowed because
		$\mu_\scal$ is symmetric:
		\[
		g_\alpha(z)=\int\frac{z\,d\mu_\scal(x)}{z^2+\alpha x^2}.
		\]
		For $x\in[-2,2]$, $z^2+\alpha x^2=0$ exactly when $z=\pm\beta_\alpha x$,
		so $\mathcal U$ is the set where $\min_{|x|\leq2}|z^2+\alpha x^2|>0$.
		This minimum is continuous in $(\alpha,z)$, so $\mathcal U$ is open, and
		the integrand is jointly holomorphic and locally bounded there, so
		$g_\alpha(z)$ is jointly holomorphic on $\mathcal U$.
	\end{proof}
	
	The zeros of the expected polynomial follow $\mu_\alpha$:
	
	\begin{lemma}[Zeros of the expected polynomial]\label{lem:hermite-zero-law}
		For $n\geq1$, let $t_{n,1},\ldots,t_{n,n}$ be the zeros of
		$\Her_n(\cdot\,;1/n)$. They are real and simple and lie in
		$(-2,2)$. For every $\alpha\in\C$, the zeros of $P_{n,\alpha}$ are
		$\beta_\alpha t_{n,j}$, counted with multiplicity, and hence lie on
		$\Sigma_\alpha$. For each fixed $\alpha$,
		\begin{equation}\label{eq:hermite-zero-law}
			\mu_{n,\alpha}
			=\frac1n\sum_{j=1}^n\delta_{\beta_\alpha t_{n,j}}
			\Rightarrow\mu_\alpha.
		\end{equation}
		Moreover, $|h_{n,\alpha}(z)|\leq\dist(z,\Sigma_\alpha)^{-1}$, and
		\begin{equation}\label{eq:hermite-transform-limit}
			h_{n,\alpha}(z)\longrightarrow g_\alpha(z)
			\qquad(n\to\infty)
		\end{equation}
		locally uniformly on $\mathcal U$.
	\end{lemma}
	
	The semicircle limit for the normalized Hermite zeros is classical;
	see, e.g., \cite[Theorem~1]{KornyikMichaletzky2016}. We include a proof
	that also gives the finite-$n$ zero bounds and joint transform convergence.
	
	\begin{proof}
		The recurrence \eqref{eq:hermite-recurrence} identifies
		$\Her_n(\cdot\,;1/n)$ as the characteristic polynomial of the symmetric
		tridiagonal matrix with zero diagonal and off-diagonal entries
		$\sqrt{j/n}$, $1\leq j<n$. Its eigenvalues are real and simple, since
		its off-diagonal entries are nonzero, and by Gershgorin's theorem they
		lie in $(-2,2)$, since its absolute row sums are less than $2$. For $\alpha\neq0$, \eqref{eq:annealed-matching} and
		\eqref{eq:hermite-scaling} give
		$P_{n,\alpha}(\beta_\alpha x)=\beta_\alpha^n\Her_n(x;1/n)$, proving
		the asserted scaling of the zeros; at $\alpha=0$, $P_{n,0}=z^n$.
		Their support on $\Sigma_\alpha$ gives
		$|h_{n,\alpha}(z)|\leq\dist(z,\Sigma_\alpha)^{-1}$.
		
		The segments $\Sigma_\alpha$ depend continuously on $\alpha$ in
		Hausdorff distance. Thus $\dist(z,\Sigma_\alpha)$ is continuous and
		positive on $\mathcal U$, and the jointly holomorphic functions
		$h_{n,\alpha}(z)$ form a locally bounded normal family there. Any
		locally uniform subsequential limit $h$ satisfies
		$\alpha h^2+zh-1=0$ by \eqref{eq:hermite-riccati} and convergence of
		derivatives. For fixed $\alpha$, the zeros lie in
		$\{|\zeta|\leq2\sqrt{|\alpha|}\}$, so
		$h_{n,\alpha}(z)=z^{-1}+O_\alpha(|z|^{-2})$ uniformly in $n$ near
		infinity. Hence $h(z)\sim z^{-1}$, as does $g_\alpha$. Subtracting their
		quadratic equations gives
		\[
		(h-g_\alpha)\bigl[z+\alpha(h+g_\alpha)\bigr]=0.
		\]
		The bracket is nonzero near infinity, so $h=g_\alpha$ there and then
		throughout the connected domain $\Omega_\alpha$ by the identity theorem.
		Every subsequential limit is therefore $g_\alpha$, proving
		\eqref{eq:hermite-transform-limit} locally uniformly on $\mathcal U$.
		
		For fixed $\alpha$, the measures $\mu_{n,\alpha}$ are supported on the
		compact segment $\Sigma_\alpha$. Any weak subsequential limit has Cauchy
		transform $g_\alpha=g_{\mu_\alpha}$ off this segment. A measure on a
		segment is determined by its Cauchy transform, as one sees by rotating
		to the real line, so the limit is $\mu_\alpha$. This proves
		\eqref{eq:hermite-zero-law}.
	\end{proof}
	
	\section{The main result}\label{sec:functional-limits}
	The limit of $p_{n,\alpha}/\E p_{n,\alpha}$ is naturally expressed in the
	coordinate $\tau_\alpha=-\alpha g_\alpha$, where $g_\alpha$ is the Cauchy
	transform of $\mu_\alpha$ from Section~\ref{sec:rotated-semicircle} and
	$\Sigma_\alpha=\operatorname{supp}\mu_\alpha$. For $\alpha\neq0$,
	$\tau_\alpha$ maps $\Omega_\alpha\cup\{\infty\}$ conformally onto a disk
	(Section~\ref{sec:disk-coordinate}), and the limit is a power series in
	$\tau_\alpha$, convergent where $|\tau_\alpha|<1$.
	Section~\ref{sec:functional-limit} states the functional limit theorem.
	
	\subsection{The disk coordinate.}\label{sec:disk-coordinate}
	Set
	\begin{equation}\label{eq:tau-definition}
		\tau_\alpha(z):=-\alpha g_\alpha(z),\qquad z\in\Omega_\alpha.
	\end{equation}
	Multiplying \eqref{eq:g-quadratic} by $\alpha$ shows that $\tau_\alpha$
	is a root of the quadratic
	\begin{equation}\label{eq:basic-quadratic}
		X^2-zX-\alpha,
	\end{equation}
	whose roots have sum $z$ and product $-\alpha$; the other root is
	$z-\tau_\alpha=1/g_\alpha$. Thus
	\begin{equation}\label{eq:two-roots}
		X^2-zX-\alpha
		=\Bigl(X-\frac1{g_\alpha(z)}\Bigr)\bigl(X-\tau_\alpha(z)\bigr).
	\end{equation}
	The next two lemmas use this quadratic to show that $\tau_\alpha$ is the
	root of smaller modulus and a conformal coordinate.
	
	\begin{lemma}[Root separation]\label{lem:quadratic-branches}
		For each $\alpha\in\C$, the two roots of \eqref{eq:basic-quadratic}
		have equal modulus if and only if $z\in\Sigma_\alpha$. For
		$z\in\Omega_\alpha$, $\tau_\alpha(z)$ is the root of smaller modulus,
		so $|\tau_\alpha(z)|<\sqrt{|\alpha|}$ when $\alpha\neq0$.
	\end{lemma}
	
	\begin{proof}
		At $\alpha=0$ the roots are $0$ and $z$, with equal modulus only at
		$z=0$. For $\alpha\neq0$, equal modulus means that both roots have
		modulus $\sqrt{|\alpha|}$. Writing one as $\beta_\alpha\omega$ with
		$|\omega|=1$, the other is
		$\beta_\alpha/\omega=\beta_\alpha\overline\omega$,
		so $z=\beta_\alpha(\omega+\overline\omega)\in\beta_\alpha[-2,2]$.
		Conversely, for $z=2\beta_\alpha\cos\varphi$ with $\varphi\in\R$, the
		roots are $\beta_\alpha e^{\pm i\varphi}$.
		
		On $\Omega_\alpha$ the moduli differ. Near infinity,
		$1/g_\alpha(z)\sim z$ is the larger root by
		Lemma~\ref{lem:semicircle-equation}; since $\Omega_\alpha$ is connected and
		the moduli never coincide there, $1/g_\alpha$ is the larger root and
		$\tau_\alpha$ the smaller one throughout. For $\alpha\neq0$, the
		product of the moduli is $|\alpha|$, so
		$|\tau_\alpha|<\sqrt{|\alpha|}$.
	\end{proof}

	\begin{lemma}[The disk coordinate]\label{lem:conformal-disk-coordinate}
		For $\alpha\neq0$, $\tau_\alpha$, extended by
		$\tau_\alpha(\infty)=0$, is a conformal bijection from
		$\Omega_\alpha\cup\{\infty\}$ onto $\{t:|t|<\sqrt{|\alpha|}\}$, with
		inverse
		\begin{equation}\label{eq:joukowski-inverse}
			J_\alpha(t):=t-\frac\alpha t,\qquad 0<|t|<\sqrt{|\alpha|},
		\end{equation}
		and $J_\alpha(0):=\infty$. Equivalently, for every $\alpha\in\C$
		(including $\alpha=0$), $g_\alpha$, extended by $g_\alpha(\infty)=0$,
		maps $\Omega_\alpha\cup\{\infty\}$ conformally onto
		$\{u\in\C:|\alpha u^2|<1\}$, with inverse
		\begin{equation}\label{eq:cauchy-inverse}
			z_\alpha(u):=u^{-1}-\alpha u,\qquad |\alpha u^2|<1,\ u\neq0,
		\end{equation}
		and $z_\alpha(0):=\infty$; for $\alpha\neq0$, $z_\alpha(u)=J_\alpha(-\alpha u)$.
	\end{lemma}
	
	\begin{proof}
		Solving $t^2-zt-\alpha=0$ ($t\neq0$) for $z$ gives $z=J_\alpha(t)$.
		For $z\in\Omega_\alpha$, the product of the roots is $-\alpha\neq0$, so
		$0<|\tau_\alpha(z)|<\sqrt{|\alpha|}$ by
		Lemma~\ref{lem:quadratic-branches} and the other root is
		$-\alpha/\tau_\alpha(z)$; the sum of the roots then gives
		$J_\alpha(\tau_\alpha(z))=z$.
		Conversely, if $0<|t|<\sqrt{|\alpha|}$ and $z:=J_\alpha(t)$, then the
		roots of \eqref{eq:basic-quadratic} are $t$ and $-\alpha/t$, with
		$|t|<\sqrt{|\alpha|}<|\alpha/t|$, so Lemma~\ref{lem:quadratic-branches}
		gives $z\in\Omega_\alpha$ and $\tau_\alpha(z)=t$. Thus $\tau_\alpha$
		and $J_\alpha$ are mutually inverse between $\Omega_\alpha$ and
		$\{0<|t|<\sqrt{|\alpha|}\}$. Since $g_\alpha(z)\sim z^{-1}$
		(Lemma~\ref{lem:semicircle-equation}), $\tau_\alpha(z)=-\alpha g_\alpha(z)\to0$
		at infinity. Riemann's removable singularity
		theorem extends $\tau_\alpha$ holomorphically by
		$\tau_\alpha(\infty)=0$, so this holomorphic bijection is conformal.
		For $\alpha\neq0$, the second statement follows by setting $t=-\alpha u$,
		since $\tau_\alpha=-\alpha g_\alpha$ and $|t|<\sqrt{|\alpha|}$ is
		equivalent to $|\alpha u^2|<1$; at $\alpha=0$, $g_0(z)=1/z$ has inverse
		$z_0(u)=1/u$ on $\C$.
	\end{proof}
	
	At $\alpha=-1$, Lemma~\ref{lem:conformal-disk-coordinate} recovers the classical picture:
	$\Sigma_{-1}=[-2,2]$, $\tau_{-1}=g_{-1}=g_{\mu_\scal}$, and
	$J_{-1}(u)=u+1/u$ is the Joukowski map \cite{Ransford1995}.
	
	\subsection{The functional limit theorem.}\label{sec:functional-limit}
	Normalize the $\alpha$-characteristic polynomial by its expectation:
	\begin{equation}\label{eq:normalized-polynomial-definition}
		F_{n,\alpha}(z):=p_{n,\alpha}(z)/P_{n,\alpha}(z).
	\end{equation}
	Since all zeros of $P_{n,\alpha}$ lie on $\Sigma_\alpha$, the
	polynomials $p_{n,\alpha}(z)$ and $P_{n,\alpha}(z)$ in $(\alpha,z)$
	define a jointly holomorphic ratio $F_n(\alpha,z):=F_{n,\alpha}(z)$ on
	$\mathcal U$, with $F_{n,\alpha}\in\mathcal H(\Omega_\alpha)$ for every
	$\alpha\in\C$.
	
	The limit has one Gaussian mode for each cycle length $\ell$: diagonal
	entries for $\ell=1$, squared off-diagonal magnitudes for $\ell=2$, and
	simple cycles of length $\ell\geq3$ (compare \eqref{eq:informal-main}).
	Their variances are
	\begin{equation}\label{eq:mode-variances}
		\sigma_1^2:=v,\qquad \sigma_2^2:=\frac\kappa2,\qquad
		\sigma_\ell^2:=\frac{1+\rho_\ell(\theta)}{\ell},\quad \ell\geq3,
	\end{equation}
	where $\kappa:=\Var(|\xi_{12}|^2)$. For $\ell\geq2$, $\rho_\ell(\theta)$
	is the average phase of an $\ell$-cycle. By Hermitian symmetry, an edge
	traversed upward ($i<j$) carries $\theta$ and one traversed downward
	carries $\overline\theta$. So, with $\operatorname{asc}(\pi)$ the number
	of cyclic ascents of $\pi\in\mathfrak S_\ell$, that is, of indices $j$
	with $\pi(j)<\pi(j+1)$ where $\pi(\ell+1):=\pi(1)$,
	\begin{equation}\label{eq:cycle-order-average}
		\rho_\ell(\theta):=\frac1{\ell!}\sum_{\pi\in\mathfrak S_\ell}
		\operatorname{Re}\!\left(\theta^{\operatorname{asc}(\pi)}
		\overline\theta^{\,\ell-\operatorname{asc}(\pi)}\right),
		\qquad \ell\geq2.
	\end{equation}
	They satisfy $|\rho_\ell(\theta)|\leq|\theta|^\ell$ and
	$\rho_\ell(\theta)=\theta^\ell$ for real $\theta$;
	Appendix~\ref{sec:eulerian-cycle-average} computes them via Eulerian
	numbers.
	
	Fix independent centered real Gaussian variables $(G_\ell)_{\ell\geq1}$
	with $\Var(G_\ell)=\sigma_\ell^2$, and form one Gaussian analytic
	function on the unit disk:
	\begin{equation}\label{eq:common-gaussian-field}
		S(t):=\sum_{\ell\geq1}G_\ell t^\ell,\qquad
		\mathcal C(t):=\E[S(t)^2]=\sum_{\ell\geq1}\sigma_\ell^2t^{2\ell},
		\qquad |t|<1.
	\end{equation}
	These series converge absolutely and locally uniformly on a single event
	of probability one.\footnote{For every $0<r<1$,
		$\E\sum_{\ell\geq1}r^\ell|G_\ell|<\infty$, since
		$\E|G_\ell|\leq\sigma_\ell$ and $\sigma_\ell^2\leq2/\ell$ for
		$\ell\geq3$. Taking rational $r\in(0,1)$ gives a single event of
		probability one on which $S$ converges absolutely and locally
		uniformly on the unit disk. The same variance bound gives this
		convergence for $\mathcal C$.}
	
	Since $S$ is defined on the unit disk, $S(\tau_\alpha(z))$ is defined
	precisely where $|\tau_\alpha(z)|<1$. Set
	\begin{equation}\label{eq:natural-joint-domain}
		\mathcal D:=\{(\alpha,z)\in\mathcal U:
		|\tau_\alpha(z)|<1\}.
	\end{equation}
	By Lemma~\ref{lem:semicircle-equation}, $\mathcal U$ is open and
	$g_\alpha(z)$ and $\tau_\alpha(z)$ are jointly holomorphic on it, so
	$\mathcal D$ is open. By Lemma~\ref{lem:quadratic-branches},
	$|\tau_\alpha|<1$ on $\Omega_\alpha$ when $|\alpha|\leq1$, so
	$\mathcal D$ contains $\{\alpha\}\times\Omega_\alpha$ for every
	$|\alpha|\leq1$ (for $|\alpha|>1$, the condition $|\tau_\alpha|<1$
	restricts $z$ to the exterior of an ellipse; see
	Appendix~\ref{app:exterior-regime}). For $(\alpha,z)\in\mathcal D$ with
	$\alpha\neq0$, define
	\begin{equation}\label{eq:limiting-exponent}
		H_\alpha(z):=\frac{S(\tau_\alpha(z))}{\alpha}
		-\frac{\mathcal C(\tau_\alpha(z))}{2\alpha^2},
		\qquad \Psi_\alpha(z):=e^{H_\alpha(z)}.
	\end{equation}
	Since $\tau_\alpha=-\alpha g_\alpha$, we have
	$\tau_\alpha^\ell/\alpha=(-1)^\ell\alpha^{\ell-1}g_\alpha^\ell$, so in
	powers of $g_\alpha$,
	\begin{equation}\label{eq:mode-expansion}
		H_\alpha(z)=\sum_{\ell\geq1}(-1)^\ell\alpha^{\ell-1}G_\ell\,g_\alpha(z)^\ell
		-\frac12\sum_{\ell\geq1}\alpha^{2\ell-2}\sigma_\ell^2\,
		g_\alpha(z)^{2\ell}.
	\end{equation}
	The right side also makes sense at $\alpha=0$, and we use it to define
	$H_0(z)=-G_1/z-v/(2z^2)$ and $\Psi_0:=e^{H_0}$. Since
	$(-1)^\ell\alpha^{\ell-1}g_\alpha^\ell=-g_\alpha\tau_\alpha^{\ell-1}$ and
	$|\tau_\alpha|<1$ on $\mathcal D$, both series converge locally uniformly
	on $\mathcal D$, so $\Psi(\alpha,z):=\Psi_\alpha(z)$ is almost surely
	jointly holomorphic there, and $F_n$ and $\Psi$ both lie in
	$\mathcal H(\mathcal D)$. The $\ell$th mode enters with the cycle weight
	of Lemma~\ref{lem:cycle-decomposition} times one factor $g_\alpha$ for
	each vertex of the cycle, as in \eqref{eq:informal-main};
	Section~\ref{sec:strategy} explains where these factors come from.
	
	\begin{theorem}[Functional limit]\label{gaussian}
		Form all polynomials from the same matrix $M_n$ and use the common
		Gaussian function $S$ in \eqref{eq:common-gaussian-field}. Then
		\begin{equation}\label{eq:joint-parameter-limit}
			F_n\Rightarrow\Psi
			\qquad\text{in }\mathcal H(\mathcal D).
		\end{equation}
		In particular, for every $\alpha\in\C$ with $|\alpha|\leq1$,
		\begin{equation}\label{eq:functional-limit}
			F_{n,\alpha}\Rightarrow\Psi_\alpha
			\qquad\text{in }\mathcal H(\Omega_\alpha).
		\end{equation}
	\end{theorem}

	\begin{remark}[Wick exponential]\label{rem:wick-exponential}
		For $X$ in the complex linear span of a centered real Gaussian family,
		its \emph{Wick exponential} is $:\!\exp(X)\!:\,=\exp(X-\tfrac12\E X^2)$
		\cite[Chapter~3]{Janson1997}, with $\E X^2$ rather than $\E|X|^2$, so
		that $\E:\!\exp(X)\!:\,=1$. Thus
		$\Psi_\alpha(z)=\,:\!\exp\bigl(S(\tau_\alpha(z))/\alpha\bigr)\!:$ for
		$\alpha\neq0$, and $\Psi_0(z)=\,:\!\exp(-G_1/z)\!:$. The
		\emph{Wick correction} $-\mathcal C(\tau_\alpha(z))/(2\alpha^2)$ makes
		$\E\Psi_\alpha(z)=1$, matching $\E F_{n,\alpha}(z)=1$.
	\end{remark}
	
	The convergence also holds in the quadratic Wasserstein distance
	$W_{2,K}$ (Section~\ref{sec:notation}), which in addition controls
	second moments.
	
	\begin{proposition}[Quadratic Wasserstein convergence]\label{prop:wasserstein-convergence}
		For every compact $\mathcal K\Subset\mathcal D$,
		\begin{equation}\label{eq:compact-wasserstein-convergence}
			W_{2,\mathcal K}(F_n,\Psi)\longrightarrow0.
		\end{equation}
		In particular, the squared norms $\|F_n\|_{\mathcal K}^2$ are
		uniformly integrable and
		$\E\|F_n\|_{\mathcal K}^2\to\E\|\Psi\|_{\mathcal K}^2<\infty$.
	\end{proposition}
	
	Theorem~\ref{gaussian} and
	Proposition~\ref{prop:wasserstein-convergence} are proved in
	Section~\ref{sec:assembly}, and Section~\ref{sec:extension-at-infinity}
	extends the convergence to $z=\infty$.
	
\begin{figure}[htbp]
		\centering
		\includegraphics[width=0.9\textwidth]{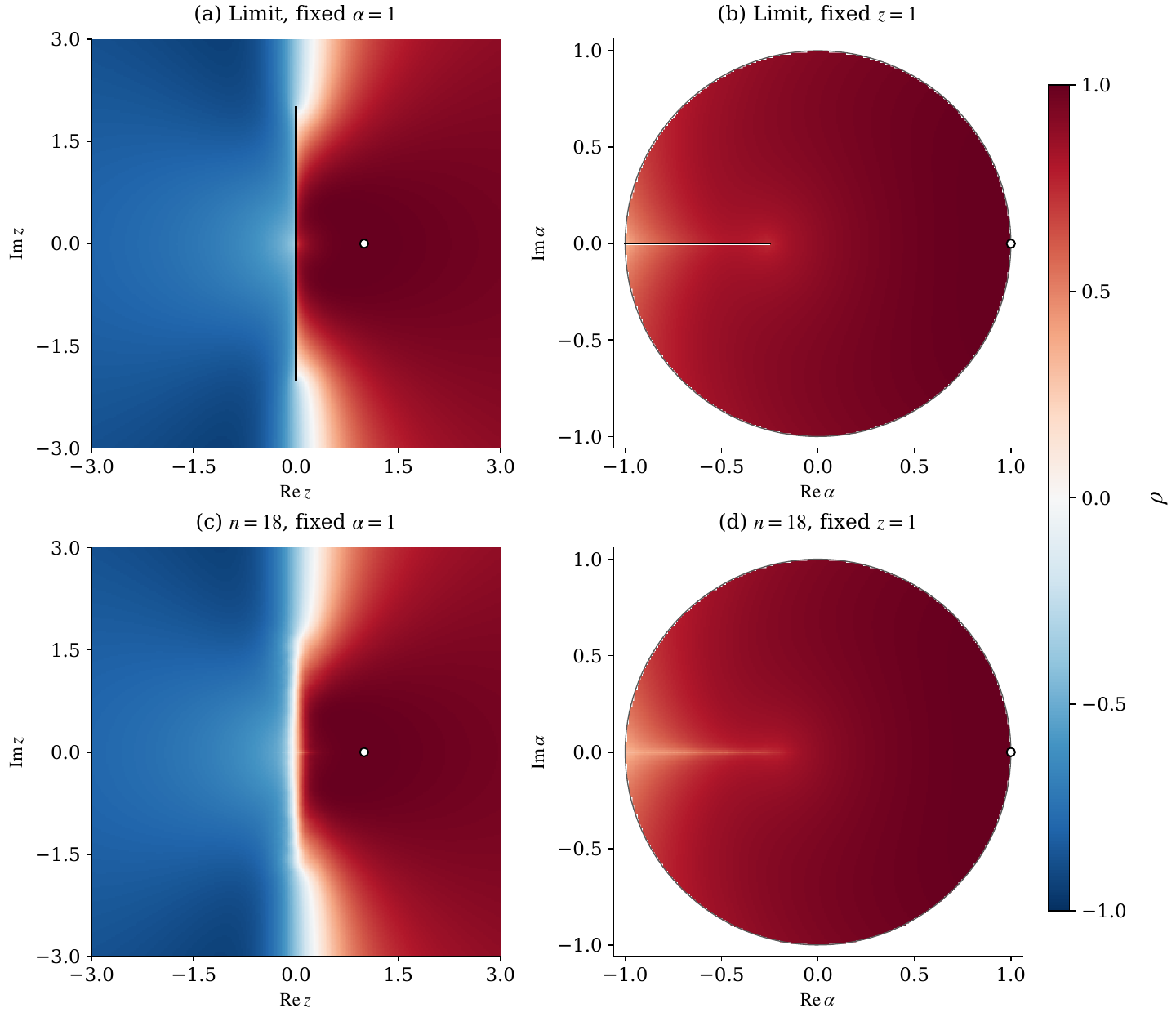}
		\caption{GUE correlations of $\log|\Psi_\alpha(z)|$ with
			$\log|\Psi_1(1)|$ in (a,b), compared with empirical correlations
			of $\log|p_{18,\alpha}(z)|$ with $\log|p_{18,1}(1)|$ in (c,d).
			The empirical panels use the same 4,096 independent GUE matrices.
			The left column fixes $\alpha=1$; the right fixes $z=1$ and varies
			$\alpha$ over the unit disk. Open circles mark the reference point.
			All panels use the same color scale, with red and blue indicating
			positive and negative correlations. Black segments mark excluded
			limiting parameters: $z\in i[-2,2]$ in (a) and
			$\alpha\in[-1,-1/4]$ in (b); the empirical panels include these
			segments. The limiting kernels are given in
			Proposition~\ref{prop:covariance-kernel}; GOE has the same limiting
			correlations.}
		\label{fig:limiting-correlation}
	\end{figure}

	Figure~\ref{fig:limiting-correlation} illustrates the dependence between
	values of the joint limit $\Psi$ through correlations of its
	logarithmic magnitude along two parameter slices, compared with
	empirical estimates from 4,096 independent $18\times18$ GUE matrices.
	The empirical panels use $\log|p_{18,\alpha}(z)|$ directly;
	normalization by $P_{18,\alpha}(z)$ leaves correlations and variances
	unchanged wherever it is defined.
	
	\subsection{Extension at infinity.}\label{sec:extension-at-infinity}
	Adjoining the point at infinity in the $z$-variable lets us control
	large zeros and work near $z=\infty$. Set
	\[
	\widehat{\mathcal D}:=\mathcal D\cup(\C\times\{\infty\})
	\subset\C\times\widehat\C,
	\]
	an open subset of $\C\times\widehat\C$; holomorphy at infinity is
	understood in the local coordinate $w=1/z$.
	
	\begin{corollary}[Extension at infinity]\label{lem:extension-at-infinity}
		\begin{enumerate}[label=(\alph*)]
			\item $F_n$ and $\Psi$ extend jointly holomorphically to
			$\widehat{\mathcal D}$, with $F_{n,\alpha}(\infty)=\Psi_\alpha(\infty)=1$.
			\item For every compact $\mathcal K\Subset\widehat{\mathcal D}$, there
			is a compact $\mathcal K_0\Subset\mathcal D$ such that
			$\|f\|_{\mathcal K}\leq\|f\|_{\mathcal K_0}$ for all
			$f\in\mathcal H(\widehat{\mathcal D})$.
			\item The convergence in Theorem~\ref{gaussian} holds in
			$\mathcal H(\widehat{\mathcal D})$, with the topology of uniform
			convergence on compact subsets.
		\end{enumerate}
	\end{corollary}
	
	\begin{proof}
		(a) In the coordinate $w=1/z$, $w^np_{n,\alpha}(1/w)$ and
		$w^nP_{n,\alpha}(1/w)$ are polynomials in $(\alpha,w)$ with constant
		term one, so their ratio $F_{n,\alpha}(1/w)$ extends holomorphically
		to $w=0$ with value $1$. For $\Psi$, $g_\alpha$ is the reciprocal of
		the larger root of $X^2-zX-\alpha$, so
		\[
		g_\alpha(1/w)=\frac{2w}{1+\sqrt{1+4\alpha w^2}}=w+O(w^3),
		\]
		which is jointly holomorphic for $|\alpha w^2|<1/4$. Substituting
		into \eqref{eq:mode-expansion} gives $H_\alpha(\infty)=0$, hence
		$\Psi_\alpha(\infty)=1$.
		
		(b) Let $A$ be the projection of $\mathcal K$ onto the $\alpha$-plane,
		and choose $r_0$ so large that
		$A\times(\{|z|\geq r_0\}\cup\{\infty\})\subset\widehat{\mathcal D}$.
		Then
		$\mathcal K_0:=\bigl(\mathcal K\cap\{|z|\leq r_0\}\bigr)
		\cup\bigl(A\times\{|z|=r_0\}\bigr)$ is compact in $\mathcal D$, and
		the maximum principle on each exterior disk
		$\{|z|\geq r_0\}\cup\{\infty\}$ gives the bound.
		
		(c) By (b), each supremum norm $\|\cdot\|_{\mathcal K}$ with
		$\mathcal K\Subset\widehat{\mathcal D}$ is bounded by a supremum norm
		$\|\cdot\|_{\mathcal K_0}$ with $\mathcal K_0\Subset\mathcal D$. Hence restriction
		identifies $\mathcal H(\widehat{\mathcal D})$ with a closed subspace of
		$\mathcal H(\mathcal D)$, carrying the same topology. By (a), $F_n$ and
		$\Psi$ lie in this subspace, so Theorem~\ref{gaussian} gives
		convergence in law there.
	\end{proof}
	
	\section{Zeros and the semicircle law}\label{sec:consequences}
	Since the limit in Theorem~\ref{gaussian} never
	vanishes, it controls the zeros of $p_{n,\alpha}$. This section derives
	zero confinement and the semicircle law. Further consequences, which
	use the structure of the exponent $H_\alpha$, are the fluctuations of the
	zeros in Section~\ref{sec:further-consequences} and the determinantal
	approximation in Section~\ref{sec:determinantal-approximation}.
	
	\subsection{Uniform zero confinement.}
	\begin{theorem}[Uniform zero confinement]\label{zerofree}
		For every $\varepsilon>0$, with probability tending to one, every
		zero $\zeta$ of $p_{n,\alpha}$ satisfies
		$\dist(\zeta,\Sigma_\alpha)<\varepsilon$, simultaneously for all
		$|\alpha|\leq1$.
	\end{theorem}
	
	\begin{proof}
		Let $B$ be the closed unit disk and write
		$N_\varepsilon(\Sigma_\alpha):=\{z\in\C:
		\dist(z,\Sigma_\alpha)<\varepsilon\}$. Hausdorff continuity of the
		segments $\Sigma_\alpha$ implies that
		\[
		K_\varepsilon:=\{(\alpha,z):\alpha\in B,\
		z\in\widehat\C\setminus N_\varepsilon(\Sigma_\alpha)\}
		\]
		is compact in $\widehat{\mathcal D}$. The limit $\Psi$ is nonzero
		on this set, including at infinity, and hence its minimum modulus
		there is positive almost surely. Corollary~\ref{lem:extension-at-infinity}
		and Portmanteau therefore give
		\[
		\mathbb P\left(\min_{K_\varepsilon}|F_n|=0\right)\longrightarrow0.
		\]
		All zeros of $P_{n,\alpha}$ lie on $\Sigma_\alpha$, so any zero of
		$p_{n,\alpha}$ outside $N_\varepsilon(\Sigma_\alpha)$ would be a zero
		of $F_n$ on $K_\varepsilon$. This proves the claim simultaneously
		for all $\alpha\in B$.
	\end{proof}
	
	\begin{remark}[Moment assumptions]\label{rem:fourth-moments}
		At $\alpha=-1$, tails $\mathbb P(|\xi_{12}|>t)\asymp t^{-r}$ with
		$2<r<4$ force
		$\|M_n\|\geq n^{-1/2}\max_{i<j}|\xi_{ij}|\to\infty$ in probability.
		Since $p_{n,-1}$ is the characteristic polynomial of $M_n$, some zero
		then leaves every neighborhood of $\Sigma_{-1}=[-2,2]$. Thus no $(4-\varepsilon)$th moment assumption suffices for the zero-free region of
		Theorem~\ref{zerofree}, although the weak semicircle law needs only
		second moments. See, e.g., \cite{LeeYin2014} for more discussion.
	\end{remark}
	
	\subsection{Cauchy transforms and the zero measure.}
	Let $\zeta_{n,1}^{(\alpha)},\ldots,\zeta_{n,n}^{(\alpha)}$ be the
	zeros of $p_{n,\alpha}$, counted with multiplicity, and set
	$\nu_{n,\alpha}:=n^{-1}\sum_{j=1}^n\delta_{\zeta_{n,j}^{(\alpha)}}$.
	Its Cauchy transform is
	\begin{equation}\label{eq:empirical-cauchy-transform}
		g_{n,\alpha}(z)
		:=g_{\nu_{n,\alpha}}(z)
		=\frac1n\sum_{j=1}^n\frac1{z-\zeta_{n,j}^{(\alpha)}}
		=\frac1n\frac{p_{n,\alpha}'(z)}{p_{n,\alpha}(z)}
	\end{equation}
	away from its poles.
	
	\begin{corollary}[Cauchy transforms and the semicircle law]\label{semicircle}
		For every compact $\mathcal K\Subset\mathcal D$,
		\begin{equation}\label{eq:empirical-transform-convergence}
			\sup_{(\alpha,z)\in\mathcal K}
			|g_{n,\alpha}(z)-g_\alpha(z)|
			\xrightarrow{\mathbb P}0,
		\end{equation}
		where the supremum is $+\infty$ if a pole lies in $\mathcal K$.
		In particular, for each fixed $|\alpha|\leq1$,
		$\nu_{n,\alpha}\Rightarrow\mu_\alpha$ in probability: for every bounded
		continuous $f:\C\to\C$,
		\begin{equation}\label{eq:fixed-alpha-zero-law}
			\int f\,d\nu_{n,\alpha}\xrightarrow{\mathbb P}\int f\,d\mu_\alpha.
		\end{equation}
	\end{corollary}
	
	\begin{proof}
		Logarithmic differentiation gives
		\begin{equation}\label{eq:empirical-transform-decomposition}
			g_{n,\alpha}
			=h_{n,\alpha}
			+\frac1n\frac{\partial_zF_{n,\alpha}}{F_{n,\alpha}}
		\end{equation}
		where the quotient is defined. The deterministic transform $h_{n,\alpha}$
		converges to $g_\alpha$
		locally uniformly on $\mathcal U$, by Lemma~\ref{lem:hermite-zero-law}. Differentiation and division
		are continuous at the nonvanishing limit in
		Theorem~\ref{gaussian}, so
		$\sup_{(\alpha,z)\in\mathcal K}|\partial_zF_{n,\alpha}(z)/F_{n,\alpha}(z)|=O_{\mathbb P}(1)$; the
		probability of a zero on $\mathcal K$ tends to zero. This proves
		\eqref{eq:empirical-transform-convergence}.
		
		For \eqref{eq:fixed-alpha-zero-law}, fix $|\alpha|\leq1$ and a
		subsequence. Theorem~\ref{zerofree} with $\varepsilon=1/k$, and
		\eqref{eq:empirical-transform-convergence} on a countable compact
		exhaustion of $\{\alpha\}\times\Omega_\alpha$, give a further
		subsequence along which, almost surely, every fixed neighborhood of
		$\Sigma_\alpha$ eventually contains all zeros and
		$g_{n,\alpha}\to g_\alpha$ locally uniformly on $\Omega_\alpha$. Along
		it, $(\nu_{n,\alpha})$ is tight, and every weak limit is supported on
		$\Sigma_\alpha$ and has Cauchy transform $g_\alpha=g_{\mu_\alpha}$ off
		$\Sigma_\alpha$. A measure on a segment is determined by its Cauchy
		transform, as in the proof of Lemma~\ref{lem:hermite-zero-law}, so
		$\nu_{n,\alpha}\Rightarrow\mu_\alpha$ almost surely along this
		subsequence. Since every subsequence has such a further subsequence,
		\eqref{eq:fixed-alpha-zero-law} follows.
	\end{proof}
	
	For each fixed parameter, Corollary~\ref{semicircle} gives the
	semicircle law in probability. It in fact holds almost surely and
	uniformly in $|\alpha|\leq1$:
	
	\begin{theorem}[Uniform almost-sure semicircle law]\label{thm:almost-sure-semicircle}
		Realize $(M_n)_{n\geq1}$ on a common probability space, with any coupling
		across $n$. Almost surely, for every bounded continuous $f:\C\to\C$,
		\begin{equation}\label{eq:uniform-empirical-zero-law}
			\sup_{|\alpha|\leq1}
			\left|\int f\,d\nu_{n,\alpha}-\int f\,d\mu_\alpha\right|
			\longrightarrow0.
		\end{equation}
	\end{theorem}
	
	The proof in Appendix~\ref{app:almost-sure-semicircle} combines the
	second-moment bounds of Section~\ref{sec:coefficients}, extended to
	$z=\infty$, with the mean-value inequality for the subharmonic
	function $\log|F_{n,\alpha}|$, using the normalization
	$F_{n,\alpha}(\infty)=1$ for every $\alpha$.
	
	\section{Strategy of the proof}\label{sec:strategy}\label{sec:heuristic}
	This section explains the proof before giving it. Readers interested
	mainly in the method can read it directly after
	Section~\ref{sec:functional-limits}.
	
	The limit depends on $z$ only through $g_\alpha(z)$: we can write
	$\Psi_\alpha=Q_\alpha\circ g_\alpha$, where (Section~\ref{sec:limit-polynomial})
	\[
	Q_\alpha(u)=\exp\Bigl(\sum_{\ell\geq1}a_{\ell,\alpha}G_\ell u^\ell
	-\frac12\sum_{\ell\geq1}a_{\ell,\alpha}^2\sigma_\ell^2u^{2\ell}\Bigr),
	\qquad a_{\ell,\alpha}=(-1)^\ell\alpha^{\ell-1}.
	\]
	Steps~1 and~2 below reduce
	$F_{n,\alpha}$ to the partition function of a gas of centered polymers,
	Steps~3 and~4 identify its limit coefficient by coefficient, and
	Step~5 combines the two. For simplicity we describe Steps~3 and~4 for
	bounded entries; Section~\ref{sec:entry-cutoff} reduces the general
	case to this one by truncation.
	
	\emph{Step 1. Exact expansion (Section~\ref{sec:expansion}).} Each
	transposition $(i\,j)$ in the cycle decomposition contributes a factor
	$\alpha|M_{ij}|^2$ with mean $\alpha/n$. Each vertex lies in $n-1$
	transpositions, so these means add up to order one and cannot be
	treated as small. Splitting each such factor into its mean and a
	centered part, and summing the means exactly over matchings, gives
	\[
	F_{n,\alpha}(z)
	=\sum_{\mathcal F}A_{n,\alpha,d(\mathcal F)}(z)
	\prod_{p\in\mathcal F}a_{d(p),\alpha}X_p,
	\]
	a sum over vertex-disjoint families $\mathcal F$ of centered polymers
	$p$ with $d(p)$ vertices each and $d(\mathcal F)$ vertices in total:
	diagonal entries, with $X_p=M_{ii}$; edges, with $X_p=|M_{ij}|^2-1/n$;
	and simple cycles, with $X_p=2\operatorname{Re}$ of the product of
	entries along the cycle. The weights $A_{n,\alpha,m}$ are deterministic
	ratios of matching polynomials, and distinct families are orthogonal.
	
	\emph{Step 2. Second moments (Section~\ref{sec:coefficients}).} Each
	$A_{n,\alpha,m}$ is a product of Cauchy transforms of the zeros of
	smaller expected polynomials, so $A_{n,\alpha,m}\to g_\alpha^m$ at fixed
	$m$. Orthogonality turns second moments into series over degrees, and
	the degree-$m$ term is bounded, uniformly in $n$ and locally uniformly
	in $(\alpha,z)$, by a constant times the corresponding term of the
	limiting series. That series converges exactly when
	$|\tau_\alpha(z)|<1$, the condition defining $\mathcal D$. Hence
	\[
	F_{n,\alpha}(z)-Q_{n,\alpha}(g_\alpha(z))\longrightarrow0
	\quad\text{in }L^2,\qquad
	Q_{n,\alpha}(u):=\sum_{\mathcal F}\prod_{p\in\mathcal F}a_{d(p),\alpha}X_pu^{d(p)}.
	\]
	
	\emph{Step 3. Gaussian polymer sums (Sections~\ref{sec:entry-cutoff}
		and~\ref{sec:power-sums}).} The sums $S_{\ell,n}$ of the polymer
	variables with $\ell$ vertices converge jointly in law to the
	independent Gaussians $G_\ell$: in their joint moments, only pairings of identical polymers
	survive in the limit.
	
	\emph{Step 4. Identifying the coefficients
		(Section~\ref{sec:exponential-limit}).} Substitute $S_{\ell,n}$ for
	$G_\ell$ in the formula for $Q_\alpha$. In each coefficient, the terms
	on distinct vertices reproduce $Q_{n,\alpha}$, and the remaining terms,
	with overlapping or repeated polymers, are orthogonal to it. The substituted formula and
	$Q_{n,\alpha}$ have the same limiting second moment, so the remaining terms vanish in
	$L^2$:
	\[
	Q_{n,\alpha}(u)-\exp\Bigl(\sum_{\ell\geq1}a_{\ell,\alpha}S_{\ell,n}u^\ell
	-\frac12\sum_{\ell\geq1}a_{\ell,\alpha}^2\sigma_\ell^2u^{2\ell}\Bigr)
	\longrightarrow0
	\qquad\text{in }L^2\text{, coefficient by coefficient},
	\]
	and by Step~3 the exponential converges in law to $Q_\alpha(u)$,
	coefficient by coefficient. For example, at $\alpha=0$ the coefficient of $u^2$ is
	$\tfrac12(S_{1,n}^2-\sum_iM_{ii}^2)$ in $Q_{n,\alpha}$ and
	$\tfrac12(S_{1,n}^2-v)$ in the formula for $Q_\alpha$: the Wick term
	$-\tfrac12v$ replaces the excluded self-pairs, two copies of the same
	polymer. Both converge in law to $(G_1^2-v)/2$.
	
	\emph{Step 5. From coefficients to functions
		(Section~\ref{sec:assembly}).} Steps~3 and~4 give the limit of each
	coefficient of $Q_{n,\alpha}$. The bound of Step~2 makes the
	high-degree coefficients uniformly small, and Step~2 also shows that
	$F_{n,\alpha}$ is close to $Q_{n,\alpha}\circ g_\alpha$. Together these
	give the functional limit. Since every error is controlled in $L^2$,
	second moments converge as well.
	
	\emph{A second view.} Appendix~\ref{app:cluster-expansion} places the
	proof in the framework of cluster expansions for polymer gases. There,
	the mean factors of Step~1 are what prevent a standard truncated
	expansion: they attach to each centered polymer with total weight of
	order one. Summing them exactly amounts to summing rooted trees, as in
	the trace method, and produces the factors $g_\alpha(z)^{d}$ of Step~2.
	What remains behaves like an ideal gas: the constraint that distinct
	polymers be vertex-disjoint becomes negligible, and repeated polymers
	give the Wick correction of Step~4.
	
	\section{Exact expansion in centered polymers}\label{sec:expansion}
	This section carries out Step~1 of Section~\ref{sec:heuristic}. From
	now on, $M=M_n$; all identities below are exact and hold for every
	$\alpha,z\in\C$. Splitting the fixed-point and transposition factors in
	Lemma~\ref{lem:cycle-decomposition} into their means and centered
	fluctuations gives
	\begin{equation}\label{eq:cycle-split}
		p_{n,\alpha}(z)
		=\sum_{\sigma\in\mathfrak S_n}
		\prod_{i\in\operatorname{Fix}(\sigma)}\!(z-M_{ii})
		\prod_{\{i,j\}\in\mathfrak m(\sigma)}\!
		\left(\frac\alpha n+\alpha\Bigl(|M_{ij}|^2-\frac1n\Bigr)\right)
		\prod_{\vec c\in C_{\geq3}(\sigma)}
		(-1)^{|\vec c|}\alpha^{|\vec c|-1}M_{\vec c},
	\end{equation}
	where the factors of longer cycles are already centered. Expanding the
	products gives centered pieces, the polymers of
	Section~\ref{sec:centered-polymers}, and mean parts on the remaining
	vertices. Summing the mean parts over matchings of those vertices
	gives the exact expansion of Proposition~\ref{prop:cycle-expansion}.
	The resulting matching polynomials $\Her_r(z;-\alpha/n)$ of the smaller
	complete graphs $K_r$ keep the edge weight $\alpha/n$ rather than
	$\alpha/r$.
	
	\subsection{Centered polymers and the exact expansion.}\label{sec:centered-polymers}
	A \emph{polymer} is one of the following
	objects $p$, with vertex set $V(p)$, degree $d(p):=|V(p)|$, and
	associated centered real random variable $X_p$:
	\begin{enumerate}[label=(\roman*)]
		\item a \emph{singleton} $\{i\}$, with $d=1$ and
		$X_{\{i\}}:=M_{ii}=\eta_i/\sqrt n$;
		\item an unordered edge $e=\{i,j\}$, with $d=2$ and
		$X_e:=|M_{ij}|^2-1/n$;
		\item an unoriented simple cycle $c$ of length $\ell\geq3$, with
		$d=\ell$ and $X_c:=M_{\vec c}+M_{\vec c^{\,-1}}=2\operatorname{Re}M_{\vec c}$,
		where $\vec c$ is either orientation of $c$; the second equality
		holds because Hermitian symmetry gives
		$M_{\vec c^{\,-1}}=\overline{M_{\vec c}}$.
	\end{enumerate}
	Let $\mathcal P_{\ell,n}$ be the set of polymers of degree $\ell$ on
	$[n]$. The cycle decomposition assigns each polymer of degree $\ell$ the
	deterministic weight
	\begin{equation}\label{eq:polymer-weights}
		a_{\ell,\alpha}:=(-1)^\ell\alpha^{\ell-1},\qquad \ell\geq1,
	\end{equation}
	with $\alpha^0:=1$. A family
	$\mathcal F$ of polymers is \emph{compatible} if its members have
	pairwise disjoint vertex sets. Put
	\[
	V(\mathcal F):=\bigcup_{p\in\mathcal F}V(p),
	\qquad
	d(\mathcal F):=\sum_{p\in\mathcal F}d(p),
	\qquad
	X_{\mathcal F}:=\prod_{p\in\mathcal F}X_p,
	\]
	write $|\mathcal F|$ for the number of polymers in $\mathcal F$, and set
	\begin{equation}\label{eq:family-weights}
		a_{\mathcal F,\alpha}:=\prod_{p\in\mathcal F}a_{d(p),\alpha}
		=(-1)^{d(\mathcal F)}\alpha^{d(\mathcal F)-|\mathcal F|}.
	\end{equation}
	\begin{proposition}[Exact matching resummation]\label{prop:cycle-expansion}
		For all $\alpha,z\in\C$,
		\begin{equation}\label{eq:exact-expansion}
			p_{n,\alpha}(z)
			=\sum_{\mathcal F\ \mathrm{compatible}}
			\Her_{n-d(\mathcal F)}(z;-\alpha/n)\,a_{\mathcal F,\alpha}X_{\mathcal F}.
		\end{equation}
		Equivalently, grouping by total degree,
		\begin{equation}\label{eq:polynomial-degree-expansion}
			p_{n,\alpha}(z)=\sum_{m=0}^n\Her_{n-m}(z;-\alpha/n)\,q_{n,\alpha,m},
			\qquad
			q_{n,\alpha,m}:=\sum_{\substack{\mathcal F\ \mathrm{compatible}\\
					d(\mathcal F)=m}}a_{\mathcal F,\alpha}X_{\mathcal F}.
		\end{equation}
	\end{proposition}
	
	\begin{proof}
		In \eqref{eq:cycle-split}, the centered parts are
		$-M_{ii}=a_{1,\alpha}X_{\{i\}}$ and
		$\alpha(|M_{ij}|^2-1/n)=a_{2,\alpha}X_e$, and the factor of a cycle
		$\vec c$ of length $\ell\geq3$ is $a_{\ell,\alpha}M_{\vec c}$. Expanding
		the products, a term is specified by $\sigma$ and a choice of mean or
		centered part at each fixed point and transposition. Such terms
		correspond bijectively to triples consisting of
		\begin{itemize}
			\item the family $\mathcal F$ of the chosen centered singletons and
			edges and the unoriented cycles of $\sigma$ of length at least
			three, which is compatible because the cycles of $\sigma$ are
			disjoint;
			\item an orientation of each cycle of $\mathcal F$;
			\item the matching $\mathfrak m_0$ of $[n]\setminus V(\mathcal F)$
			formed by the transpositions at which the mean was chosen.
		\end{itemize}
		The term equals $a_{\mathcal F,\alpha}$ times the variables $X_p$ of
		the singletons and edges of $\mathcal F$, the monomials $M_{\vec c}$ of
		its oriented cycles, and the mean weight
		$z^{n-d(\mathcal F)-2|\mathfrak m_0|}(\alpha/n)^{|\mathfrak m_0|}$ of
		$\mathfrak m_0$. Summing over the two orientations of each cycle
		replaces $M_{\vec c}$ by $M_{\vec c}+M_{\vec c^{\,-1}}=X_c$, and summing
		the mean weights over all matchings $\mathfrak m_0$ gives the matching
		polynomial $\Her_{n-d(\mathcal F)}(z;-\alpha/n)$ of the free vertices.
		This proves \eqref{eq:exact-expansion}, and
		\eqref{eq:polynomial-degree-expansion} follows by grouping the
		families by $d(\mathcal F)$.
	\end{proof}
	
	\subsection{Orthogonality and the degree isometry.}
	Distinct compatible families are orthogonal, so second moments of the
	expansion reduce to sums of $\E X_{\mathcal F}^2$ over single
	families. We compute these exactly; the resulting degree isometry
	drives the second-moment arguments of Section~\ref{sec:coefficients}.
	
	The independent \emph{matrix variables} are the diagonal entries
	$\eta_i$ and the off-diagonal variables $\xi_e$, one for each unordered
	edge $e$. A singleton $\{i\}$ uses $\eta_i$; an edge or cycle polymer
	uses $\xi_e$ for each of its edges $e$. For a matrix variable
	$\epsilon$, let $\mathcal M_\epsilon$ be the $\sigma$-algebra generated
	by all other matrix variables.
	
	\begin{lemma}[Conditional centering]\label{lem:conditional-centering}
		If a polymer $p$ uses the matrix variable $\epsilon$, then
		$\E[X_p\mid\mathcal M_\epsilon]=0$.
	\end{lemma}
	
	\begin{proof}
		A singleton gives $X_{\{i\}}=\eta_i/\sqrt n$, and an edge gives
		$X_e=(|\xi_e|^2-1)/n$; both depend only on $\epsilon$ and have mean
		zero. A simple cycle $c$ uses $\epsilon=\xi_e$, $e=\{i,j\}$, once in
		each orientation, so $X_c=M_{ij}Z+M_{ji}\overline Z$ with $Z$
		$\mathcal M_\epsilon$-measurable, and $\E M_{ij}=\E M_{ji}=0$.
	\end{proof}
	
	\begin{lemma}[Orthogonality]\label{lem:family-orthogonality}
		If $\mathcal F\neq\mathcal G$ are compatible families, then
		$\E[X_{\mathcal F}X_{\mathcal G}]=0$. In particular,
		$\E X_{\mathcal F}=0$ for every compatible $\mathcal F\neq\varnothing$.
	\end{lemma}
	
	\begin{proof}
		The matrix variables used by a compatible
		family $\mathcal F$ determine it: the singletons are read off from
		the diagonal entries, and the other polymers are the connected
		components of the graph formed by the edges, which by compatibility
		are isolated edges and simple cycles. Thus, if
		$\mathcal F\neq\mathcal G$, we may assume after interchanging the two
		families that some matrix variable $\epsilon$ used by $\mathcal F$ is
		not used by $\mathcal G$. By compatibility, $\epsilon$ is used by
		exactly one polymer $p\in\mathcal F$, so all factors of
		$X_{\mathcal F}X_{\mathcal G}$ other than $X_p$ are
		$\mathcal M_\epsilon$-measurable. Lemma~\ref{lem:conditional-centering}
		gives $\E[X_{\mathcal F}X_{\mathcal G}]=0$.
	\end{proof}
	
	\begin{lemma}[Polymer variances]\label{lem:cycle-variance-average}
		For every $B\subset[n]$ of size $\ell\geq1$, the centered variables
		satisfy
		\begin{equation}\label{eq:cycle-variance-block}
			\sum_{\substack{p\in\mathcal P_{\ell,n}\\V(p)=B}}\E X_p^2
			=\frac{\ell!}{n^\ell}\sigma_\ell^2.
		\end{equation}
	\end{lemma}
	
	\begin{proof}
		For an oriented simple cycle $\vec c$ of length $\ell\geq3$, let
		$\operatorname{asc}(\vec c)$ be the number of its upward edges
		$i\to j$ with $i<j$, that is, its cyclic ascents. Independence
		gives
		\[
		\E M_{\vec c}^2=n^{-\ell}\theta^{\operatorname{asc}(\vec c)}
		\overline\theta^{\,\ell-\operatorname{asc}(\vec c)},
		\qquad \E|M_{\vec c}|^2=n^{-\ell}.
		\]
		Since $X_c=M_{\vec c}+\overline{M_{\vec c}}$,
		\[
		\E X_c^2=2n^{-\ell}
		\left(1+\operatorname{Re}
		\bigl(\theta^{\operatorname{asc}(\vec c)}
		\overline\theta^{\,\ell-\operatorname{asc}(\vec c)}\bigr)\right).
		\]
		Relabel $B$ increasingly as $[\ell]$, which preserves ascents. Each
		oriented cycle on $[\ell]$ arises from exactly $\ell$ permutations
		$\pi\in\mathfrak S_\ell$, its rotations, all with the same number of
		cyclic ascents. Averaging over oriented cycles is therefore averaging
		over $\mathfrak S_\ell$ in \eqref{eq:cycle-order-average}. Reversal
		exchanges $\operatorname{asc}$ and $\ell-\operatorname{asc}$ and
		conjugates the weight, so the real part is unchanged, and its average
		over unoriented cycles is also $\rho_\ell(\theta)$.
		On a fixed set $B$ of $\ell\geq3$ vertices, there are
		$(\ell-1)!$ oriented cycles and half as many unoriented cycles.
		Summing the displayed second moments over the latter gives
		$(\ell-1)!n^{-\ell}(1+\rho_\ell(\theta))$, which is
		\eqref{eq:cycle-variance-block}. At $\ell=2$, the centered edge
		variable has variance $\kappa n^{-2}=2!n^{-2}\sigma_2^2$, and at
		$\ell=1$ the diagonal variable has variance $vn^{-1}=1!\,n^{-1}\sigma_1^2$.
	\end{proof}
	
	The summands of the \emph{polymer sums}
	$S_{\ell,n}:=\sum_{p\in\mathcal P_{\ell,n}}X_p$ are pairwise orthogonal
	by Lemma~\ref{lem:family-orthogonality}, so summing
	\eqref{eq:cycle-variance-block} over the $\binom n\ell$ sets $B$ gives
	\[
	\Var(S_{\ell,n})=\frac{(n)_\ell}{n^\ell}\,\sigma_\ell^2\longrightarrow\sigma_\ell^2.
	\]
	They are the finite-$n$ counterparts of the coefficients $G_\ell$ of
	the Gaussian series $S$ in \eqref{eq:common-gaussian-field}.
	
	For $\alpha,\gamma\in\C$, define coefficients
	$c_m(\alpha,\gamma)$, which Lemma~\ref{lem:orthogonality} identifies
	as the second moments of the degree-$m$ parts of the expansion, by
	\begin{equation}\label{eq:degree-second-moment}
		\sum_{m\geq0}c_m(\alpha,\gamma)s^m
		:=\exp\Bigl(\sum_{\ell\geq1}\sigma_\ell^2
		(\alpha\overline\gamma)^{\ell-1}s^\ell\Bigr),
	\end{equation}
	with $0^0:=1$, and write $c_m(\alpha):=c_m(\alpha,\alpha)$ for the
	nonnegative diagonal coefficients; for example, $c_m(0)=v^m/m!$.
	
	\begin{lemma}[Degree isometry]\label{lem:orthogonality}
		\begin{enumerate}[label=(\alph*)]
			\item For all $\alpha,\gamma\in\C$ and $0\leq m,m'\leq n$,
			\begin{equation}\label{eq:degree-isometry}
				\E[q_{n,\alpha,m}\overline{q_{n,\gamma,m'}}]
				=\mathbf 1_{\{m=m'\}}\frac{(n)_m}{n^m}c_m(\alpha,\gamma).
			\end{equation}
			\item For all $\alpha,\gamma,z,w\in\C$,
			\begin{equation}\label{eq:finite-n-polynomial-moment}
				\E\bigl[p_{n,\alpha}(z)\overline{p_{n,\gamma}(w)}\bigr]
				=\sum_{m=0}^n\frac{(n)_m}{n^m}c_m(\alpha,\gamma)
				\Her_{n-m}(z;-\alpha/n)\overline{\Her_{n-m}(w;-\gamma/n)}.
			\end{equation}
		\end{enumerate}
	\end{lemma}
	
	\begin{proof}
		By Lemma~\ref{lem:family-orthogonality}, only identical families contribute
		to $\E[q_{n,\alpha,m}\overline{q_{n,\gamma,m'}}]$, so different
		degrees are orthogonal. For a compatible family, independence gives
		$\E X_{\mathcal F}^2=\prod_{p\in\mathcal F}\E X_p^2$. A family with
		$k_\ell$ polymers of degree $\ell$, where $\sum_\ell\ell k_\ell=m$,
		is obtained by choosing pairwise disjoint blocks in $[n]$,
		$k_\ell$ of size $\ell$ for each $\ell$, and a polymer of degree
		$\ell$ on each block of size $\ell$. There are
		\[
		\frac{(n)_m}{\prod_{\ell}(\ell!)^{k_\ell}k_\ell!}
		\]
		such choices of blocks, and \eqref{eq:cycle-variance-block} sums the variances
		of the polymers on a block of size $\ell$ to $\ell!\sigma_\ell^2/n^\ell$.
		By \eqref{eq:family-weights}, each degree-$\ell$ polymer contributes
		the weight $(\alpha\overline\gamma)^{\ell-1}$ to
		$a_{\mathcal F,\alpha}\overline{a_{\mathcal F,\gamma}}$. Hence
		\[
		\E[q_{n,\alpha,m}\overline{q_{n,\gamma,m'}}]
		=\mathbf 1_{\{m=m'\}}\frac{(n)_m}{n^m}
		\sum_{\substack{(k_\ell)_{\ell\geq1}\\\sum_\ell\ell k_\ell=m}}
		\prod_{\ell\geq1}
		\frac{\bigl(\sigma_\ell^2(\alpha\overline\gamma)^{\ell-1}\bigr)^{k_\ell}}{k_\ell!}.
		\]
		Expanding the exponential in \eqref{eq:degree-second-moment} as a
		product over $\ell\geq1$ shows that the last sum is the coefficient of
		$s^m$ there, namely $c_m(\alpha,\gamma)$, proving (a). Part (b) follows by applying (a)
		to two copies of \eqref{eq:polynomial-degree-expansion}.
	\end{proof}
	
	Equivalently, for deterministic coefficients $a_m$,
	\[
	\E\Bigl|\sum_{m=0}^na_mq_{n,\alpha,m}\Bigr|^2
	=\sum_{m=0}^n|a_m|^2\frac{(n)_m}{n^m}c_m(\alpha):
	\]
	the map $(a_m)\mapsto\sum_ma_mq_{n,\alpha,m}$ is an isometry from a
	weighted $\ell^2$ space into $L^2$.
	
	Since $c_m(1)=c_m(-1)$ and $\Her_k(iz;1/n)=i^k\Her_k(z;-1/n)$,
	\eqref{eq:finite-n-polynomial-moment} gives
	\[
	\E\bigl[p_{n,1}(z)\overline{p_{n,1}(w)}\bigr]
	=\E\bigl[p_{n,-1}(iz)\overline{p_{n,-1}(iw)}\bigr].
	\]
	Fyodorov proved this identity for the GUE \cite[Eq.~(1.15)]{Fyodorov2006},
	and G{\"o}tze and K{\"o}sters extended it to their Hermitian Wigner ensemble
	\cite[Prop.~5.4]{GotzeKosters2009}.

	\section{Second moments and the centered polymer gas}\label{sec:coefficients}
	This section carries out Step~2 of Section~\ref{sec:heuristic}.
	We first find the limit of the matching ratios at each fixed degree.
	A resolvent representation then bounds all degrees uniformly in $n$.
	These bounds give convergence of the second moments of $F_{n,\alpha}$
	on $\mathcal D$ (Proposition~\ref{prop:second-moments}) and, applied
	to $F_{n,\alpha}(z)-Q_{n,\alpha}(g_\alpha(z))$, approximate
	$F_{n,\alpha}$ by the centered polymer gas on compact subsets of $\mathcal D$
	(Proposition~\ref{prop:joint-polymer-approximation}).
	
	\subsection{Matching ratios.}
	\label{sec:hermite-ratios}\label{sec:joint-bounds}
	By Lemma~\ref{lem:hermite-zero-law}, $P_{n,\alpha}$ has no zeros on
	$\Omega_\alpha$, so we may normalize. Dividing the matching polynomial
	$\Her_{n-m}(\cdot\,;-\alpha/n)$ of $K_{n-m}$ by $P_{n,\alpha}$ gives the
	\emph{matching ratios}
	\begin{equation}\label{eq:matching-ratio-definition}
		A_{n,\alpha,m}(z)
		:=\frac{\Her_{n-m}(z;-\alpha/n)}{P_{n,\alpha}(z)}
		=\frac{\Her_{n-m}(z;-\alpha/n)}{\Her_n(z;-\alpha/n)}
		=\prod_{j=0}^{m-1}
		\frac{\Her_{n-j-1}(z;-\alpha/n)}{\Her_{n-j}(z;-\alpha/n)},
	\end{equation}
	for $0\leq m\leq n$, by Lemma~\ref{lem:annealed-polynomial}. These
	compare the weight of the mean matchings on $n-m$ free vertices
	with that on all $n$ vertices. At $\alpha=0$, $\Her_r(z;0)=z^r$, so the
	matching ratios already equal their limits: $A_{n,0,m}=g_0^m$. For
	$z\in\Omega_\alpha$, dividing \eqref{eq:polynomial-degree-expansion} by
	$P_{n,\alpha}$ gives
	\begin{equation}\label{eq:degree-expansion}
		F_{n,\alpha}(z)=\sum_{m=0}^nA_{n,\alpha,m}(z)q_{n,\alpha,m}.
	\end{equation}
	
	\begin{lemma}[Fixed-degree limit]\label{lem:fixed-degree-limit}
		For every fixed $m\geq0$,
		\begin{equation}\label{eq:fixed-degree-ratio-limit}
			A_{n,\alpha,m}(z)\longrightarrow g_\alpha(z)^m
			\qquad\text{locally uniformly on }\mathcal U.
		\end{equation}
	\end{lemma}
	
	\begin{proof}
		For $1\leq r\leq n$ put $\alpha_r:=\alpha r/n$, so that
		$-\alpha/n=-\alpha_r/r$. With $r=n-j$, the $j$th factor of the
		product in \eqref{eq:matching-ratio-definition} is
		\[
		\frac{\Her_{r-1}(z;-\alpha_r/r)}{\Her_r(z;-\alpha_r/r)}
		=h_{r,\alpha_r}(z),
		\]
		the Cauchy transform of the zeros of the smaller expected polynomial
		$P_{r,\alpha_r}$. Fix a compact $\mathcal K\Subset\mathcal U$ and a
		compact neighborhood $\mathcal K'\Subset\mathcal U$ of it. For each
		fixed $j<m$, we have $|\alpha_{n-j}-\alpha|\leq j|\alpha|/n\to0$, so
		$(\alpha_{n-j},z)\in\mathcal K'$ for all $(\alpha,z)\in\mathcal K$
		once $n$ is large. By \eqref{eq:hermite-transform-limit},
		$h_{r,\beta}(z)\to g_\beta(z)$ uniformly on $\mathcal K'$, and $g$ is
		uniformly continuous there, so each factor converges to $g_\alpha(z)$
		uniformly on $\mathcal K$. A product of $m$ uniformly convergent
		factors with bounded limits converges uniformly.
	\end{proof}
	
	Lemma~\ref{lem:fixed-degree-limit} controls each term of the degree
	expansion \eqref{eq:degree-expansion} separately, but the expansion
	contains degrees up to $n$. To bound all degrees at once, we express the
	matching ratios as entries of the resolvent of a Jacobi matrix.
	
	Let $J_n$ be the real symmetric tridiagonal $n\times n$ matrix with zero
	diagonal and off-diagonal entries $(J_n)_{j,j+1}=(J_n)_{j+1,j}=\sqrt{j/n}$,
	$1\leq j<n$. As in the proof of Lemma~\ref{lem:hermite-zero-law},
	$\det(\zeta I-J_n)=\Her_n(\zeta;1/n)$ and the spectrum of $J_n$ lies in
	$(-2,2)$, so the \emph{resolvent} $(\zeta I-J_n)^{-1}$ is defined for
	$\zeta\in\C\setminus(-2,2)$. For $\beta^2=-\alpha$,
	\eqref{eq:hermite-scaling} with $c=\beta$ gives
	$\beta^r\Her_r(z/\beta;1/n)=\Her_r(z;-\alpha/n)$; in particular
	$P_{n,\alpha}(z)=\det(zI-\beta J_n)$.
	
	\begin{lemma}[Resolvent form of the matching ratios]\label{lem:joint-coefficients}
		Let $\alpha\neq0$, let $\beta$ be a square root of $-\alpha$, let
		$z\in\Omega_\alpha$, and put $\zeta:=z/\beta$. For $1\leq m\leq n$,
		\begin{equation}\label{eq:resolvent-matching-ratio}
			\bigl[(\zeta I-J_n)^{-1}\bigr]_{n-m+1,n}
			=\beta^{m}\sqrt{\frac{(n)_m}{n^m}}\,A_{n,\alpha,m}(z).
		\end{equation}
	\end{lemma}
	
	\begin{proof}
		Put $b_j:=\sqrt{j/n}$ and $H_r:=\Her_r(\zeta;1/n)
		=\beta^{-r}\Her_r(z;-\alpha/n)$, by \eqref{eq:hermite-scaling}. In
		particular $H_n=\det(\zeta I-J_n)=\beta^{-n}P_{n,\alpha}(z)\neq0$.
		By Cramer's rule, the $i$th entry of the last column of
		$(\zeta I-J_n)^{-1}$ is $(-1)^{i+n}D_i/H_n$, where $D_i$ is
		the determinant of $\zeta I-J_n$ with row $n$ and column $i$
		deleted. Since $J_n$ is tridiagonal, this matrix is block lower
		triangular. For instance, for $n=5$ and $i=3$,
		\[
		\zeta I-J_5=\begin{pmatrix}
			\zeta&-b_1&0&0&0\\
			-b_1&\zeta&-b_2&0&0\\
			0&-b_2&\zeta&-b_3&0\\
			0&0&-b_3&\zeta&-b_4\\
			0&0&0&-b_4&\zeta
		\end{pmatrix},
		\qquad
		D_3=\det\left(\begin{array}{cc|cc}
			\zeta&-b_1&0&0\\
			-b_1&\zeta&0&0\\\hline
			0&-b_2&-b_3&0\\
			0&0&\zeta&-b_4
		\end{array}\right).
		\]
		In general, the upper left block is the leading $(i-1)\times(i-1)$
		block of $\zeta I-J_n$, whose determinant is $H_{i-1}$ by the same
		recurrence as for $\det(\zeta I-J_n)$, and the lower right block is
		triangular with diagonal entries $-b_i,\dots,-b_{n-1}$. Hence
		$(-1)^{i+n}D_i=H_{i-1}\prod_{j=i}^{n-1}b_j$. For $i=n-m+1$ the
		entry is
		\[
		\frac{H_{n-m}}{H_n}\prod_{n-m<j<n}b_j
		=\beta^m\sqrt{\frac{(n)_m}{n^m}}\,
		\frac{\Her_{n-m}(z;-\alpha/n)}{\Her_n(z;-\alpha/n)},
		\]
		since $(n)_m/n^m=\prod_{n-m<j\leq n}j/n$ and the $j=n$ factor is $1$.
		This is \eqref{eq:resolvent-matching-ratio} by
		\eqref{eq:matching-ratio-definition}; it is the last-column case of
		the classical formula for the inverse of a tridiagonal matrix
		\cite{Usmani1994}.
	\end{proof}
	
	The identity \eqref{eq:resolvent-matching-ratio} reduces bounds on the
	matching ratios to bounds on resolvent entries at distance $m-1$ from
	the diagonal. Since $J_n$ is tridiagonal, $(J_n^k)_{i,j}=0$ whenever
	$|i-j|>k$, so the $(n-m+1,n)$ entry of a polynomial in $J_n$ of degree
	less than $m-1$ vanishes. If the resolvent is expanded in polynomials
	that stay bounded on the spectrum, only the terms of degree at least
	$m-1$ contribute, and each degree costs a factor of order $|g_\alpha|$.
	This is the polynomial-approximation argument of Demko, Moss and Smith
	\cite{DemkoMossSmith1984} for the decay of inverses of band matrices;
	here the Chebyshev expansion of the resolvent is explicit, which gives
	the rate exactly.
	
	\begin{lemma}[Uniform bound on the matching ratios]\label{lem:matching-ratio-bound}
		For each compact $\mathcal K\Subset\mathcal U$ there is
		$C_{\mathcal K}<\infty$ such that, for $(\alpha,z)\in\mathcal K$,
		$n\geq1$ and $0\leq m\leq n$,
		\begin{equation}\label{eq:weighted-matching-ratio}
			\frac{(n)_m}{n^m}|A_{n,\alpha,m}(z)|^2
			\leq C_{\mathcal K}|g_\alpha(z)|^{2m}.
		\end{equation}
	\end{lemma}
	
	\begin{proof}
		The case $m=0$ is trivial, and at $\alpha=0$,
		$A_{n,0,m}=g_0^m$ and $(n)_m/n^m\leq1$. Let $\alpha\neq0$ and
		$1\leq m\leq n$, and let $\beta$ and $\zeta$ be as in
		Lemma~\ref{lem:joint-coefficients}. We use the Chebyshev polynomials
		$T_k(\cos\vartheta)=\cos(k\vartheta)$. These are polynomials of degree
		$k$, by the recurrence $T_{k+1}(x)=2xT_k(x)-T_{k-1}(x)$, and satisfy
		$|T_k|\leq1$ on $[-1,1]$. For $y=2\cos\vartheta$ and $|\omega|<1$, summing the
		geometric series in $\omega e^{\pm i\vartheta}$ gives
		\begin{equation}\label{eq:chebyshev-kernel}
			\frac{1-\omega^2}{1-\omega y+\omega^2}
			=\frac1{1-\omega e^{i\vartheta}}+\frac1{1-\omega e^{-i\vartheta}}-1
			=1+2\sum_{k\geq1}T_k(y/2)\,\omega^k.
		\end{equation}
		Write $u=g_\alpha(z)$ and $\omega:=\beta u$, so that
		$|\omega|^2=|\alpha u^2|<1$ by Lemma~\ref{lem:conformal-disk-coordinate}.
		By \eqref{eq:cauchy-inverse}, $z=u^{-1}+\beta^2u$, so
		$\zeta=\omega+\omega^{-1}=J_{-1}(\omega)$, the Joukowski map of
		Section~\ref{sec:disk-coordinate}, and
		$\zeta-y=(1-\omega y+\omega^2)/\omega$.
		Multiplying \eqref{eq:chebyshev-kernel} by $\omega/(1-\omega^2)$ and
		applying the spectral theorem, with $\|T_k(J_n/2)\|\leq1$, gives
		\[
		(\zeta I-J_n)^{-1}
		=\frac{\omega}{1-\omega^2}\Bigl(I+2\sum_{k\geq1}\omega^k\,
		T_k(J_n/2)\Bigr),
		\]
		with the series converging in operator norm. Only the terms with
		$k\geq m-1$ contribute to the $(n-m+1,n)$ entry, so
		\eqref{eq:resolvent-matching-ratio} gives
		\[
		|\beta|^m\sqrt{\frac{(n)_m}{n^m}}\,|A_{n,\alpha,m}(z)|
		\leq\frac{2|\omega|}{|1-\omega^2|}\sum_{k\geq m-1}|\omega|^k
		=\frac{2|\omega|^m}{|1-\omega^2|(1-|\omega|)}.
		\]
		Since $|\omega|^m=|\beta|^m|u|^m$ and
		$|1-\omega^2|=|1+\alpha u^2|\geq1-|\alpha u^2|$,
		\[
		\sqrt{\frac{(n)_m}{n^m}}\,|A_{n,\alpha,m}(z)|
		\leq\frac{2|u|^m}{(1-|\alpha u^2|)(1-|\beta u|)},\qquad 1\leq m\leq n.
		\]
		On $\mathcal K$, $\sup|\alpha g_\alpha(z)^2|<1$ by
		compactness, which gives \eqref{eq:weighted-matching-ratio}.
	\end{proof}
	
	\subsection{Second-moment convergence.}\label{sec:uniform-tails}\label{sec:second-moments}
	By \eqref{eq:degree-expansion} and Lemma~\ref{lem:orthogonality}(a),
	for $(\alpha,z)\in\mathcal U$,
	\begin{equation}\label{eq:second-moment-identity}
		\E|F_{n,\alpha}(z)|^2
		=\sum_{m=0}^n\frac{(n)_m}{n^m}c_m(\alpha)|A_{n,\alpha,m}(z)|^2.
	\end{equation}
	Lemma~\ref{lem:fixed-degree-limit} gives the limit of each term. To
	identify the limit, for $(\alpha,z)\in\mathcal U$ put
	\begin{equation}\label{eq:limit-second-moment}
		V_\alpha(z):=\sum_{m\geq0}c_m(\alpha)|g_\alpha(z)|^{2m}
		=\exp\Bigl(\sum_{\ell\geq1}\sigma_\ell^2|\alpha|^{2\ell-2}
		|g_\alpha(z)|^{2\ell}\Bigr)\in(0,\infty],
	\end{equation}
	the second equality by \eqref{eq:degree-second-moment} with
	$s=|g_\alpha(z)|^2$.
	To bound the tails, allow $u$ to vary in place of $g_\alpha(z)$ and put
	$\mathcal V:=\{(\alpha,u)\in\C^2:|\alpha u|<1\}$.
	
	\begin{lemma}[Degree tails]\label{lem:degree-tails}
		For each compact $K\Subset\mathcal V$ there are $b>1$ and
		$C_K<\infty$ such that
		\begin{equation}\label{eq:joint-degree-tail}
			\sum_{m>L}c_m(\alpha)|u|^{2m}
			\leq C_K b^{-L},\qquad (\alpha,u)\in K,\quad L\geq0.
		\end{equation}
	\end{lemma}
	
	\begin{proof}
		Choose $b>1$ with $b\sup_K|\alpha u|^2<1$. Since the $c_m(\alpha)$ are
		nonnegative, \eqref{eq:degree-second-moment} gives
		\[
		\sum_{m>L}c_m(\alpha)|u|^{2m}
		\leq b^{-L}\sum_{m\geq0}c_m(\alpha)(b|u|^2)^m
		=b^{-L}\exp\!\Bigl(b|u|^2\sum_{\ell\geq1}\sigma_\ell^2
		(b|\alpha u|^2)^{\ell-1}\Bigr),
		\]
		and the exponential factor is bounded on $K$.
	\end{proof}
	
	For $u=g_\alpha(z)$ we have $|\alpha u|=|\tau_\alpha(z)|$, so the
	bound applies locally uniformly on $\mathcal D$.
	
	\begin{proposition}[Second-moment convergence]\label{prop:second-moments}
		$V_\alpha(z)<\infty$ on $\mathcal D$, and
		$\E|F_{n,\alpha}(z)|^2\to V_\alpha(z)$ locally uniformly on
		$\mathcal D$.
	\end{proposition}
	
	\begin{proof}
		On each compact $\mathcal K\Subset\mathcal D$,
		\eqref{eq:weighted-matching-ratio} and Lemma~\ref{lem:degree-tails}
		bound the tails beyond degree $L$ of both \eqref{eq:second-moment-identity}
		and \eqref{eq:limit-second-moment} by $C_{\mathcal K}b^{-L}$,
		uniformly in $n$. In particular $V_\alpha$ is finite there. Each
		fixed summand converges uniformly by \eqref{eq:fixed-degree-ratio-limit}.
		Hence, for every $L$,
		\[
		\limsup_{n\to\infty}\sup_{\mathcal K}\bigl|\E|F_{n,\alpha}(z)|^2-V_\alpha(z)\bigr|
		\leq2C_{\mathcal K}b^{-L},
		\]
		and letting $L\to\infty$ proves the claim.
	\end{proof}
	
	\begin{remark}[The domain as the region of finite second moment]
		\label{rem:second-moment-domain}
		The domain $\mathcal D$ is exactly the region of $\mathcal U$ where
		the second moments of $F_{n,\alpha}(z)$ stay bounded; for
		$|\alpha|>1$, its boundary in each $z$-slice is the ellipse of
		Appendix~\ref{app:exterior-regime}.
		
		On $\mathcal D$, Proposition~\ref{prop:second-moments} bounds the
		second moments. On $\mathcal U\setminus\mathcal D$, the terms of
		\eqref{eq:second-moment-identity} are nonnegative and converge
		termwise by Lemma~\ref{lem:fixed-degree-limit}, which holds on all of
		$\mathcal U$, so Fatou's lemma gives
		$\liminf_n\E|F_{n,\alpha}(z)|^2\geq V_\alpha(z)$, and $V_\alpha(z)=\infty$
		there. Indeed, for $\alpha\neq0$ the exponent in
		\eqref{eq:limit-second-moment} is $\mathcal C(|\tau_\alpha(z)|)/|\alpha|^2$,
		with the series interpreted in $[0,\infty]$ (at $\alpha=0$ it is
		$v/|z|^2$, and $\tau_0=0$). Let $B_N$ be the Hermitian matrix with
		zero diagonal, $\theta$ above it and $\overline\theta$ below it. Then
		$\operatorname{Tr}(B_N^\ell)$ sums
		$\theta^{\#\mathrm{up}}\overline\theta^{\,\#\mathrm{down}}$ over closed
		walks of length $\ell$, all but $O(N^{\ell-1})$ of which visit $\ell$
		distinct indices. Since the trace is real,
		\eqref{eq:cycle-order-average} gives
		$\rho_\ell(\theta)=\lim_NN^{-\ell}\operatorname{Tr}(B_N^\ell)$, so
		$\rho_{2k}(\theta)\geq0$. Hence $\sigma_{2k}^2\geq1/(2k)$ for $k\geq2$,
		and $\mathcal C(r)\geq\sum_{k\geq2}r^{4k}/(2k)=\infty$ for $r\geq1$.
	\end{remark}
	
\begin{figure}[tbp]
		\centering
		\includegraphics[width=0.9\textwidth]{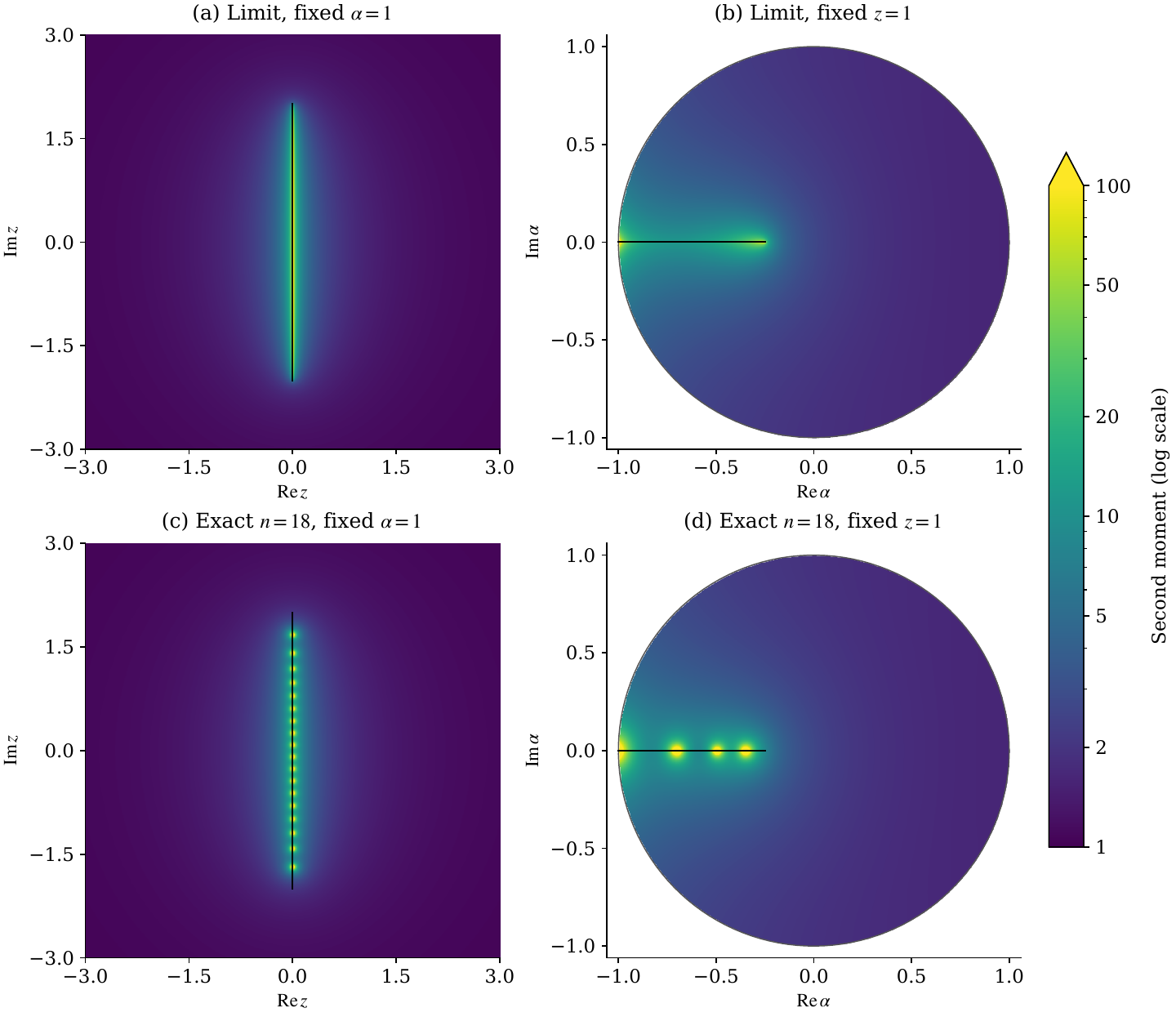}
		\caption{GUE limiting second moments $V_\alpha(z)=\E|\Psi_\alpha(z)|^2$ in (a,b), compared
			with the exact finite-$n$ moments $\E|F_{18,\alpha}(z)|^2$ in (c,d).
			The left column fixes $\alpha=1$; the right fixes $z=1$ and varies
			$\alpha$ over the unit disk. All panels share a logarithmic color
			scale, saturating above $100$. Black segments mark excluded
			parameters in both rows: $z\in i[-2,2]$ in (a,c) and
			$\alpha\in[-1,-1/4]$ in (b,d). The finite-$n$ values are exact
			ensemble moments computed from the matching expansion and
			orthogonality. The localized peaks in (c,d) occur near zeros of the
			normalizing expected polynomial $P_{18,\alpha}(z)$.}
		\label{fig:second-moments}
	\end{figure}

	Figure~\ref{fig:second-moments} compares $V_\alpha(z)$ with the exact
	finite-$n$ second moments for GUE, on the same parameter slices as
	Figure~\ref{fig:limiting-correlation}; the finite-$n$ moments are
	computed exactly from the orthogonality identity
	\eqref{eq:second-moment-identity}.
	
	\subsection{The centered polymer gas.}\label{sec:compact-polymer-approximation}\label{sec:polymer-polynomial}
	Replacing the matching ratios by their fixed-degree limits suggests
	approximating $F_{n,\alpha}(z)$ by $Q_{n,\alpha}(g_\alpha(z))$. The
	matching ratio $A_{n,\alpha,d(\mathcal F)}$, which multiplies
	$a_{\mathcal F,\alpha}X_{\mathcal F}$ after dividing
	\eqref{eq:exact-expansion} by $P_{n,\alpha}$, comes from the vertices left free
	by the whole family and does not generally factor over its polymers,
	whereas its limit $g_\alpha^{d(\mathcal F)}=\prod_{p\in\mathcal F}g_\alpha^{d(p)}$
	does. The limit therefore assigns each polymer its own activity, and
	the \emph{polymer polynomial}
	\begin{equation}\label{eq:polymer-polynomial}
		Q_{n,\alpha}(u)
		:=\sum_{m=0}^nq_{n,\alpha,m}u^m
		=\sum_{\mathcal F\ \mathrm{compatible}}
		a_{\mathcal F,\alpha}X_{\mathcal F}u^{d(\mathcal F)}
	\end{equation}
	is the \emph{partition function} of a gas. Polymer $p$
	has \emph{activity} $a_{d(p),\alpha}X_pu^{d(p)}$. A \emph{configuration} is a
	compatible family of polymers: the \emph{hard-core interaction}
	forbids two polymers from sharing a vertex. The partition function sums the product of the activities over all
	allowed configurations, with the empty configuration contributing $1$,
	and the variable $u$ counts occupied vertices. Section~\ref{sec:fixed-order}
	shows that, because the polymers are centered, the hard-core interaction
	between distinct polymers becomes negligible as $n\to\infty$---the gas
	becomes asymptotically ideal (see also
	Appendix~\ref{app:cluster-expansion}), which is why $Q_{n,\alpha}$
	converges to an exponential.
	
	The second-moment identities above are pointwise. To obtain the
	uniform statement of Proposition~\ref{prop:joint-polymer-approximation},
	we use the mean-value inequality to convert pointwise second moments
	into bounds on supremum norms.
	
	\begin{lemma}[Supremum norms]\label{lem:sup-norm}
		Let $U\subset\C^d$ be open, $d\geq1$, and let $K\subset U$ be
		compact. There are a compact set $K_+\subset U$ and $C<\infty$,
		depending only on $K$ and $U$, such that
		\[
		\|f\|_K^2\leq C\int_{K_+}|f(w)|^2\,dA(w),
		\qquad f\in\mathcal H(U),
		\]
		where $A$ is Lebesgue measure on $\C^d$. Consequently, every random
		holomorphic function $f$ on $U$ satisfies
		\[
		\E\|f\|_K^2\leq C\,A(K_+)\sup_{w\in K_+}\E|f(w)|^2.
		\]
	\end{lemma}
	
	\begin{proof}
		Choose $\delta>0$ such that the set $K_+$ of points whose coordinates
		are all within distance $\delta$ of those of some point of $K$ is
		contained in $U$. For $w\in K$, the closed polydisk
		$P(w,\delta):=\{w'\in\C^d:|w'_k-w_k|\leq\delta\text{ for }1\leq k\leq d\}$
		lies in $K_+$. Since $|f|^2$ is subharmonic in each
		variable separately, iterating the mean-value inequality gives
		\[
		|f(w)|^2\leq\frac1{(\pi\delta^2)^d}\int_{P(w,\delta)}|f|^2\,dA
		\leq\frac1{(\pi\delta^2)^d}\int_{K_+}|f|^2\,dA.
		\]
		Taking the supremum over $w\in K$ proves the deterministic bound.
		Taking expectations and applying Tonelli's theorem gives the
		second bound.
	\end{proof}

	\begin{proposition}[Polymer approximation]\label{prop:joint-polymer-approximation}
		For every compact $\mathcal K\Subset\mathcal D$,
		\begin{equation}\label{eq:joint-polymer-approximation}
			\E\sup_{(\alpha,z)\in\mathcal K}
			\bigl|F_{n,\alpha}(z)-Q_{n,\alpha}(g_\alpha(z))\bigr|^2
			\longrightarrow0.
		\end{equation}
	\end{proposition}
	
	\begin{proof}
		By \eqref{eq:degree-expansion} and Lemma~\ref{lem:orthogonality}(a),
		\begin{equation}\label{eq:replacement-identity}
			\E\bigl|F_{n,\alpha}(z)-Q_{n,\alpha}(g_\alpha(z))\bigr|^2
			=\sum_{m=0}^n\frac{(n)_m}{n^m}c_m(\alpha)
			\bigl|A_{n,\alpha,m}(z)-g_\alpha(z)^m\bigr|^2.
		\end{equation}
		On each compact $\mathcal K\Subset\mathcal D$, every fixed summand tends
		uniformly to zero by Lemma~\ref{lem:fixed-degree-limit}. Since
		$(n)_m/n^m\leq1$, \eqref{eq:weighted-matching-ratio} bounds the
		$m$th summand by $2(C_{\mathcal K}+1)c_m(\alpha)|g_\alpha(z)|^{2m}$.
		Its tails vanish uniformly by Lemma~\ref{lem:degree-tails}. Thus the
		pointwise second moments tend to zero locally uniformly. Applying
		Lemma~\ref{lem:sup-norm} to the holomorphic difference gives
		\eqref{eq:joint-polymer-approximation}.
	\end{proof}
	
	\section{Gaussian polymer sums and the functional limit}\label{sec:fixed-order}
	This section proves Theorem~\ref{gaussian} and
	Proposition~\ref{prop:wasserstein-convergence}.
	Section~\ref{sec:limit-polynomial} reduces them to convergence of the
	coefficients of $Q_{n,\alpha}$ (Proposition~\ref{prop:coefficient-limit}).
	Section~\ref{sec:entry-cutoff} reduces that to bounded entries, and
	Sections~\ref{sec:power-sums} and~\ref{sec:exponential-limit} prove it
	there, using the joint moments of the polymer sums and orthogonality.
	Section~\ref{sec:assembly} combines coefficient convergence with the
	geometric degree bound of Lemma~\ref{lem:degree-tails} to prove the
	functional and quadratic Wasserstein limits.
	
	\subsection{Limiting function and coefficient convergence.}\label{sec:limit-polynomial}
	The limit depends on $z$ only through $g_\alpha(z)$. For
	$(\alpha,u)\in\mathcal V$, define
	\begin{equation}\label{eq:limit-polynomial}
		Q_\alpha(u):=\exp\Bigl(\sum_{\ell\geq1}a_{\ell,\alpha}G_\ell u^\ell
		-\frac12\sum_{\ell\geq1}a_{\ell,\alpha}^2\sigma_\ell^2u^{2\ell}\Bigr).
	\end{equation}
	Since $a_{\ell,\alpha}u^\ell=-u(-\alpha u)^{\ell-1}$, both series
	converge absolutely and locally uniformly on $\{|\alpha u|<1\}$, on the
	probability-one event on which $S$ is defined
	(Section~\ref{sec:functional-limit}), so $Q_\alpha(u)$ is jointly
	holomorphic. By the expansion \eqref{eq:mode-expansion} of $H_\alpha$
	in powers of $g_\alpha$, $\Psi_\alpha=Q_\alpha\circ g_\alpha$ on
	$\mathcal D$.
	
	Put $b_{j,m}:=[\alpha^j u^m]Q_\alpha$ and
	\begin{equation}\label{eq:finite-coefficient}
		b_{n,j,m}:=[\alpha^j u^m]Q_{n,\alpha}
		=(-1)^m\sum_{\substack{\mathcal F\ \mathrm{compatible}\\
				d(\mathcal F)=m,\ |\mathcal F|=m-j}}X_{\mathcal F}.
	\end{equation}
	The second equality holds because, in \eqref{eq:polymer-polynomial}, a
	family of degree $m$ has weight
	$a_{\mathcal F,\alpha}=(-1)^m\alpha^{m-|\mathcal F|}$ by
	\eqref{eq:family-weights}, so it contributes to $\alpha^j$ exactly when
	it has $m-j$ polymers. In particular
	$q_{n,\alpha,m}=\sum_j\alpha^jb_{n,j,m}$, and for each $m$ only
	$0\leq j\leq m$ occur, by \eqref{eq:finite-coefficient} for
	$Q_{n,\alpha}$ and \eqref{eq:limit-polynomial} for $Q_\alpha$.
	
	\begin{proposition}[Coefficient limit]\label{prop:coefficient-limit}
		For $L\geq0$, put $b_n^{(\leq L)}:=(b_{n,j,m})_{0\leq j\leq m\leq L}$
		and $b^{(\leq L)}:=(b_{j,m})_{0\leq j\leq m\leq L}$. Then
		\[
		W_2\bigl(b_n^{(\leq L)},b^{(\leq L)}\bigr)\longrightarrow0.
		\]
		In particular,
		$b_n^{(\leq L)}\Rightarrow b^{(\leq L)}$.
	\end{proposition}
	
	The proof occupies Sections~\ref{sec:entry-cutoff}--\ref{sec:exponential-limit}.

	\subsection{Truncation.}\label{sec:entry-cutoff}
	The moment computation of Section~\ref{sec:power-sums} needs bounded
	entries. To reduce to that case, we truncate the entries at level $R$
	and recenter and renormalize them. Lemma~\ref{lem:truncation} below
	controls the error, and Step~2 of the proof of
	Proposition~\ref{prop:coefficient-limit} removes the truncation. The truncation is harmless for two reasons.
	The truncated entries have the same means and variances, so the
	expected polynomial $P_{n,\alpha}$ is unchanged. And the truncation
	errors of distinct families are orthogonal, so the error in each
	coefficient of $Q_{n,\alpha}$ stays small uniformly in $n$.
	
	With $\xi:=\xi_{12}$, put
	\[
	\zeta_R:=\xi\mathbf1_{\{|\xi|\leq R\}}
	-\E[\xi\mathbf1_{\{|\xi|\leq R\}}],\qquad
	s_R^2:=\E|\zeta_R|^2,\qquad
	\xi^{(R)}:=\zeta_R/s_R.
	\]
	Since $\E\xi=0$, we have
	$\|\zeta_R-\xi\|_{L^4}\leq2\|\xi\mathbf1_{\{|\xi|>R\}}\|_{L^4}\to0$.
	Hence $s_R\to1$, so $\xi^{(R)}$ is well-defined and bounded for
	large $R$, and $\|\xi^{(R)}-\xi\|_{L^4}\to0$. By construction,
	$\E\xi^{(R)}=0$ and $\E|\xi^{(R)}|^2=1$, so the truncated entries
	satisfy the assumptions of Section~\ref{sec:intro}.
	
	Similarly, with $\eta:=\eta_1$, let $\eta^{(R)}$ be the centered
	truncation $\eta\mathbf1_{\{|\eta|\leq R\}}-\E[\eta\mathbf1_{\{|\eta|\leq R\}}]$,
	rescaled to have variance $v$ for large $R$. (If $v=0$, then $\eta=0$
	and we take $\eta^{(R)}=0$.) The same argument in $L^2$ gives
	$\|\eta^{(R)}-\eta\|_{L^2}\to0$.
	
	For every unordered edge $e$ and every vertex $i$, apply these
	truncations to $\xi_e$ and $\eta_i$, and let $M_n^{(R)}$ be the
	Hermitian matrix with off-diagonal entries $\xi_e^{(R)}/\sqrt n$ and
	diagonal entries $\eta_i^{(R)}/\sqrt n$. Let $X_p^{(R)}$,
	$X_{\mathcal F}^{(R)}$ and $Q_{n,\alpha}^{(R)}$ be its centered
	variables, family products and polymer polynomial, and put
	$b^{(R)}_{n,j,m}:=[\alpha^ju^m]Q^{(R)}_{n,\alpha}$. The corresponding mode variances
	$(\sigma_\ell^{(R)})^2$ are given by \eqref{eq:mode-variances} with
	$\theta_R:=\E[(\xi^{(R)})^2]$ and $\kappa_R:=\Var(|\xi^{(R)}|^2)$ in
	place of $\theta$ and $\kappa$.
	
	\begin{lemma}[Truncation]\label{lem:truncation}
		As $R\to\infty$:
		\begin{enumerate}[label=(\alph*)]
			\item $\sup_n\E\bigl|b_{n,j,m}-b^{(R)}_{n,j,m}\bigr|^2\to0$ for all
			$0\leq j\leq m$;
			\item $(\sigma_\ell^{(R)})^2\to\sigma_\ell^2$ for all $\ell\geq1$.
		\end{enumerate}
	\end{lemma}
	
	\begin{proof}
		(a) The case $m=0$ is trivial, so let $m\geq1$. Put
		$\Delta_{\mathcal F}:=X_{\mathcal F}-X_{\mathcal F}^{(R)}$. By
		\eqref{eq:finite-coefficient} for both matrices,
		\[
		D:=b_{n,j,m}-b^{(R)}_{n,j,m}
		=(-1)^m\sum_{\substack{\mathcal F\ \mathrm{compatible}\\
				d(\mathcal F)=m,\ |\mathcal F|=m-j}}\Delta_{\mathcal F}.
		\]
		We show that the terms are orthogonal and each is small. As in Lemma~\ref{lem:family-orthogonality}, for distinct compatible
		families $\mathcal F,\mathcal G$ there is a matrix variable $\epsilon$
		used by just one of them, say $\mathcal F$. Given all other entries,
		$\Delta_{\mathcal G}$ is fixed. The pair $(\epsilon,\epsilon^{(R)})$
		is independent of those entries, and
		Lemma~\ref{lem:conditional-centering}, applied to both matrices, gives
		conditional mean zero for $X_{\mathcal F}$ and for
		$X_{\mathcal F}^{(R)}$. Hence $\E[\Delta_{\mathcal F}\Delta_{\mathcal G}]=0$.
		Thus $\E|D|^2=\sum_{\mathcal F}\E|\Delta_{\mathcal F}|^2$,
		summed over the families in $D$.
		
		For each such family, $n^{m/2}X_{\mathcal F}$ is a polynomial in
		the unscaled entries. In each monomial a diagonal entry occurs at
		most once, and each off-diagonal entry and its conjugate have total
		degree at most two. The $L^2$ convergence of the diagonal entries and the
		$L^4$ convergence of the off-diagonal entries give $L^2$ convergence
		of each monomial, by independence across entries. Hence
		$n^m\E|\Delta_{\mathcal F}|^2\to0$ as $R\to\infty$. This quantity
		depends only on $\mathcal F$ up to an order-preserving relabeling of
		its vertices, and for fixed $m$ there are finitely many such
		classes, so the convergence is uniform in $n$ and $\mathcal F$.
		Since there are $O_m(n^m)$ families,
		$\E|D|^2\leq C_m\sup_{\mathcal F}n^m\E|\Delta_{\mathcal F}|^2$,
		which tends to zero uniformly in $n$. This proves (a).
		
		(b) The $L^4$ convergence gives $\theta_R\to\theta$ and
		$\kappa_R\to\kappa$, the latter by
		\[
		\||\xi^{(R)}|^2-|\xi|^2\|_{L^2}
		\leq\|\xi^{(R)}-\xi\|_{L^4}
		(\|\xi^{(R)}\|_{L^4}+\|\xi\|_{L^4})\longrightarrow0.
		\]
		Since each $\rho_\ell$ is a polynomial in $\theta$ and
		$\overline\theta$, (b) follows.
	\end{proof}

	\subsection{Gaussian moments of the polymer sums.}\label{sec:power-sums}
	For this section, we will
	assume that $|\eta_i|,|\xi_{ij}|\leq R$ almost surely. 
	Recall the polymer sums
	\[ 	S_{\ell,n}:=\sum_{p\in\mathcal P_{\ell,n}}X_p \]
	from Section~\ref{sec:expansion}. Boundedness is needed
	here because the moments of the $S_{\ell,n}$ involve products of many
	polymer variables, and hence moments of the entries of all orders.
	
	We use the following counting principle. For a fixed number
	of factors, there are finitely many overlap patterns on $N$
	vertices, and each has $O(n^N)$ labelings. If each term is
	$O(n^{-D})$, the terms using $N$ vertices therefore contribute
	$O(n^{N-D})$ in total. The next lemma identifies the Gaussian pairings
	as the patterns attaining $N=D$ in the moments of the polymer sums.
	
	\begin{lemma}[Pairing patterns]\label{configuration}
		Assume $|\eta_i|,|\xi_{ij}|\leq R$ and fix $s,L\in\mathbb N$. For
		each $n$, let $p_1,\ldots,p_s$ be polymers of degrees at most $L$.
		Let $H$ be the graph formed by all edges of $p_1,\ldots,p_s$,
		together with the vertices $i$ of the singleton polymers $\{i\}$,
		and put
		\[
		N:=|V(H)|,\qquad D:=\frac12\sum_{j=1}^s d(p_j).
		\]
		Then
		\[
		\left|\E\prod_{j=1}^sX_{p_j}\right|\leq C_{s,L,R}n^{-D}.
		\]
		If the expectation is nonzero, then $N\leq D$, with equality
		possible only if every connected component of $H$ is the graph of a
		single polymer $p$, which occurs through exactly two factors $X_p$
		(so $s$ is even and $p_1,\ldots,p_s$ partition into $s/2$ pairs of
		identical polymers on pairwise disjoint vertex sets).
	\end{lemma}
	\begin{proof}
		Boundedness of the entries gives
		\[
		|X_{\{i\}}|\leq Rn^{-1/2},\qquad
		|X_e|\leq (R^2+1)n^{-1},\qquad
		|X_c|\leq 2R^{d(c)}n^{-d(c)/2}.
		\]
		Multiplying proves the bound. Suppose the expectation is nonzero.
		By Lemma~\ref{lem:conditional-centering}, every matrix variable used
		by one factor must occur in another factor. A factor containing a
		vertex $v$ uses a matrix variable at $v$, namely its diagonal entry or
		an incident edge, and every factor using that variable contains $v$.
		Thus each occupied vertex $v$ belongs to at least two factors. Writing $k_v$ for their
		number and counting vertex occurrences gives
		\[
		2D=\sum_{j=1}^s d(p_j)=\sum_{v\in V(H)}k_v\geq 2N.
		\]
		In the equality case exactly two factors meet at each vertex. They
		use the same incident matrix variables, since any variable used by
		only one of them would occur just once in the product. Along each
		edge the same two factors therefore meet at both endpoints. By
		connectedness this pair covers the entire component with identical
		edges. At an isolated vertex the two factors both use its diagonal
		variable. Thus each component consists of two copies of one polymer.
	\end{proof}
	
	\begin{proposition}[Gaussian moments of the polymer sums]\label{prop:polymer-sum-moments}
		For bounded entries and every fixed $L\geq1$,
		\[
		\E\prod_{\ell=1}^LS_{\ell,n}^{k_\ell}\longrightarrow
		\E\prod_{\ell=1}^LG_\ell^{k_\ell}
		\qquad\text{for all integers }k_1,\ldots,k_L\geq0,
		\]
		and $(S_{1,n},\ldots,S_{L,n})\Rightarrow(G_1,\ldots,G_L)$.
	\end{proposition}
	\begin{proof}
		For nonnegative integers $k_1,\ldots,k_L$, expand the mixed moment
		$\E\prod_{\ell=1}^L S_{\ell,n}^{k_\ell}$.
		By Lemma~\ref{configuration} and the counting principle, the leading
		terms pair identical polymers, with distinct pairs on disjoint vertex
		sets. The remaining terms contribute $o(1)$ in total. The two factors
		in each pair have the same degree, so the moment tends to
		zero if any $k_\ell$ is odd.
		
		Suppose $k_\ell=2r_\ell$ for every $\ell$. Each leading configuration
		determines its pairing uniquely, so we may fix a pairing of the
		positions within each degree; put $D=\sum_{\ell=1}^L\ell r_\ell$.
		There are $(n)_D/\prod_\ell(\ell!)^{r_\ell}$ ways to assign disjoint
		vertex sets to these pairs. By \eqref{eq:cycle-variance-block}, the
		polymer variances on each set of size $\ell$ sum to
		$\ell!\,\sigma_\ell^2/n^\ell$. Independence across the sets gives
		the total contribution of this pairing as
		\[
		\frac{(n)_D}{n^D}\prod_{\ell=1}^L(\sigma_\ell^2)^{r_\ell}.
		\]
		There are $(2r_\ell-1)!!$ pairings in each degree, and
		$\E G_\ell^{2r_\ell}=(2r_\ell-1)!!\,\sigma_\ell^{2r_\ell}$, while odd
		Gaussian moments vanish. Since $(n)_D/n^D\to1$, the limiting mixed
		moments are therefore exactly $\prod_{\ell=1}^L\E G_\ell^{k_\ell}$ for
		the independent Gaussian modes.
		A Gaussian vector is determined by its moments, so
		$(S_{1,n},\ldots,S_{L,n})\Rightarrow(G_1,\ldots,G_L)$.
	\end{proof}
	
	\subsection{Coefficients of the Wick exponential.}\label{sec:exponential-limit}
	We now prove Proposition~\ref{prop:coefficient-limit} by comparing
	$b_{n,j,m}$ with polynomials in the polymer sums. First, the second
	moments of the coefficients are explicit.
	
	\begin{lemma}[Coefficient second moments]\label{lem:coefficient-second-moment}
		\begin{enumerate}[label=(\alph*)]
			\item For $(\alpha,u),(\gamma,t)\in\mathcal V$,
			\begin{equation}\label{eq:wick-coefficient-kernel}
				\E\bigl[Q_\alpha(u)\overline{Q_\gamma(t)}\bigr]
				=\sum_{k\geq0}c_k(\alpha,\gamma)(u\overline t)^k.
			\end{equation}
			In particular, the coefficients of different $u$-degrees of $Q$
			are orthogonal.
			\item For $0\leq j\leq m\leq n$,
			\begin{equation}\label{eq:coefficient-second-moment}
				\E|b_{n,j,m}|^2=\frac{(n)_m}{n^m}\,\E|b_{j,m}|^2.
			\end{equation}
		\end{enumerate}
	\end{lemma}
	
	\begin{proof}
		(a) Put $Z:=\sum_{\ell\geq1}a_{\ell,\alpha}G_\ell u^\ell$ and
		$Z':=\sum_{\ell\geq1}a_{\ell,\gamma}G_\ell t^\ell$, so that
		$Q_\alpha(u)=e^{Z-\frac12\E Z^2}$ and $Q_\gamma(t)=e^{Z'-\frac12\E Z'^2}$.
		Since $Z+\overline{Z'}$ is a complex linear combination of the
		$G_\ell$, $\E e^{Z+\overline{Z'}}=e^{\frac12\E(Z+\overline{Z'})^2}$; the
		squares cancel against the normalizations, and
		$a_{\ell,\alpha}\overline{a_{\ell,\gamma}}=(\alpha\overline\gamma)^{\ell-1}$, so
		\[
		\E\bigl[Q_\alpha(u)\overline{Q_\gamma(t)}\bigr]
		=e^{\E Z\overline{Z'}}
		=\exp\Bigl(\sum_{\ell\geq1}\sigma_\ell^2(\alpha\overline\gamma)^{\ell-1}
		(u\overline t)^\ell\Bigr),
		\]
		which is (a) by the definition \eqref{eq:degree-second-moment} of
		$c_k(\alpha,\gamma)$.
		
		(b) By (a), the coefficients of $u^m$ in $Q_\alpha$ and $Q_\gamma$
		have covariance $c_m(\alpha,\gamma)$. For the corresponding
		coefficients $q_{n,\alpha,m}$ of $Q_{n,\alpha}$,
		Lemma~\ref{lem:orthogonality}(a) gives the same covariance times
		$(n)_m/n^m$. Take the coefficient of
		$\alpha^j\overline\gamma^{\,j}$ in both.
	\end{proof}
	
	To compare $b_{n,j,m}$ with $b_{j,m}$, we substitute the polymer sums
	for the Gaussian modes. For $0\leq j\leq m$, define the polynomial
	\[
	b_{j,m}(x):=[\alpha^j u^m]\exp\Bigl(\sum_{\ell=1}^m
	a_{\ell,\alpha}x_\ell u^\ell
	-\frac12\sum_{\ell=1}^m a_{\ell,\alpha}^2\sigma_\ell^2u^{2\ell}\Bigr),
	\]
	so that $b_{j,m}(G)=b_{j,m}$ for $G:=(G_\ell)_{\ell\geq1}$, and put
	$S_n:=(S_{\ell,n})_{\ell\geq1}$.
	
	\begin{lemma}[Orthogonal decomposition]\label{lem:orthogonal-decomposition}
		Assume $|\eta_i|,|\xi_{ij}|\leq R$ almost surely, and let
		$0\leq j\leq m\leq n$. Then
		$\E\bigl[(b_{j,m}(S_n)-b_{n,j,m})\,b_{n,j,m}\bigr]=0$. Consequently,
		by Pythagoras,
		\begin{equation}\label{eq:wick-coefficient-projection}
			\E|b_{j,m}(S_n)-b_{n,j,m}|^2
			=\E|b_{j,m}(S_n)|^2-\E|b_{n,j,m}|^2.
		\end{equation}
	\end{lemma}
	
	\begin{proof}
		\emph{(i) Decomposition.} Recall from \eqref{eq:finite-coefficient} that
		\[
		b_{n,j,m}=(-1)^m\sum_{\substack{\mathcal F\ \mathrm{compatible}\\
				d(\mathcal F)=m,\ |\mathcal F|=m-j}}X_{\mathcal F}.
		\]
		The exponential defining $b_{j,m}(S_n)$ factors as
		\[
		\exp\Bigl(\sum_{\ell=1}^ma_{\ell,\alpha}S_{\ell,n}u^\ell\Bigr)\cdot
		\exp\Bigl(-\frac12\sum_{\ell=1}^ma_{\ell,\alpha}^2\sigma_\ell^2u^{2\ell}\Bigr),
		\]
		and we call the factors from the second exponential \emph{Wick
			factors}. Expanding $S_{\ell,n}=\sum_{p\in\mathcal P_{\ell,n}}X_p$ and
		using $\prod_ia_{d(p_i),\alpha}=(-1)^{\sum_id(p_i)}\alpha^{\sum_id(p_i)-k}$,
		the first exponential has, for each $m'\geq0$,
		\[
		[u^{m'}]\exp\Bigl(\sum_\ell a_{\ell,\alpha}S_{\ell,n}u^\ell\Bigr)
		=\sum_{k\geq0}\frac1{k!}
		\sum_{\substack{(p_1,\ldots,p_k)\\ \sum_id(p_i)=m'}}
		(-1)^{m'}\alpha^{m'-k}\,X_{p_1}\cdots X_{p_k},
		\]
		the inner sum over ordered $k$-tuples of polymers. Wick factors
		contribute $u$-degree but no polymers, so every term contributing to
		$[\alpha^ju^m]$ is a constant times
		$X_{p_1}\cdots X_{p_k}$, which we say \emph{occupies} the vertex set
		$\bigcup_iV(p_i)$, with $\sum_id(p_i)\leq m$, strictly if a Wick
		factor is used. Hence it occupies
		$\bigl|\bigcup_iV(p_i)\bigr|\leq\sum_id(p_i)\leq m$ vertices, with
		equality exactly when no Wick factor is used and $p_1,\ldots,p_k$ are
		pairwise vertex-disjoint; then $k=m-j$, and these tuples are the $k!$
		orderings of the compatible families in the display above, which
		cancel the factor $1/k!$. So the terms occupying $m$ vertices give
		exactly $b_{n,j,m}$, and $T_n:=b_{j,m}(S_n)-b_{n,j,m}$ is the sum of
		the terms occupying fewer than $m$ vertices.
		
		\emph{(ii) Orthogonality.} We show that, for polymers $p_1,\ldots,p_k$
		and a compatible family $\mathcal F$,
		\begin{equation}\label{eq:missed-vertex-orthogonality}
			\E\bigl[X_{p_1}\cdots X_{p_k}\,X_{\mathcal F}\bigr]=0
			\qquad\text{whenever }V(\mathcal F)\not\subset\textstyle\bigcup_iV(p_i).
		\end{equation}
		This applies to every term of $T_n$ and every family $\mathcal F$ in
		$b_{n,j,m}$, since $|V(\mathcal F)|=m$ while the term occupies fewer
		than $m$ vertices. To prove \eqref{eq:missed-vertex-orthogonality},
		pick $v\in V(\mathcal F)\setminus\bigcup_iV(p_i)$, and let $\epsilon$ be a matrix variable at $v$ used
		by $\mathcal F$, namely the diagonal entry or an edge of the polymer
		of $\mathcal F$ containing $v$. Every polymer using $\epsilon$
		contains $v$, so no $p_i$ uses $\epsilon$, and exactly one polymer of
		$\mathcal F$ does. Conditioning on all other matrix variables,
		Lemma~\ref{lem:conditional-centering} gives
		\eqref{eq:missed-vertex-orthogonality}. Summing over terms and families
		gives $\E[T_nb_{n,j,m}]=0$, hence
		$\E|b_{j,m}(S_n)|^2=\E|b_{n,j,m}|^2+\E|T_n|^2$, which is
		\eqref{eq:wick-coefficient-projection}.
	\end{proof}
	
	\begin{proof}[Proof of Proposition~\ref{prop:coefficient-limit}]
		\emph{Step 1: bounded entries.} Assume
		$|\eta_i|,|\xi_{ij}|\leq R$ almost surely, and fix $j,m$. By
		Lemma~\ref{lem:orthogonal-decomposition}, it suffices to show that the
		two second moments on the right of
		\eqref{eq:wick-coefficient-projection} have the same limit. Both tend
		to $\E|b_{j,m}|^2$: the first by Proposition~\ref{prop:polymer-sum-moments},
		since $b_{j,m}$ is a polynomial, and the second by
		\eqref{eq:coefficient-second-moment}, since $(n)_m/n^m\to1$. Hence
		$b_{j,m}(S_n)-b_{n,j,m}\to0$ in $L^2$. For $m\leq L$, the
		$b_{j,m}(S_n)$ are polynomials in $(S_{1,n},\ldots,S_{L,n})$, so by
		Proposition~\ref{prop:polymer-sum-moments} they converge jointly in law
		to $b_{j,m}(G)=b_{j,m}$, and by the $L^2$ estimate so do the
		$b_{n,j,m}$.
		
		Summing \eqref{eq:coefficient-second-moment} over
		$0\leq j\leq m\leq L$ gives
		$\E\|b_n^{(\leq L)}\|^2\to\E\|b^{(\leq L)}\|^2$, so
		$W_2(b_n^{(\leq L)},b^{(\leq L)})\to0$ by the characterization of
		$W_2$-convergence in Section~\ref{sec:notation}.
		
		\emph{Step 2: general entries.} We remove the boundedness by
		truncating at level $R$ (Section~\ref{sec:entry-cutoff}). Fix $L\geq1$, and let $b_n^{(R,\leq L)}$ be the vector
		$(b^{(R)}_{n,j,m})_{0\leq j\leq m\leq L}$. Using independent standard Gaussians
		$N_\ell$, couple $G_\ell=\sigma_\ell N_\ell$ and
		$G_\ell^{(R)}=\sigma_\ell^{(R)}N_\ell$, and let $b^{(R,\leq L)}$ be
		$b^{(\leq L)}$ with $G_\ell$ and $\sigma_\ell^2$ replaced by
		$G_\ell^{(R)}$ and $(\sigma_\ell^{(R)})^2$. Then
		\[
		W_2\bigl(b_n^{(\leq L)},b^{(\leq L)}\bigr)
		\leq\bigl\|b_n^{(\leq L)}-b_n^{(R,\leq L)}\bigr\|_{L^2}
		+W_2\bigl(b_n^{(R,\leq L)},b^{(R,\leq L)}\bigr)
		+\bigl\|b^{(R,\leq L)}-b^{(\leq L)}\bigr\|_{L^2}.
		\]
		The first term tends to zero as
		$R\to\infty$, uniformly in $n$, by Lemma~\ref{lem:truncation}(a). For
		each large integer $R$, the second tends to zero as $n\to\infty$, by
		Step~1 applied to the truncated entries. The third tends to zero as
		$R\to\infty$: the entries of $b^{(R,\leq L)}$ are polynomials in
		finitely many $N_\ell$ whose coefficients converge by
		Lemma~\ref{lem:truncation}(b), and Gaussians have moments of all
		orders. Letting first $n\to\infty$ and then $R\to\infty$ proves
		Proposition~\ref{prop:coefficient-limit}.
	\end{proof}
	
	\begin{corollary}[Polymer sums]\label{cor:polymer-sums-general}
		For every $\ell\geq1$, $S_{\ell,n}=(-1)^\ell b_{n,\ell-1,\ell}$ and
		$G_\ell=(-1)^\ell b_{\ell-1,\ell}$. Hence
		$(S_{\ell,n})_{\ell\leq L}\Rightarrow(G_\ell)_{\ell\leq L}$ jointly in
		law for every $L$, under the assumptions of Section~\ref{sec:intro}.
	\end{corollary}
	
	\begin{proof}
		By \eqref{eq:finite-coefficient} with $m=\ell$ and $j=\ell-1$, the
		families in $b_{n,\ell-1,\ell}$ are single polymers of degree $\ell$,
		so $b_{n,\ell-1,\ell}=(-1)^\ell S_{\ell,n}$. In $Q_\alpha$, a product of
		$k$ terms $a_{\ell_i,\alpha}G_{\ell_i}u^{\ell_i}$ has $u$-degree
		exceeding its $\alpha$-degree by $k$, and each factor from the second
		sum adds $2$, so only $a_{\ell,\alpha}G_\ell u^\ell$ contributes to
		$\alpha^{\ell-1}u^\ell$ and $b_{\ell-1,\ell}=(-1)^\ell G_\ell$.
		Proposition~\ref{prop:coefficient-limit} gives the joint limit.
	\end{proof}
	
	Boundedness of the entries is thus needed only for the convergence of
	higher moments in Proposition~\ref{prop:polymer-sum-moments}.
	
	\subsection{Functional and Wasserstein convergence.}\label{sec:assembly}
	\label{sec:quadratic-wasserstein}\label{sec:inner-products}
	Write $Q_n(\alpha,u):=Q_{n,\alpha}(u)$ and $Q(\alpha,u):=Q_\alpha(u)$,
	and let $Q_n^{(\leq L)}$ and $Q^{(\leq L)}$ be their Taylor polynomials
	through degree $L$ in $u$. Recall that
	Proposition~\ref{prop:wasserstein-convergence} asserts
	$W_{2,\mathcal K}(F_n,\Psi)\to0$ for every compact $\mathcal K\Subset\mathcal D$,
	together with convergence of $\E\|F_n\|_{\mathcal K}^2$. We first prove the
	analogous statement for $Q_n$ and $Q$, and then transfer it to $F_n$
	by the polymer approximation.
	
	\begin{proof}[Proof of Proposition~\ref{prop:wasserstein-convergence}]
		\emph{Step 1: $Q_n\to Q$.} Fix $L$ and a compact $K\Subset\mathcal V$. For both $Q_n$ and $Q$,
		the coefficients of different degrees are orthogonal, by
		Lemmas~\ref{lem:orthogonality}(a)
		and~\ref{lem:coefficient-second-moment}(a). Since $(n)_m/n^m\leq1$,
		the pointwise $L^2$ tail beyond degree $L$ of either one is bounded
		by $\sum_{m>L}c_m(\alpha)|u|^{2m}$, which is at most $Cb^{-L}$ on
		compact sets by Lemma~\ref{lem:degree-tails}. Applying this on the
		enlargement of $K$ from Lemma~\ref{lem:sup-norm} gives, for some
		$b_K>1$,
		\begin{equation}\label{eq:polymer-taylor-tail}
			\sup_n\E\|Q_n-Q_n^{(\leq L)}\|_K^2
			+\E\|Q-Q^{(\leq L)}\|_K^2\leq C_K b_K^{-L}.
		\end{equation}
		By Proposition~\ref{prop:coefficient-limit},
		$W_2(b_n^{(\leq L)},b^{(\leq L)})\to0$. The map taking these
		coefficients to their Taylor polynomial is linear and bounded into
		$C(K)$, so $W_{2,K}(Q_n^{(\leq L)},Q^{(\leq L)})\to0$. By the triangle
		inequality and
		\eqref{eq:polymer-taylor-tail},
		\[
		W_{2,K}(Q_n,Q)\leq 2\bigl(C_Kb_K^{-L}\bigr)^{1/2}
		+W_{2,K}\bigl(Q_n^{(\leq L)},Q^{(\leq L)}\bigr),
		\]
		and letting first $n\to\infty$ and then $L\to\infty$ gives
		$W_{2,K}(Q_n,Q)\to0$.
		
		\emph{Step 2: $F_n\to\Psi$.} For $\mathcal K\Subset\mathcal D$, composition with
		$\phi(\alpha,z):=(\alpha,g_\alpha(z))$ maps $C(\phi(\mathcal K))$ to
		$C(\mathcal K)$ with norm at most one, and $Q\circ\phi=\Psi$, so
		$W_{2,\mathcal K}(Q_n\circ\phi,\Psi)\to0$. Proposition~\ref{prop:joint-polymer-approximation}
		gives $\E\|F_n-Q_n\circ\phi\|_{\mathcal K}^2\to0$, proving
		\eqref{eq:compact-wasserstein-convergence} by the triangle inequality.
		
		The moment conclusions follow from the characterization of
		$W_2$-convergence in Section~\ref{sec:notation}.
	\end{proof}
	
	Convergence in law in $\mathcal H(\mathcal D)$ follows from
	Proposition~\ref{prop:wasserstein-convergence} by a general principle.
	
	\begin{lemma}[From compact $W_2$ to weak convergence]\label{lem:w2-weak}
		Let $U\subset\C^2$ be open, and let $X_n$ and $X$ be random elements of
		$\mathcal H(U)$ such that $W_{2,K}(X_n,X)\to0$ for every compact
		$K\Subset U$. Then
		$X_n\Rightarrow X$ in $\mathcal H(U)$.
	\end{lemma}
	
	\begin{proof}
		By the portmanteau theorem \cite[Theorem~2.1]{Billingsley1999}, it
		suffices to show $\limsup_n\mathbb P(X_n\in A)\leq\mathbb P(X\in A)$ for
		closed $A\subset\mathcal H(U)$. Fix compact sets $K_j\subset U$ with
		$K_j\subset\operatorname{int}K_{j+1}$ and $\bigcup_jK_j=U$, write
		$r_jf:=f|_{K_j}$, and let $A_j$ be the closure of $r_j(A)$ in $C(K_j)$.
		
		\emph{Step 1: fixed $j$.} Convergence in $W_{2,K_j}$ gives $r_jX_n\Rightarrow r_jX$
		in $C(K_j)$ (Section~\ref{sec:notation}). Since $X_n\in A$ implies
		$r_jX_n\in A_j$, and $A_j$ is closed,
		\[
		\limsup_{n\to\infty}\mathbb P(X_n\in A)
		\leq\limsup_{n\to\infty}\mathbb P(r_jX_n\in A_j)
		\leq\mathbb P(r_jX\in A_j).
		\]
		
		\emph{Step 2: letting $j\to\infty$.} The events $\{r_jX\in A_j\}$ decrease, since
		restriction from $K_{j+1}$ to $K_j$ is continuous and maps $A_{j+1}$
		into $A_j$. Their intersection is $\{X\in A\}$: if $r_jf\in A_j$ for
		all $j$, choose $f_j\in A$ with $\|f_j-f\|_{K_j}<1/j$; every compact
		subset of $U$ lies in some $K_j$, so $f_j\to f$ in $\mathcal H(U)$ and
		$f\in A$. Letting $j\to\infty$ proves the claim.
	\end{proof}
	
	\begin{proof}[Proof of Theorem~\ref{gaussian}]
		Apply Lemma~\ref{lem:w2-weak} with $U=\mathcal D$, $X_n=F_n$ and
		$X=\Psi$, using Proposition~\ref{prop:wasserstein-convergence}. This
		gives $F_n\Rightarrow\Psi$ in $\mathcal H(\mathcal D)$, and restricting
		to the slices $\{\alpha\}\times\Omega_\alpha$ with $|\alpha|\leq1$
		gives \eqref{eq:functional-limit}.
	\end{proof}

	\section{Fluctuations of the zeros}\label{sec:further-consequences}
	The functional limit describes the fluctuations of the zeros through
	two operations: logarithmic differentiation gives the fluctuations of
	the Cauchy transform, and contour integration gives analytic statistics
	of the zeros. This section can be read directly after
	Sections~\ref{sec:functional-limits} and~\ref{sec:consequences}.
	
	\subsection{Fluctuations of Cauchy transforms.}
	Recall that $g_{n,\alpha}$ and $h_{n,\alpha}$ are the Cauchy
	transforms of the empirical zero measures $\nu_{n,\alpha}$ of
	$p_{n,\alpha}$ and $\mu_{n,\alpha}$ of $P_{n,\alpha}$. Logarithmic
	differentiation of $p_{n,\alpha}=P_{n,\alpha}F_{n,\alpha}$, as in
	\eqref{eq:empirical-transform-decomposition}, splits the
	Cauchy-transform fluctuations into a deterministic centering bias and
	a random term:
	\begin{equation}\label{eq:transform-fluctuation-split}
		n\bigl(g_{n,\alpha}(z)-g_\alpha(z)\bigr)
		=n\bigl(h_{n,\alpha}(z)-g_\alpha(z)\bigr)
		+\frac{\partial_z F_{n,\alpha}(z)}{F_{n,\alpha}(z)}.
	\end{equation}
	The second term converges to $H_\alpha'(z)$ by Theorem~\ref{gaussian}.
	The first term has a finite limit coming from the Riccati equation for
	$h_{n,\alpha}$.
	
	\begin{lemma}[Deterministic centering correction]\label{correction}
		Locally uniformly in $(\alpha,z)\in\mathcal U$,
		\begin{equation}\label{eq:annealed-transform-correction}
			n\bigl(h_{n,\alpha}(z)-g_\alpha(z)\bigr)\longrightarrow g^{(1)}_\alpha(z),
		\end{equation}
		where
		\begin{equation}\label{eq:zero-transform-bias}
			g^{(1)}_\alpha(z)
			:=-\frac{\alpha g_\alpha'(z)}{z+2\alpha g_\alpha(z)}
			=\frac{\alpha g_\alpha(z)^3}
			{(1+\alpha g_\alpha(z)^2)^2},
			\qquad g^{(1)}_0:=0.
		\end{equation}
	\end{lemma}
	
	\begin{proof}
		Write $h=h_{n,\alpha}$ and $g=g_\alpha$. Subtracting the quadratic
		$\alpha g^2+zg-1=0$ \eqref{eq:g-quadratic} from the Riccati equation
		$\frac\alpha nh'+\alpha h^2+zh-1=0$ \eqref{eq:hermite-riccati} gives
		\begin{equation}\label{eq:hermite-riccati-difference}
			n(h-g)\bigl(z+\alpha(h+g)\bigr)=-\alpha h'.
		\end{equation}
		From \eqref{eq:g-quadratic} and its derivative we have
		\[
		z+2\alpha g_\alpha=\frac{1+\alpha g_\alpha^2}{g_\alpha}=\sqrt{z^2+4\alpha},
		\qquad
		g_\alpha'=-\frac{g_\alpha^2}{1+\alpha g_\alpha^2},
		\]
		with the branch of $\sqrt{z^2+4\alpha}$ asymptotic to $z$ at infinity;
		these give the second expression in \eqref{eq:zero-transform-bias}. By
		\eqref{eq:hermite-transform-limit} in
		Lemma~\ref{lem:hermite-zero-law}, $h_{n,\alpha}\to g_\alpha$
		locally uniformly on $\mathcal U$, also after differentiation in $z$,
		and $\sqrt{z^2+4\alpha}$ is bounded away from zero on compact subsets
		of $\mathcal U$. Dividing \eqref{eq:hermite-riccati-difference} by
		$z+\alpha(h+g)$ and letting $n\to\infty$ gives
		\eqref{eq:annealed-transform-correction}.
	\end{proof}
	
	\begin{theorem}[Joint Cauchy-transform fluctuations]
		\label{thm:cauchy-transform-fluctuations}
		On every compact $\mathcal K\Subset\mathcal D$, the Cauchy-transform
		fluctuations satisfy
		\begin{equation}\label{eq:empirical-transform-clt}
			n\bigl(g_{n,\alpha}(z)-g_\alpha(z)\bigr)
			\Rightarrow g^{(1)}_\alpha(z)+H_\alpha'(z)
			\qquad\text{in }C(\mathcal K),
		\end{equation}
		and in particular
		\begin{equation}\label{eq:empirical-transform-rate}
			\sup_{(\alpha,z)\in\mathcal K}|g_{n,\alpha}(z)-g_\alpha(z)|
			=O_{\mathbb P}(n^{-1}),
		\end{equation}
		with the convention that the left side of
		\eqref{eq:empirical-transform-clt} is $0$ and the supremum in
		\eqref{eq:empirical-transform-rate} is $+\infty$ on the event (of
		vanishing probability) that $p_{n,\alpha}$ vanishes on $\mathcal K$.
	\end{theorem}
	
	\begin{proof}
		By Theorem~\ref{gaussian}, $F_n\Rightarrow\Psi$ in
		$\mathcal H(\mathcal D)$. Since $\min_{\mathcal K}|\Psi|>0$ almost
		surely, the probability that $F_n$ has a zero on $\mathcal K$ tends
		to zero. At any function $F$ nonvanishing on $\mathcal K$, the map
		$F\mapsto(\partial_zF/F)|_{\mathcal K}$ is continuous into
		$C(\mathcal K)$ by Cauchy's formula. Hence
		$\partial_z F_n/F_n\Rightarrow H_\alpha'$ in $C(\mathcal K)$,
		with the quotient set to zero on the exceptional event.
		Lemma~\ref{correction} gives $n(h_{n,\alpha}-g_\alpha)\to g^{(1)}_\alpha$
		uniformly on $\mathcal K$. Combining these limits in
		\eqref{eq:transform-fluctuation-split} proves
		\eqref{eq:empirical-transform-clt}, and tightness in
		$C(\mathcal K)$ gives \eqref{eq:empirical-transform-rate}.
	\end{proof}
	
	\subsection{Analytic linear statistics.}
	Integrating the Cauchy-transform limit against an analytic test
	function describes sums of that function over the zeros. Subtracting
	the leading semicircle contribution leaves fluctuations of order one.
	Let $U\subset\C$ be an open neighborhood of $\Sigma_\alpha$.
	For $f\in\mathcal H(U)$, put
	\[
	\mathcal L_{n,\alpha}(f)
	:=\sum_{j=1}^n f(\zeta_{n,j}^{(\alpha)})-n\int f\,d\mu_\alpha,
	\]
	and set $\mathcal L_{n,\alpha}(f):=0$ if some zero lies outside
	$U$, an event whose probability tends to zero. A positively
	oriented piecewise smooth Jordan curve $\Gamma=\partial D$ is
	\emph{admissible} if $\Sigma_\alpha\subset D$ and
	$\overline D\Subset U$.
	
	The correction $g^{(1)}_\alpha$ produces a deterministic bias in these
	statistics, defined for $f$ analytic near $\Sigma_\alpha$ by
	\begin{equation}\label{eq:hermite-bias}
		B_\alpha(f)
		:=\frac12\int_{-2}^2f(\beta_\alpha t)
		\frac{dt}{\pi\sqrt{4-t^2}}
		-\frac14\bigl(f(2\beta_\alpha)+f(-2\beta_\alpha)\bigr).
	\end{equation}
	
	\begin{lemma}[Hermite bias]\label{lem:hermite-bias}
		Fix $\alpha\in\C$, $f\in\mathcal H(U)$ and an admissible contour
		$\Gamma$. Then
		\begin{equation}\label{eq:bias-contour-identity}
			\frac1{2\pi i}\int_\Gamma f(z)g^{(1)}_\alpha(z)\,dz=B_\alpha(f),
		\end{equation}
		and consequently
		\[
		\sum_{j=1}^n f(\beta_\alpha t_{n,j})-n\int f\,d\mu_\alpha
		\longrightarrow B_\alpha(f).
		\]
	\end{lemma}
	
	\begin{proof}
		With $\sqrt{z^2+4\alpha}=z+2\alpha g_\alpha$, as in the proof of
		Lemma~\ref{correction},
		$1+\alpha g_\alpha^2=g_\alpha\sqrt{z^2+4\alpha}$
		and $\alpha g_\alpha=(\sqrt{z^2+4\alpha}-z)/2$, so
		\[
		g^{(1)}_\alpha(z)
		=\frac{\alpha g_\alpha^3}{(1+\alpha g_\alpha^2)^2}
		=\frac{\alpha g_\alpha}{z^2+4\alpha}
		=\frac1{2\sqrt{z^2+4\alpha}}-\frac{z}{2(z^2+4\alpha)}.
		\]
		Using the arcsine Cauchy transform and a partial fraction,
		\[
		\int_{-2}^2\frac{dt}{(w-t)\,\pi\sqrt{4-t^2}}=\frac1{\sqrt{w^2-4}},
		\qquad
		\frac14\Bigl(\frac1{z-2\beta_\alpha}+\frac1{z+2\beta_\alpha}\Bigr)
		=\frac{z}{2(z^2+4\alpha)},
		\]
		the first at $w=z/\beta_\alpha$ for $\alpha\neq0$ (and trivially at
		$\alpha=0$), we recognize $g^{(1)}_\alpha$ as the
		Cauchy transform of $(\beta_\alpha\id)_\#\mu^{(1)}$, where
		\[
		\mu^{(1)}(dt)
		=\frac{\mathbf1_{(-2,2)}(t)}{2\pi\sqrt{4-t^2}}\,dt
		-\frac14\delta_{-2}(dt)-\frac14\delta_2(dt).
		\]
		Thus, compared with $n$ copies of the semicircle law, the points $t_{n,j}$ carry a deficit of $\frac14$ at each edge $\pm2$, balanced by
		half an arcsine law in the bulk; $\mu^{(1)}$ has total mass zero.
		Cauchy's formula gives the left side of
		\eqref{eq:bias-contour-identity} as
		$\int f\,d((\beta_\alpha\id)_\#\mu^{(1)})$, which is $B_\alpha(f)$ by
		\eqref{eq:hermite-bias}. Finally, since the zeros of $P_{n,\alpha}$ lie
		on $\Sigma_\alpha$, the residue theorem gives
		\[
		\frac1{2\pi i}\int_\Gamma f\,n(h_{n,\alpha}-g_\alpha)\,dz
		=\sum_jf(\beta_\alpha t_{n,j})-n\int f\,d\mu_\alpha,
		\]
		and
		Lemma~\ref{correction} on $\Gamma$ proves the last claim.
	\end{proof}
	
	\begin{remark}[The Hermite bias as a frozen log-gas correction]\label{rem:hermite-bias-beta}
		For the Gaussian $\beta$-ensembles, normalized so that the semicircle
		law lives on $[-2,2]$, the mean of a linear statistic satisfies
		$\E\sum_jf(\lambda_j)-n\int f\,d\mu_\scal\to(2/\beta-1)\bigl[\tfrac14(f(2)+f(-2))
		-\tfrac12\int f\,d\mu_{\mathrm{arc}}\bigr]$, with $\mu_{\mathrm{arc}}$
		the arcsine law \cite{Johansson1998},
		\cite[Theorem~1.2 and Section~2.1]{DumitriuEdelman2006}. As
		$\beta\to\infty$ the eigenvalues freeze at the zeros of the Hermite
		polynomial \cite{DumitriuEdelman2005}, and the factor becomes $-1$:
		this is $B_{-1}(f)$.
	\end{remark}

	\begin{theorem}[Analytic zero statistics]\label{clt}
		Let $|\alpha|\leq1$ and $f_1,\ldots,f_r\in\mathcal H(U)$. Then
		\[
		(\mathcal L_{n,\alpha}(f_1),\ldots,\mathcal L_{n,\alpha}(f_r))
		\Rightarrow(\mathcal L_\alpha(f_1),\ldots,\mathcal L_\alpha(f_r))
		\qquad\text{in }\C^r,
		\]
		where the limit is the random vector, with jointly Gaussian real and
		imaginary parts, given by
		\[
		\mathcal L_\alpha(f)
		=B_\alpha(f)+\frac1{2\pi i}\int_\Gamma f(z)H_\alpha'(z)\,dz
		\]
		for any admissible contour, with $B_\alpha$ as in
		\eqref{eq:hermite-bias}.
		The convergence holds jointly for finitely many parameters
		$|\alpha|\leq1$, each with its own $U$ and test functions.
	\end{theorem}
	
	\begin{proof}
		Fix an admissible contour $\Gamma=\partial D$. With probability
		tending to one, every zero of $p_{n,\alpha}$ lies in $D$, by
		Theorem~\ref{zerofree}. On this event, the residue theorem gives
		\[
		\mathcal L_{n,\alpha}(f)
		=\frac1{2\pi i}\int_\Gamma
		f(z)n(g_{n,\alpha}(z)-g_\alpha(z))\,dz.
		\]
		Apply Theorem~\ref{thm:cauchy-transform-fluctuations} on
		$\{\alpha\}\times\Gamma$ and use continuity of the contour integrals
		for $f_1,\ldots,f_r$. For several parameters, take the union of the
		corresponding sets $\{\alpha_j\}\times\Gamma_j$.
		By \eqref{eq:bias-contour-identity}, the integral of
		$fg^{(1)}_\alpha$ is $B_\alpha(f)$, giving the asserted limit.
	\end{proof}
	
	\subsection{Dickson expansion and the explicit limit.}\label{sec:dickson-statistics}\label{sec:dickson-explicit-limit}
	For $|\alpha|\leq1$, the contour formula of Theorem~\ref{clt} can be
	evaluated explicitly. In the Dickson basis below, each positive degree
	isolates one Gaussian mode. Recall from
	Section~\ref{sec:disk-coordinate} that $1/g_\alpha$ and
	$\tau_\alpha=-\alpha g_\alpha$ are the two roots of
	\begin{equation}\label{eq:recap-roots}
		X^2-zX-\alpha
		=\Bigl(X-\frac1{g_\alpha(z)}\Bigr)\bigl(X-\tau_\alpha(z)\bigr),
		\qquad z\in\Omega_\alpha,
	\end{equation}
	with $|\tau_\alpha|<\sqrt{|\alpha|}$ for $\alpha\neq0$, and that
	$J_\alpha(w)=w-\alpha/w$ inverts $\tau_\alpha$ on $\{|w|<\sqrt{|\alpha|}\}$
	(and $1/g_\alpha$ on $\{|w|>\sqrt{|\alpha|}\}$).
	
	\emph{Dickson polynomials and coefficients.} For $\ell\geq0$, let
	\begin{equation}\label{eq:dickson-roots}
		\mathcal T_{\ell,\alpha}(x):=x_+^\ell+x_-^\ell,
		\qquad x_\pm\text{ the roots of }X^2-xX-\alpha.
	\end{equation}
	These are the Dickson polynomials $D_k(x,-\alpha)$
	\cite{LidlMullenTurnwald1993}. As power sums of the roots of
	$X^2-xX-\alpha$, they satisfy
	\begin{equation}\label{eq:dickson-recurrence}
		\mathcal T_{0,\alpha}=2,\qquad \mathcal T_{1,\alpha}=x,\qquad
		\mathcal T_{\ell,\alpha}
		=x\mathcal T_{\ell-1,\alpha}+\alpha\mathcal T_{\ell-2,\alpha}
		\quad(\ell\geq2).
	\end{equation}
	Thus $\mathcal T_{\ell,\alpha}$ is monic of degree $\ell$ in $x$, with
	coefficients polynomial in $\alpha$, and $\mathcal T_{\ell,0}(x)=x^\ell$
	for $\ell\geq1$. For $\alpha\neq0$ and $x=2\beta_\alpha\cos\vartheta$, the roots in
	\eqref{eq:dickson-roots} are $\beta_\alpha e^{\pm i\vartheta}$, so
	$\mathcal T_{\ell,\alpha}(x)=2\beta_\alpha^\ell T_\ell(x/(2\beta_\alpha))$,
	where $T_\ell(\cos\vartheta)=\cos\ell\vartheta$ are the Chebyshev
	polynomials; at $\alpha=-1$, $\mathcal T_{\ell,-1}(x)=2T_\ell(x/2)$. At $x=z\in\Omega_\alpha$, the roots in
	\eqref{eq:dickson-roots} are $1/g_\alpha(z)$ and
	$\tau_\alpha(z)=-\alpha g_\alpha(z)$ by \eqref{eq:recap-roots}, so
	\begin{equation}\label{eq:dickson-key-identity}
		\mathcal T_{\ell,\alpha}(z)
		=g_\alpha(z)^{-\ell}+\bigl(-\alpha g_\alpha(z)\bigr)^\ell,
		\qquad z\in\Omega_\alpha,\ \ell\geq1.
	\end{equation}
	
	\emph{Expansions in the Dickson basis.} For $f$ analytic near
	$\Sigma_\alpha$, we expand
	\begin{equation}\label{eq:analytic-chebyshev-expansion}
		f(z)=f_{0,\alpha}
		+\sum_{\ell\geq1}f_{\ell,\alpha}\mathcal T_{\ell,\alpha}(z),
	\end{equation}
	where for $\alpha\neq0$ the coefficients are given by
	\[
	f_{\ell,\alpha}
	=\frac1{\pi\beta_\alpha^\ell}
	\int_0^\pi f(2\beta_\alpha\cos\vartheta)\cos(\ell\vartheta)\,d\vartheta,
	\qquad \ell\geq0,
	\]
	while at $\alpha=0$ the expansion is the Taylor series at $0$, with
	$f_{\ell,0}=f^{(\ell)}(0)/\ell!$.
	
	\begin{lemma}[Dickson expansion and coefficient extraction]\label{lem:dickson-biorthogonality}
		Let $\alpha\in\C$ and let $f$ be analytic on a neighborhood $U$ of
		$\Sigma_\alpha$.
		\begin{enumerate}[label=(\alph*)]
			\item The expansion \eqref{eq:analytic-chebyshev-expansion} converges
			locally uniformly on a neighborhood of $\Sigma_\alpha$ (and on
			every confocal ellipse whose closed interior lies in $U$), with
			$|\beta_\alpha^\ell f_{\ell,\alpha}|$ decaying geometrically.
			\item For every admissible contour $\Gamma$,
			\begin{equation}\label{eq:chebyshev-contour-extraction}
				\frac1{2\pi i}\int_\Gamma f'(z)g_\alpha(z)^m\,dz
				=mf_{m,\alpha},\qquad m\geq1.
			\end{equation}
		\end{enumerate}
	\end{lemma}
	
	\begin{proof}
		At $\alpha=0$, $\mathcal T_{m,0}(z)=z^m$ for $m\geq1$ and
		$g_0(z)=1/z$, so the claims follow from the Taylor expansion of $f$
		and Cauchy's coefficient formula for $f'$. Now suppose $\alpha\ne0$
		and put $\beta=\beta_\alpha$. For $z=J_\alpha(w)=w-\alpha/w$, the
		two roots in \eqref{eq:dickson-roots} are $w$ and $-\alpha/w$.
		
		Choose $R>|\beta|$ so that the ellipse $J_\alpha(\{|w|=R\})$ and
		its interior lie in $U$. Their preimage is the closed annulus
		$|\alpha|/R\leq|w|\leq R$, on a neighborhood of which
		$F(w)=f(J_\alpha(w))$ has a Laurent expansion $\sum_{k\in\mathbb Z}a_kw^k$.
		Since $F(w)=F(-\alpha/w)$, uniqueness gives $a_{-m}=(-\alpha)^m a_m$.
		Pairing the terms of indices $\pm m$ gives $a_m\mathcal T_{m,\alpha}(z)$.
		The Laurent series converges uniformly on the closed annulus, so the
		paired series converges uniformly on the filled ellipse. The coefficient
		formula on $w=\beta e^{i\vartheta}$ gives the cosine integral defining
		$f_{m,\alpha}$, hence $a_m=f_{m,\alpha}$ for $m\geq0$. On $|w|=R$
		the same formula gives $|\beta^m a_m|\leq\max_{|w|=R}|F(w)|\,(|\beta|/R)^m$.
		This proves (a).
		
		Given an admissible $\Gamma$, choose such an ellipse inside it and
		deform $\Gamma$ to the ellipse by Cauchy's theorem. On $|w|=R$ we
		have $g_\alpha(J_\alpha(w))=1/w$, so for $m\geq1$,
		\[
		\frac1{2\pi i}\int_\Gamma f'(z)g_\alpha(z)^m\,dz
		=\frac1{2\pi i}\int_{|w|=R}w^{-m}F'(w)\,dw
		=m a_m.\qedhere
		\]
	\end{proof}
	
	\begin{remark}[Faber polynomials]\label{rem:faber-polynomials}
		By \eqref{eq:dickson-key-identity}, the $\mathcal T_{\ell,\alpha}$ are
		the \emph{monic Faber polynomials} of $\Sigma_\alpha$
		\cite[Section~3]{Liesen2001}: the polynomial parts at infinity of
		$g_\alpha(z)^{-\ell}$, the $\ell$th power of the root in
		\eqref{eq:recap-roots} that grows at infinity. In the coordinate $g=g_\alpha(z)$, the term
		$g^{-\ell}$ pairs only with $g^\ell$ in the expansion
		\eqref{eq:mode-expansion} of $H_\alpha$,
		while $(-\alpha g)^\ell$ is analytic at $g=0$; this biorthogonality
		is why $\mathcal T_{\ell,\alpha}$ isolates $G_\ell$.
	\end{remark}
	
	The explicit limit involves the deterministic means
	\begin{equation}\label{eq:chebyshev-means}
		m_{2r-1,\alpha}:=0,\qquad
		m_{2r,\alpha}:=-(-\alpha)^r+r\sigma_r^2\alpha^{2r-2},
		\qquad r\geq1,
	\end{equation}
	with powers of exponent zero interpreted as $1$.
	
	\begin{corollary}[Explicit limit and Dickson diagonalization]
		\label{cor:chebyshev-diagonalization}
		Fix $|\alpha|\leq1$.
		\begin{enumerate}[label=(\alph*)]
			\item For every $\ell\geq1$, the limit in Theorem~\ref{clt}
			satisfies
			\[
			\mathcal L_\alpha(\mathcal T_{\ell,\alpha})
			=m_{\ell,\alpha}+\ell(-\alpha)^{\ell-1}G_\ell.
			\]
			The limits for distinct $\ell$ are independent, and the same $G_\ell$ governs
			every parameter $\alpha$ in the joint limit.
			\item For $f$ analytic near $\Sigma_\alpha$,
			\begin{equation}\label{eq:statistic-chebyshev-full-law}
				\mathcal L_\alpha(f)
				=\sum_{\ell\geq1}m_{\ell,\alpha}f_{\ell,\alpha}
				+\sum_{\ell\geq1}\ell(-\alpha)^{\ell-1}f_{\ell,\alpha}G_\ell.
			\end{equation}
			The first series is the mean and converges absolutely; the
			second is $\mathcal L_\alpha(f)-\E\mathcal L_\alpha(f)$
			and converges in $L^2$. In particular,
			\[
			\E|\mathcal L_\alpha(f)-\E\mathcal L_\alpha(f)|^2
			=\sum_{\ell\geq1}\ell^2\sigma_\ell^2|\alpha|^{2\ell-2}|f_{\ell,\alpha}|^2;
			\]
			at $\alpha=0$ only the first mode remains.
		\end{enumerate}
	\end{corollary}
	
	\begin{proof}[Proof of Corollary~\ref{cor:chebyshev-diagonalization}]
		Write $u=g_\alpha(z)$. The Hermite bias is itself a derivative: using
		$u'=-u^2/(1+\alpha u^2)$,
		\begin{equation}\label{eq:hermite-bias-antiderivative}
			-\frac12\partial_z\log(1+\alpha u^2)
			=\frac{\alpha u^3}{(1+\alpha u^2)^2}=g^{(1)}_\alpha(z),
		\end{equation}
		with the logarithm given by its power series since $|\alpha u^2|<1$.
		We collect the random exponent \eqref{eq:mode-expansion} and the bias
		into the series
		\begin{equation}\label{eq:dickson-generating-limit}
			\widetilde H_\alpha(u)
			:=\sum_{\ell\geq1}a_{\ell,\alpha}G_\ell u^\ell
			-\frac12\sum_{\ell\geq1}a_{\ell,\alpha}^2\sigma_\ell^2u^{2\ell}
			-\frac12\log(1+\alpha u^2),
		\end{equation}
		so that $H_\alpha(z)-\frac12\log(1+\alpha g_\alpha(z)^2)
		=\widetilde H_\alpha(g_\alpha(z))$. Since
		$-\ell a_{\ell,\alpha}=\ell(-\alpha)^{\ell-1}$ by
		\eqref{eq:polymer-weights}, expanding $\log(1+\alpha u^2)$ gives
		\begin{equation}\label{eq:limiting-dickson-coefficient}
			-\ell[u^\ell]\widetilde H_\alpha(u)
			=m_{\ell,\alpha}+\ell(-\alpha)^{\ell-1}G_\ell,\qquad\ell\geq1.
		\end{equation}
		
		Choose an admissible $\Gamma$. By \eqref{eq:bias-contour-identity},
		\eqref{eq:hermite-bias-antiderivative}, and the definition of
		$\widetilde H_\alpha$,
		\[
		\mathcal L_\alpha(f)
		=\frac1{2\pi i}\int_\Gamma f(z)\bigl(g^{(1)}_\alpha(z)+H_\alpha'(z)\bigr)\,dz
		=\frac1{2\pi i}\int_\Gamma f(z)\,\partial_z\widetilde H_\alpha(g_\alpha(z))\,dz.
		\]
		For
		$R=\max_\Gamma|g_\alpha|$ we have $|\alpha|R<1$, since
		$|\tau_\alpha|<1$ on $\Omega_\alpha$ (Lemma~\ref{lem:quadratic-branches}),
		and $|\alpha|R^2<1$ by Lemma~\ref{lem:conformal-disk-coordinate}. At radius
		$R$ the deterministic series in $\widetilde H_\alpha$ converge
		absolutely, and the Gaussian series converges absolutely in $L^2$
		since $\|G_\ell\|_2=\sigma_\ell$ and
		$|a_{\ell,\alpha}|R^\ell=R(|\alpha|R)^{\ell-1}$. So integration by
		parts and termwise use of \eqref{eq:chebyshev-contour-extraction} give
		\[
		\mathcal L_\alpha(f)
		=-\frac1{2\pi i}\int_\Gamma f'(z)\widetilde H_\alpha(g_\alpha(z))\,dz
		=-\sum_{\ell\geq1}\ell f_{\ell,\alpha}[u^\ell]\widetilde H_\alpha(u),
		\]
		with the stated series convergence.
		By \eqref{eq:limiting-dickson-coefficient}, this gives (b), including
		its variance formula by independence of the $G_\ell$. For
		$f=\mathcal T_{\ell,\alpha}$ only coefficient $\ell$ survives, which
		gives the identity in (a); the same $G_\ell$ enter $H_\alpha$ for every
		$\alpha$, which gives the rest of (a).
	\end{proof}
	
	For example, $\mathcal T_{1,\alpha}(z)=z$ and
	$\mathcal T_{2,\alpha}(z)=z^2+2\alpha$, so
	\[
	\mathcal L_\alpha(z)=G_1,\qquad
	\mathcal L_\alpha(z^2+2\alpha)=\alpha+v-2\alpha G_2.
	\]
	The two terms in the even means \eqref{eq:chebyshev-means} come from
	the Hermite bias $B_\alpha$ and the Wick correction in
	$H_\alpha$, respectively. In particular,
	$m_{4,\alpha}=(\kappa-1)\alpha^2$; at $\alpha=-1$, standard GOE has mean
	$1$ in every positive even degree, whereas standard GUE has mean $0$.
	
	\begin{corollary}[Dickson statistics from the spectrum]\label{cor:dickson-from-spectrum}
		Take $U=\C$. For every $\ell\geq1$ and every compact $B\subset\C$,
		\[
		\sup_{\alpha\in B}\Bigl|\mathcal L_{n,\alpha}(\mathcal T_{\ell,\alpha})-m_{\ell,\alpha}
		-(-\alpha)^{\ell-1}\bigl(\mathcal L_{n,-1}(\mathcal T_{\ell,-1})-m_{\ell,-1}\bigr)\Bigr|
		\xrightarrow{\mathbb P}0.
		\]
	\end{corollary}
	
	\begin{proof}
		Let $D_n(\alpha)$ be the difference in the corollary.
		
		\emph{Step 1: pointwise convergence on the unit disk.} Fix
		$|\alpha|\leq1$. Theorem~\ref{clt}, applied jointly at $\alpha$ and
		$-1$ with $U=\C$, and Corollary~\ref{cor:chebyshev-diagonalization}(a)
		give
		\[
		\Bigl(\mathcal L_{n,\alpha}(\mathcal T_{\ell,\alpha})-m_{\ell,\alpha},\
		\mathcal L_{n,-1}(\mathcal T_{\ell,-1})-m_{\ell,-1}\Bigr)
		\Rightarrow\bigl(\ell(-\alpha)^{\ell-1}G_\ell,\ \ell G_\ell\bigr),
		\]
		with the same $G_\ell$ in both components. Hence $D_n(\alpha)$
		converges in law to $0$, and so in probability.
		
		\emph{Step 2: uniformity in $\alpha$.} We claim
		that $D_n$ is a polynomial in $\alpha$ of degree at most $\ell-1$.
		The polynomial $p_{n,\alpha}$ is monic. In the coefficient of
		$z^{n-m}$, each permutation fixes the $n-m$ positions where $z$ is
		chosen and has at least one cycle on the other $m$, so this
		coefficient has $\alpha$-degree at most $m-1$. By Newton's
		identities, the $k$th power sum of the zeros has $\alpha$-degree at
		most $k-1$ for $k\geq1$. Since $\mathcal T_{\ell,\alpha}(x)$ is a
		combination of the terms $\alpha^jx^{\ell-2j}$, the sum of
		$\mathcal T_{\ell,\alpha}$ over the zeros has $\alpha$-degree at most
		$\ell-1$. So do $n\int\mathcal T_{\ell,\alpha}\,d\mu_\alpha$, because
		$\int x^{2k}\,d\mu_\alpha$ is a multiple of $(-\alpha)^k$ and odd
		moments vanish, and $m_{\ell,\alpha}$, by \eqref{eq:chebyshev-means}.
		This proves the claim. Now fix distinct $\alpha_1,\dots,\alpha_\ell$
		in the unit disk, with Lagrange basis polynomials $\lambda_k$. Then
		$D_n=\sum_{k=1}^\ell D_n(\alpha_k)\lambda_k$, so
		\[
		\sup_{\alpha\in B}|D_n(\alpha)|
		\leq\max_{1\leq k\leq\ell}|D_n(\alpha_k)|
		\sum_{k=1}^\ell\sup_{\alpha\in B}|\lambda_k(\alpha)|
		\xrightarrow{\mathbb P}0.\qedhere
		\]
	\end{proof}
	
	Each Dickson fluctuation of the $\alpha$-zeros is asymptotically
	$(-\alpha)^{\ell-1}$ times the centered Chebyshev statistic of the
	eigenvalues of the same matrix, uniformly in $\alpha$ on compact sets,
	including the exterior regime $|\alpha|>1$ of
	Appendix~\ref{app:exterior-regime}.
	
	\begin{remark}[Comparison with the polymer sums]\label{rem:dickson-polymer-sum}
		For $|\alpha|\leq1$, the Dickson statistics are also asymptotically
		explicit functions of the entries:
		\[
		\mathcal L_{n,\alpha}(\mathcal T_{\ell,\alpha})-m_{\ell,\alpha}
		-\ell(-\alpha)^{\ell-1}S_{\ell,n}\xrightarrow{\mathbb P}0.
		\]
		Indeed, Step~1
		of the proof of Proposition~\ref{prop:wasserstein-convergence} and
		Lemma~\ref{lem:w2-weak} give $Q_n\Rightarrow Q$ in
		$\mathcal H(\mathcal V)$, and with $\phi(\alpha,z)=(\alpha,g_\alpha(z))$,
		continuity of $q\mapsto(q,q\circ\phi)$ and
		Proposition~\ref{prop:joint-polymer-approximation} give
		$(Q_n,F_n)\Rightarrow(Q,\Psi)$. By Corollary~\ref{cor:polymer-sums-general},
		$(-1)^\ell S_{\ell,n}$ and $(-1)^\ell G_\ell$ are the coefficients of
		$\alpha^{\ell-1}u^\ell$ in $Q_n$ and $Q$, while the proof of
		Theorem~\ref{clt} writes $\mathcal L_{n,\alpha}(\mathcal T_{\ell,\alpha})$,
		with probability tending to one, as a continuous functional of $F_n$
		plus a deterministic term, equal at $\Psi$ to
		$m_{\ell,\alpha}+\ell(-\alpha)^{\ell-1}G_\ell$. So the two quantities
		converge jointly to the same variable.
	\end{remark}
	
	Appendix~\ref{app:analytic-statistics} gives closed formulas for the
	covariance of the limiting field across parameters.
	
	\section{Determinantal approximation}\label{sec:determinantal-approximation}
	We now approximate $p_{n,\alpha}$ by a power of an ordinary
	determinant of the same matrix, evaluated at a transformed spectral
	parameter. Like Section~\ref{sec:further-consequences}, this section
	can be read directly after Sections~\ref{sec:functional-limits}
	and~\ref{sec:consequences}.
	
	\subsection{The approximation theorem.}\label{sec:determinantal-theorem}
	\emph{The idea.} Recall the disk coordinate
	$\tau_\alpha(z)=-\alpha g_\alpha(z)$ from
	\eqref{eq:tau-definition}. For $\alpha\neq0$, set
	\begin{equation}\label{eq:determinant-comparison-map}
		w_\alpha(z):=J_{-1}(\tau_\alpha(z))
		=\tau_\alpha(z)+\tau_\alpha(z)^{-1}.
	\end{equation}
	For $(\alpha,z)\in\mathcal D$ with $\alpha\neq0$, we have
	$0<|\tau_\alpha(z)|<1$. Since $J_{-1}$ inverts $\tau_{-1}=g_{-1}$ on
	the punctured unit disk (Lemma~\ref{lem:conformal-disk-coordinate}),
	$w_\alpha(z)\notin[-2,2]$ and
	\[
	\tau_{-1}(w_\alpha(z))=\tau_\alpha(z):
	\]
	$w_\alpha(z)$ is the point whose disk coordinate for $\alpha=-1$ equals
	$\tau=\tau_\alpha(z)$. Their limiting exponents are built from the same
	Gaussian series $S$ of \eqref{eq:common-gaussian-field}: by
	\eqref{eq:limiting-exponent},
	\[
	H_\alpha(z)=\frac{S(\tau)}{\alpha}-\frac{\mathcal C(\tau)}{2\alpha^2},
	\qquad
	\mathcal J(\tau):=H_{-1}(w_\alpha(z))=-S(\tau)-\tfrac12\mathcal C(\tau).
	\]
	Raising the normalized determinant $F_{n,-1}(w_\alpha(z))$ to the power
	$-1/\alpha$ therefore reproduces the Gaussian part $S(\tau)/\alpha$
	exactly, and leaves only a deterministic mismatch in the Wick terms:
	\begin{equation}\label{eq:determinant-mode-cancellation}
		H_\alpha(z)+\frac1\alpha\mathcal J(\tau_\alpha(z))
		=-\frac{1+\alpha}{2\alpha^2}\mathcal C(\tau_\alpha(z)).
	\end{equation}
	This suggests
	\[
	p_{n,\alpha}(z)\approx P_{n,\alpha}(z)
	\left(\frac{P_{n,-1}(w_\alpha(z))}{\det(w_\alpha(z)I-M_n)}\right)^{1/\alpha}
	e^{-(1+\alpha)\mathcal C(\tau_\alpha(z))/(2\alpha^2)}.
	\]
	The deterministic \emph{correction factor}
	$e^{-(1+\alpha)\mathcal C(\tau_\alpha)/(2\alpha^2)}$ equals $1$ at
	$\alpha=-1$.
	
	\emph{The logarithm.} To define the power $1/\alpha$ consistently, we
	take a logarithm of the normalized determinant in the disk
	coordinate. Let $\lambda_1,\ldots,\lambda_n$ be the eigenvalues of
	$M_n$, and $y_1,\ldots,y_n\in(-2,2)$ the zeros of $P_{n,-1}$. For
	$0<|t|<1$,
	\begin{equation}\label{eq:disk-normalized-determinant}
		F_{n,-1}(t+t^{-1})
		=\prod_{j=1}^n\frac{1-\lambda_jt+t^2}{1-y_jt+t^2}.
	\end{equation}
	Since $y_j\in(-2,2)$, each denominator vanishes only on the unit circle, so the right side
	is holomorphic on the unit disk and equals $1$ at $t=0$. Let
	\begin{equation}\label{eq:determinant-normalized-logarithm}
		\mathcal J_n(t):=\log\prod_{j=1}^n\frac{1-\lambda_jt+t^2}{1-y_jt+t^2},
		\qquad \mathcal J_n(0)=0,
	\end{equation}
	be its analytic logarithm on the largest centered zero-free disk
	contained in the unit disk.
	
	\begin{theorem}[Joint determinantal approximation]
		\label{thm:determinantal-approximation}
		For $(\alpha,z)\in\mathcal D$ with $\alpha\neq0$, define
		\begin{equation}\label{eq:alpha-determinant-power}
			\widehat p_{n,\alpha}(z)
			:=P_{n,\alpha}(z)
			\left(\frac{P_{n,-1}(w_\alpha(z))}
			{\det(w_\alpha(z)I-M_n)}\right)^{1/\alpha}
			e^{-(1+\alpha)\mathcal C(\tau_\alpha(z))/(2\alpha^2)},
		\end{equation}
		where the power means $\exp(-\mathcal J_n(\tau_\alpha(z))/\alpha)$,
		and extend $\widehat p_{n,\alpha}$ to $\alpha=0$ by continuity, that is, by
		\begin{equation}\label{eq:alpha-zero-determinant-approximation}
			\widehat p_{n,0}(z)
			=z^n\exp\!\left[-\frac{\operatorname{Tr}M_n}{z}
			-\frac{v}{2z^2}\right],\qquad z\neq0.
		\end{equation}
		Then, for
		every compact $\mathcal K\Subset\mathcal D$,
		\begin{equation}\label{eq:alpha-relative-approximation}
			\sup_{(\alpha,z)\in\mathcal K}
			\left|\frac{\widehat p_{n,\alpha}(z)}{p_{n,\alpha}(z)}-1\right|
			\xrightarrow{\mathbb P}0,
		\end{equation}
		where the supremum is $+\infty$ if $\mathcal J_n$ is undefined at some
		$\tau_\alpha(z)$ or $p_{n,\alpha}$ vanishes somewhere on $\mathcal K$;
		these events have probability tending to zero.
	\end{theorem}
	
	\begin{proof}
		\emph{Step 1: the determinant ratio.} By
		\eqref{eq:disk-normalized-determinant}, $t\mapsto F_{n,-1}(J_{-1}(t))$
		is holomorphic on the unit disk and equals $1$ at $t=0$, where
		$J_{-1}(0)=\infty$. By Corollary~\ref{lem:extension-at-infinity},
		jointly with $F_n$ on $\mathcal D$, it converges in law on compact
		subsets of the unit disk to $e^{\mathcal J(t)}$, where
		$\mathcal J(t)=-S(t)-\tfrac12\mathcal C(t)$. The
		limit is nonvanishing, so, for each $r<1$, the probability that
		$\mathcal J_n$ is defined on $|t|\leq r$ tends to one.
		
		\emph{Step 2: the exponent.} Let $\tau=\tau_\alpha(z)$. By definition,
		$\widehat p_{n,\alpha}(z)=P_{n,\alpha}(z)\exp E_n(\alpha,z)$ with
		$E_n:=-\mathcal J_n(\tau)/\alpha-(1+\alpha)\mathcal C(\tau)/(2\alpha^2)$.
		Since $-\tau/\alpha=g_\alpha(z)$,
		\begin{equation}\label{eq:determinant-removable-exponent}
			E_n(\alpha,z)
			=g_\alpha(z)\frac{\mathcal J_n(\tau)}{\tau}
			-\frac{1+\alpha}{2}g_\alpha(z)^2
			\frac{\mathcal C(\tau)}{\tau^2},
		\end{equation}
		which has no singular factor in $\alpha$. At $\tau=0$, the two quotients
		extend continuously with values $\mathcal J_n'(0)=-\operatorname{Tr}M_n$
		and $v$. The first follows by differentiating
		\eqref{eq:determinant-normalized-logarithm} at $t=0$:
		$\mathcal J_n'(0)=\sum_j(y_j-\lambda_j)=-\operatorname{Tr}M_n$, since
		the zeros $y_j$ of $P_{n,-1}$ are symmetric about $0$. Since $g_0(z)=1/z$, this shows that
		\eqref{eq:alpha-zero-determinant-approximation} is the continuous
		extension to $\alpha=0$.
		
		Fix a compact $\mathcal K\Subset\mathcal D$ and choose $r$ with
		$\sup_{\mathcal K}|\tau_\alpha(z)|<r<1$. For $\Phi$ holomorphic on
		the unit disk with $\Phi(0)=1$ and no zeros on $|t|\leq r$, its logarithm vanishing at zero satisfies
		\[
		\frac{\log\Phi(t)}{t}=\int_0^1\frac{\Phi'(st)}{\Phi(st)}\,ds,
		\]
		with value $\Phi'(0)$ at zero. If $\Phi_k\to\Phi$ locally uniformly
		with $\Phi_k(0)=1$, then $\Phi_k'/\Phi_k\to\Phi'/\Phi$ uniformly on
		this disk, and the integral gives the same convergence for the quotients.
		Step~1 and the continuous mapping theorem applied to \eqref{eq:determinant-removable-exponent}
		now give $(F_n,E_n)\Rightarrow(\Psi,H)$ in $C(\mathcal K)^2$, with the limiting exponent
		identified by \eqref{eq:determinant-mode-cancellation} and continuity at $\alpha=0$.
		
		\emph{Step 3: conclusion.} The quotient
		$\widehat p_{n,\alpha}/p_{n,\alpha}$ is $e^{E_n}/F_{n,\alpha}$.
		Exponentiating and dividing by the nonvanishing limit $e^{H_\alpha}$ gives
		\[
		\frac{\widehat p_{n,\alpha}(z)}{p_{n,\alpha}(z)}
		\Rightarrow1\qquad\text{in }C(\mathcal K).
		\]
		Convergence in distribution to a deterministic function is
		convergence in probability. The exceptional events in the statement
		have probability tending to zero: the logarithm is defined by Step~1,
		and $p_{n,\alpha}$ has no zero on $\mathcal K$ with probability tending
		to one by Theorem~\ref{gaussian}, since $e^{H_\alpha}$ has no
		zeros and $P_{n,\alpha}$ has none on $\mathcal K\subset\mathcal U$.
	\end{proof}
	
	The theorem gives no rate. Guarantees for a fixed $n$, or for
	parameters approaching the boundary of $\mathcal D$, require further
	estimates.
	
	\subsection{Special cases and computation.}\label{sec:determinantal-special-cases}
	At $\alpha=-1$, $\tau_{-1}=g_{-1}$ gives
	$w_{-1}(z)=z$ and the correction factor is $1$, so
	$\widehat p_{n,-1}(z)=\det(zI-M_n)$ exactly, wherever the logarithm is
	defined. At $\alpha=0$, the continuous extension
	\eqref{eq:alpha-zero-determinant-approximation} keeps the first two terms of
	$\log\bigl(z^{-n}p_{n,0}(z)\bigr)=\sum_i\log(1-\eta_i/(\sqrt n\,z))$,
	with $\frac1n\sum_i\eta_i^2$ replaced by its mean $v$.
	At $\alpha=1$, the power is $1$ and no logarithm is needed, so the
	theorem becomes a statement about the permanental polynomial.
	
	\begin{corollary}[Permanental polynomial]\label{cor:permanent-determinant}
		Let $w_1(z):=-g_1(z)-g_1(z)^{-1}=-\sqrt{z^2+4}$, where
		$\sqrt{z^2+4}$ is the branch on $\C\setminus i[-2,2]$ asymptotic to
		$z$ at infinity. For every compact
		$K\Subset\C\setminus i[-2,2]$,
		\begin{equation}\label{eq:permanent-determinant-approximation}
			\sup_{z\in K}\left|
			\frac{e^{-\mathcal C(g_1(z))}P_{n,1}(z)P_{n,-1}(w_1(z))}
			{\per(zI-M_n)\det(w_1(z)I-M_n)}-1\right|
			\xrightarrow{\mathbb P}0,
		\end{equation}
		with the supremum $+\infty$ if a denominator vanishes on $K$. For
		standard GOE and GUE, including their usual random diagonals,
		$e^{-\mathcal C(t)}$ is respectively $(1-t^2)^2$ and $1-t^2$.
	\end{corollary}
	
	\begin{proof}
		At $\alpha=1$, the exponent $(1+\alpha)/(2\alpha^2)$ equals $1$,
		$\tau_1=-g_1$ and $w_1$ is as stated, and $\mathcal C$ is even. By
		\eqref{eq:g-quadratic}, $z=g_1^{-1}-g_1$, so
		$(g_1+g_1^{-1})^2=z^2+4$, and $g_1(z)\sim z^{-1}$ fixes the branch. The
		power in \eqref{eq:alpha-determinant-power} is
		$\exp(-\mathcal J_n(\tau_1(z)))=P_{n,-1}(w_1(z))/\det(w_1(z)I-M_n)$
		wherever $\mathcal J_n$ is defined, and $p_{n,1}(z)=\per(zI-M_n)$.
		Apply Theorem~\ref{thm:determinantal-approximation} to
		$\{1\}\times K\subset\mathcal D$. The GOE and GUE values follow from
		\eqref{eq:permanent-determinant-correction} in
		Appendix~\ref{sec:variance-series}.
	\end{proof}
	
	For real $z>0$, the permanent of $zI-M_n$ is thus compared with the
	determinant at the real point $w_1(z)=-\sqrt{z^2+4}<-2$, outside the
	limiting spectrum $[-2,2]$. The corollary does not cover the unshifted
	permanent $\per(M_n)$, since $z=0$ lies on $\Sigma_1=i[-2,2]$.
	
	\begin{remark}[Computation]\label{rem:determinantal-computation}
		For fixed $(\alpha,z)$, $\widehat p_{n,\alpha}(z)$ can be computed in
		$O(n^3)$ operations from \eqref{eq:alpha-determinant-power}, which has
		three ingredients.
		\begin{enumerate}[label=(\roman*)]
			\item \emph{Expected polynomials.} $P_{n,\alpha}(z)$ and
			$P_{n,-1}(w)$ follow in $O(n)$ operations from the recurrence
			\eqref{eq:hermite-recurrence}.
			\item \emph{The logarithm $\mathcal J_n$.} On the disk where
			$\mathcal J_n$ is defined, it is the sum over $j$ of
			$\operatorname{Log}(1-\lambda_jt+t^2)-\operatorname{Log}(1-y_jt+t^2)$
			with principal logarithms; for real $x$, $1-xt+t^2$ takes no values
			in $(-\infty,0]$ on this disk, so each summand is holomorphic there.
			The eigenvalues $\lambda_j$ of $M_n$ cost $O(n^3)$ operations once,
			and the $y_j$ are the eigenvalues of the Jacobi matrix $J_n$
			(Section~\ref{sec:joint-bounds}); each further evaluation then
			costs $O(n)$.
			\item \emph{The correction factor.} $\mathcal C$ depends on the
			entry law only through $\theta$, $\kappa$ and $v$, by
			\eqref{eq:determinant-variance-closed-form}; if these are unknown,
			empirical means over the entries of $M_n$ are consistent
			substitutes.
		\end{enumerate}
	\end{remark}
	
\begin{figure}[tbp]
		\centering
		\includegraphics[width=0.9\textwidth]{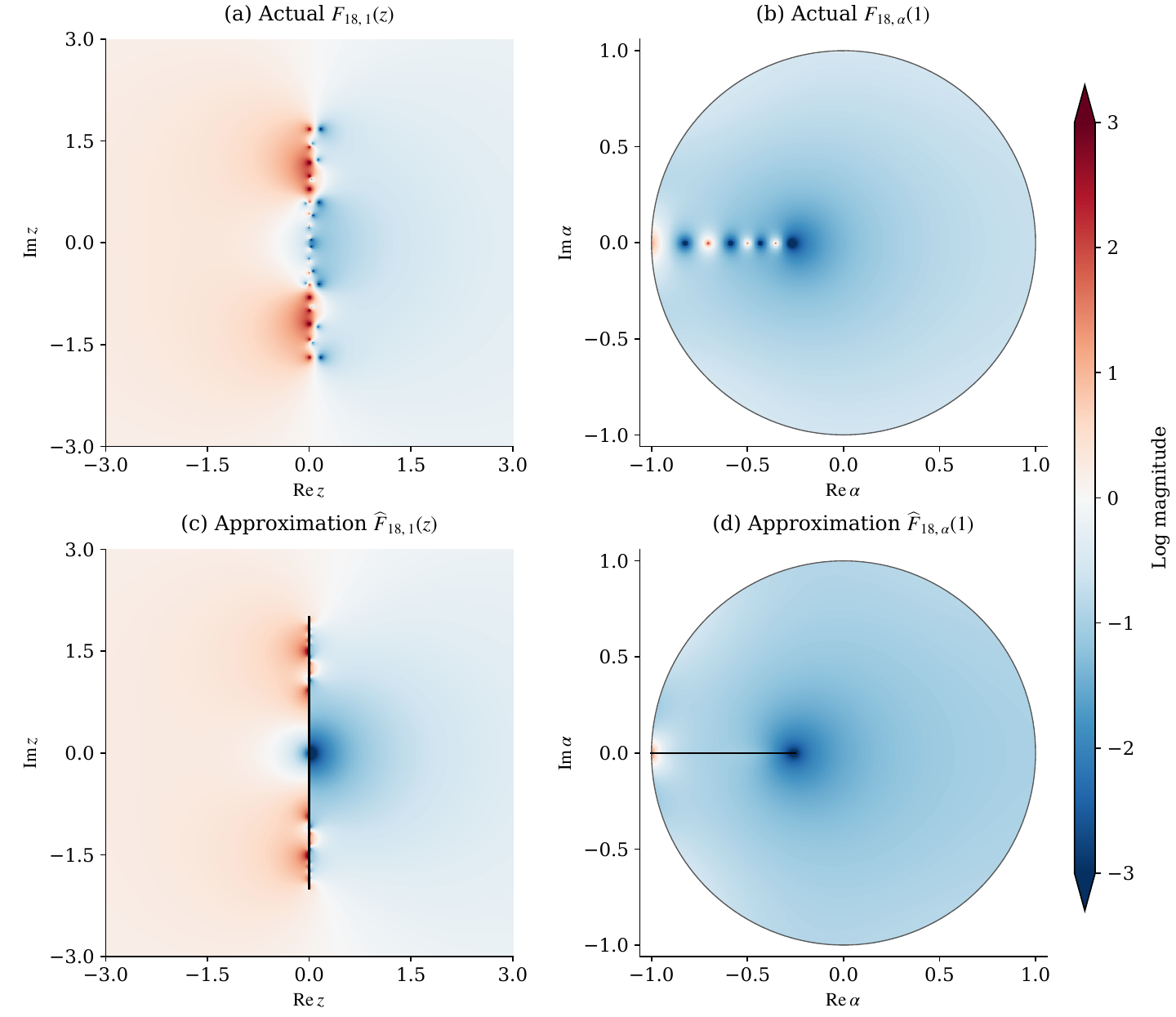}
		\caption{A single GUE sample at $n=18$: logarithmic magnitudes of
			the actual normalized functions $F_{18,\alpha}$ in (a,b) and their
			determinantal approximations $\widehat F_{18,\alpha}$ in (c,d).
			The left column fixes $\alpha=1$; the right fixes $z=1$ and varies
			$\alpha$ over the unit disk. All panels use the same color scale,
			saturating outside $[-3,3]$. Black segments in the bottom row mark
			parameters excluded from the approximation: $z\in i[-2,2]$ in (c)
			and $\alpha\in[-1,-1/4]$ in (d). Panels (a,b) include
			these segments wherever $P_{18,\alpha}(z)\neq0$.}
		\label{fig:determinantal-sample}
	\end{figure}

	Figure~\ref{fig:determinantal-sample} compares $F_{18,\alpha}$ with
	the normalized approximation
	$\widehat F_{18,\alpha}:=\widehat p_{18,\alpha}/P_{18,\alpha}$
	for one of the GUE matrices used in Figure~\ref{fig:limiting-correlation}.
	For this sample, $\operatorname{spec}(M_{18})\subset(-2,2)$, so
	$\mathcal J_{18}$ is defined throughout the unit disk.
	
	\subsection{Spectral determination.}\label{sec:spectral-determination}
	The approximation $\widehat p_{n,\alpha}(z)$ depends on $M_n$ only
	through its eigenvalues. Together with the uniform integrability
	supplied by Proposition~\ref{prop:wasserstein-convergence}, this shows
	that the normalized polynomial is asymptotically a function of the
	spectrum, with an error tending to zero in $L^2$.
	
	\begin{proposition}[Spectral determination]\label{prop:spectral-determination}
		For every compact $\mathcal K\Subset\mathcal D$,
		\[
		\E\sup_{(\alpha,z)\in\mathcal K}\bigl|F_{n,\alpha}(z)
		-\E[F_{n,\alpha}(z)\mid\operatorname{spec}M_n]\bigr|^2
		\longrightarrow0.
		\]
	\end{proposition}
	
	\begin{proof}
		Two draws that are conditionally independent given the spectrum share
		the same determinantal approximation, so they are close in
		probability; uniform integrability turns this into an $L^2$ bound.
		
		Write $\Lambda_n=\operatorname{spec}M_n$ and regard $F_n$ as a random
		element of $C(\mathcal K)$ with the supremum norm $\|\cdot\|$.
		Given $\Lambda_n$, draw $F_n^\ast$ independently from the conditional
		law of $F_n$. Both draws share the spectral approximation
		$\widehat F_n=\widehat p_{n,\alpha}/P_{n,\alpha}$, set to zero on
		$\mathcal K$ if the normalized logarithm is undefined anywhere on
		$\tau(\mathcal K)=\{\tau_\alpha(z):(\alpha,z)\in\mathcal K\}$.
		By Theorem~\ref{thm:determinantal-approximation} and tightness of
		$\|F_n\|$, we have $\|F_n-\widehat F_n\|\to0$ in probability.
		The same holds for $F_n^\ast$, since $(\Lambda_n,F_n^\ast)$ has the
		law of $(\Lambda_n,F_n)$. Hence $\|F_n-F_n^\ast\|\to0$ in probability.
		
		Proposition~\ref{prop:wasserstein-convergence} gives uniform
		integrability of $\|F_n\|^2$, and $F_n^\ast$ has the same law. Thus
		$\|F_n-F_n^\ast\|^2\leq2\|F_n\|^2+2\|F_n^\ast\|^2$ is uniformly
		integrable. Conditional independence gives
		$\E[F_n^\ast\mid F_n,\Lambda_n]=\E[F_n\mid\Lambda_n]$, so conditional
		Jensen yields
		\[
		\E\bigl\|F_n-\E[F_n\mid\Lambda_n]\bigr\|^2
		\leq\E\|F_n-F_n^\ast\|^2\longrightarrow0.\qedhere
		\]
	\end{proof}
	
	\begin{remark}[Basis dependence]\label{rem:spectral-concentration}
		The $\alpha$-determinant depends on the choice of basis, but
		Proposition~\ref{prop:spectral-determination} shows that its leading
		fluctuations do not. For instance, let $M_n$ be GOE (respectively
		GUE) and let $U_n$ be deterministic orthogonal (respectively unitary)
		matrices. Then $M_n$ and $U_nM_nU_n^*$ have the same law and the same
		spectrum, hence the same approximation $\widehat p_{n,\alpha}$, and
		Theorem~\ref{thm:determinantal-approximation} applied to both gives,
		for every compact $\mathcal K\Subset\mathcal D$,
		\[
		\sup_{(\alpha,z)\in\mathcal K}
		\left|\frac{\det_\alpha(zI-U_nM_nU_n^*)}
		{\det_\alpha(zI-M_n)}-1\right|
		\xrightarrow{\mathbb P}0.
		\]
	\end{remark}
	
	\appendix
	\section{Uniform almost-sure convergence of the zero measures}\label{app:almost-sure-semicircle}
	Exterior second moments control $\log^+|F_{n,\alpha}|$ almost surely.
	The mean-value inequality at $z=\infty$ turns this bound into
	$L^1_{\mathrm{loc}}$ convergence of $n^{-1}\log|F_{n,\alpha}|$.
	The following lemma then transfers convergence from the expected zeros
	to the random ones, proving Theorem~\ref{thm:almost-sure-semicircle}.
	
	\begin{lemma}[Zero measures from logarithmic $L^1$ convergence]
		\label{lem:log-ratio-zero-law}
		Let $\Sigma\subset\C$ be a compact line segment (or a point), and let
		$\mu$ be a probability measure supported on $\Sigma$. Let $p_j,P_j$ be
		monic polynomials of degree $n_j\to\infty$, with empirical zero measures
		$\nu_j$ and $\mu_j$. If $\mu_j\Rightarrow\mu$ and
		\[
		\frac1{n_j}\log\left|\frac{p_j}{P_j}\right|\longrightarrow0
		\qquad\text{in }L^1_{\mathrm{loc}}(\C\setminus\Sigma),
		\]
		then $\nu_j\Rightarrow\mu$.
	\end{lemma}
	
	Lemma~\ref{lem:log-ratio-zero-law} is a variant of the
	logarithmic-potential method for empirical measures; compare the
	replacement principle of \cite{TaoVuKrishnapur2010}. Here convergence
	is assumed only off $\Sigma$, which is why the proof uses test functions
	harmonic near $\Sigma$ together with polynomial density on the segment.
	
	\begin{proof}
		The distributional identity $\Delta\log|z-\zeta|=2\pi\delta_\zeta$
		\cite{Ransford1995,SaffTotik1997} gives, for every $\varphi\in C_c^\infty(\C)$ that is harmonic on a
		neighborhood of $\Sigma$,
		\[
		\int\varphi\,d(\nu_j-\mu_j)
		=\frac1{2\pi n_j}\int_{\C}\log\left|\frac{p_j(z)}{P_j(z)}\right|
		\Delta\varphi(z)\,dA(z)\longrightarrow0,
		\]
		because $\Delta\varphi$ has compact support in $\C\setminus\Sigma$.
		
		Given any neighborhood $U$ of $\Sigma$, choose a smooth cutoff
		$0\leq\chi\leq1$ supported in $U$ and equal to $1$ on a smaller
		neighborhood of $\Sigma$. Since $\chi$ is harmonic near $\Sigma$,
		$\int\chi\,d(\nu_j-\mu_j)\to0$. As $\mu_j\Rightarrow\mu$ and
		$\chi\equiv1$ on $\supp\mu\subseteq\Sigma$, we obtain
		$\int\chi\,d\nu_j\to\int\chi\,d\mu=1$. Thus $(\nu_j)$ is tight and
		every subsequential weak limit $\nu$ is supported on $\Sigma$.
		
		For each integer $k\geq0$, the function $\varphi_k(z):=\chi(z)z^k$ is
		also harmonic near $\Sigma$, so
		$\int\varphi_k\,d\nu_j-\int\varphi_k\,d\mu_j\to0$. Passing to the limit
		along a subsequence with $\nu_j\Rightarrow\nu$ yields
		$\int z^k\,d\nu(z)=\int z^k\,d\mu(z)$ for every $k\geq0$. Rotating
		$\Sigma$ to a real interval shows that polynomials are dense in
		$C(\Sigma)$ (constants suffice when $\Sigma$ is a point), so
		$\nu=\mu$. Hence $\nu_j\Rightarrow\mu$.
	\end{proof}
	
	\begin{proof}[Proof of Theorem~\ref{thm:almost-sure-semicircle}]
		\emph{Step 1: subexponential exterior growth.} For each compact
		$K\Subset\mathcal D$, Proposition~\ref{prop:second-moments} and
		Lemma~\ref{lem:sup-norm} give $\sup_n\E\|F_n\|_K^2<\infty$.
		Corollary~\ref{lem:extension-at-infinity}(b) extends this to every compact
		$K\Subset\widehat{\mathcal D}$, so
		$\mathbb P(\|F_n\|_K>n)\leq C_K/n^2$. Applying Borel--Cantelli on a
		countable compact exhaustion gives a single event of probability one on which
		\begin{equation}\label{eq:exterior-subexponential-bound}
			\frac1n\log^+\|F_n\|_K\longrightarrow0
			\qquad\text{for every compact }K\Subset\widehat{\mathcal D}.
		\end{equation}
		Fix an outcome in this event. We show that for every subsequence
		$n_j\to\infty$ and every sequence $|\alpha_j|\leq1$ with
		$\alpha_j\to a$, we have $\nu_{n_j,\alpha_j}\Rightarrow\mu_a$.
		
		\emph{Step 2: the logarithm in $L^1_{\mathrm{loc}}(\Omega_a)$.}
		Fix a compact $L\Subset\Omega_a$ and use the coordinate $u=g_a(z)$
		from Lemma~\ref{lem:conformal-disk-coordinate}. Choose $r>\sup_L|g_a|$
		with $|a|r^2<1$. The image $K_r=z_a(\{|u|\leq r\})$ is compact in
		$\Omega_a\cup\{\infty\}$. Since $|a|\leq1$ and $|a|r^2<1$ give
		$|a|r<1$, we have $\{a\}\times K_r\subset\widehat{\mathcal D}$. As $\widehat{\mathcal D}$ is open, the sets $\{\alpha_j\}\times K_r$
		lie in a common compact subset of $\widehat{\mathcal D}$ for all large $j$.
		
		Put $D_j(u)=F_{n_j,\alpha_j}(z_a(u))$. As $D_j(0)=1$, the mean-value
		inequality gives a nonnegative disk average for $\log|D_j|$.
		Using $|x|=2x^+-x$ and \eqref{eq:exterior-subexponential-bound} gives
		\[
		\frac{1}{\pi r^2 n_j}\int_{|u|<r}\bigl|\log|D_j(u)|\bigr|\,dA(u)
		\leq\frac{2}{n_j}\log^+\|D_j\|_{\{|u|\leq r\}}\longrightarrow0.
		\]
		The fixed change of variables $z=z_a(u)$ has bounded Jacobian on
		$g_a(L)\Subset\{0<|u|<r\}$. This gives convergence in $L^1(L)$
		for every such $L$, hence
		\begin{equation}\label{eq:exterior-log-L1}
			\frac1{n_j}\log|F_{n_j,\alpha_j}|\longrightarrow0
			\qquad\text{in }L^1_{\mathrm{loc}}(\Omega_a).
		\end{equation}
		
		\emph{Step 3: convergence and uniformity.}
		Lemma~\ref{lem:hermite-zero-law} gives $\mu_{n_j,-1}\Rightarrow\mu_\scal$
		on $[-2,2]$. Since $\mu_{n,-1}$ and $\mu_\scal$ are symmetric, the
		choice of square root does not matter, so we may choose square roots
		$\beta_j$ of $-\alpha_j$ converging to a square root of $-a$; scaling by $\beta_j$ then gives
		$\mu_{n_j,\alpha_j}\Rightarrow\mu_a$ and $\mu_{\alpha_j}\Rightarrow\mu_a$.
		Combining this with \eqref{eq:exterior-log-L1} and
		Lemma~\ref{lem:log-ratio-zero-law} (with $\Sigma=\Sigma_a$ and
		$\mu=\mu_a$) gives $\nu_{n_j,\alpha_j}\Rightarrow\mu_a$.
		
		Finally, if \eqref{eq:uniform-empirical-zero-law} failed for some
		bounded continuous $f$, we could select a subsequence $n_j\to\infty$
		and parameters $|\alpha_j|\leq1$ with discrepancy bounded below. Passing
		to a further subsequence with $\alpha_j\to a$ by compactness of the
		closed unit disk, both $\int f\,d\nu_{n_j,\alpha_j}$ and
		$\int f\,d\mu_{\alpha_j}$ converge to $\int f\,d\mu_a$, a
		contradiction.
	\end{proof}
	
	\begin{remark}[A weighted-norm route at fixed $\alpha$]\label{rem:weighted-norm-route}
		For fixed $0<|\alpha|\leq1$, the fixed-parameter zero law also follows
		from the weighted-norm criterion of Feng and Yao
		\cite[Theorem~4.4]{FengYao2026}. With
		$\|f\|_n^2:=\int_\R|f(x)|^2e^{-nx^2/2}\,dx$ and
		$\tilde p_{n,\alpha}(x):=\beta_\alpha^{-n}p_{n,\alpha}(\beta_\alpha x)$,
		the expected polynomials $\Her_r(\cdot\,;1/n)$ are orthogonal for this
		norm, and Lemma~\ref{lem:orthogonality} gives
		\begin{equation}\label{eq:integrated-second-moment}
			\frac{\E\|\tilde p_{n,\alpha}\|_n^2}{\|\Her_n(\cdot\,;1/n)\|_n^2}
			=\sum_{m=0}^n|\alpha|^{-m}c_m(\alpha).
		\end{equation}
		Evaluating \eqref{eq:degree-second-moment} at
		$s=n/((n+1)|\alpha|)$ bounds this by $C_\alpha(n+1)^2$, the polynomial
		growth that criterion requires.
	\end{remark}
	
	\section{The covariance kernel of the limiting field}\label{app:analytic-statistics}\label{sec:analytic-statistics}
	This appendix gives closed formulas for the covariance of the limiting
	field; none of them are used in the proofs of the main theorems. We
	express the individual cycle averages $\rho_\ell(\theta)$ in terms of
	Eulerian numbers. Splitting a permutation at its largest entry gives
	a differential equation for the generating function of weighted permutations.
	Solving it yields the cycle-average series and the covariance kernel
	across parameters for complex $\theta$. The kernel gives the
	correlations plotted in Figure~\ref{fig:limiting-correlation}, and the
	same computation gives the correction factor of
	Theorem~\ref{thm:determinantal-approximation} in closed form.
	
	\subsection{Cycle variance averages.}\label{sec:eulerian-cycle-average}
	Let $A(m,k)$ be the Eulerian number counting permutations of $m$
	elements with exactly $k$ ascents \cite[Chapter~1]{Petersen2015}.
	
	\begin{lemma}[Eulerian form of the cycle averages]\label{lem:eulerian-cycle-average}
		For $\ell\geq2$, the average in \eqref{eq:cycle-order-average} is
		\begin{equation}\label{eq:eulerian-cycle-average}
			\rho_\ell(\theta)=\frac1{(\ell-1)!}
			\sum_{j=1}^{\ell-1}A(\ell-1,j-1)
			\theta^j\overline\theta^{\,\ell-j}.
		\end{equation}
		In particular, $\rho_2(\theta)=|\theta|^2$.
	\end{lemma}
	
	\begin{proof}
		Rotate a cyclic ordering so that its smallest vertex comes first. Its
		remaining $\ell-1$ vertices form a permutation with $j-1$ ascents,
		where $j$ is the number of cyclic ascents: the first edge is upward
		and the last edge is downward. Averaging over these $(\ell-1)!$
		orderings proves \eqref{eq:eulerian-cycle-average}, initially with a
		real part in the sum. Reversing the ordering exchanges $j$ and
		$\ell-j$, so the sum is real and the real part can be omitted.
	\end{proof}
	
	For example,
	\[
	\rho_3(\theta)=|\theta|^2\operatorname{Re}\theta,\qquad
	\rho_4(\theta)=|\theta|^2
	\left((\operatorname{Re}\theta)^2
	+\frac{(\operatorname{Im}\theta)^2}{3}\right).
	\]
	The same Eulerian weights appear in the Chebyshev fluctuations of
	complex Wigner matrices; see
	\cite[equation~(1.4) and Section~6]{SchenkerSchulzBaldes2007}.
	
	\subsection{Variance series and correction factor.}\label{sec:variance-series}
	We now sum the cycle weights to obtain a closed form for the variance
	series $\mathcal C$, which gives the correction factor of
	Theorem~\ref{thm:determinantal-approximation} explicitly.
	
	Write $\theta_1:=\operatorname{Re}\theta$ and $\theta_2:=\operatorname{Im}\theta$, and
	set
	\begin{equation}\label{eq:theta-trigonometric-factor}
		\Delta_\theta(t):=
		\begin{cases}
			\cos(\theta_2t)-\dfrac{\theta_1}{\theta_2}\sin(\theta_2t),&\theta_2\neq0,\\
			1-\theta_1t,&\theta_2=0.
		\end{cases}
	\end{equation}
	The apparent singularity at $\theta_2=0$ is removable.
	
	\begin{lemma}[Cycle-average generating function]\label{lem:cycle-average-generating-function}
		The function $\Delta_\theta$ has no zeros in $|t|<1$, and its
		logarithm normalized at zero satisfies
		\begin{equation}\label{eq:cycle-average-generating-function}
			\sum_{\ell\geq2}\frac{\rho_\ell(\theta)}\ell t^\ell
			=-\log\Delta_\theta(t)-\theta_1t,\qquad |t|<1.
		\end{equation}
	\end{lemma}
	
	\begin{proof}
		Let
		\[
		\mathcal A(t):=\sum_{m\geq1}\frac{t^m}{m!}\sum_{k=0}^{m-1}A(m,k)\,
		\theta^k\overline\theta^{\,m-1-k}
		\]
		be the exponential generating function of nonempty permutations,
		with each ascent weighted by $\theta$ and each descent by
		$\overline\theta$. Split a permutation of $\{1,\ldots,m\}$ at the
		position of $m$. The entries on each side form a permutation of a
		subset, a nonempty left side adds an ascent into $m$, and a nonempty
		right side adds a descent out of $m$. Hence, as formal power series,
		\[
		\mathcal A'=(1+\theta\mathcal A)(1+\overline\theta\mathcal A),
		\qquad\mathcal A(0)=0.
		\]
		Let $R(t):=\sum_{\ell\geq2}\rho_\ell(\theta)t^\ell/\ell$, which is
		holomorphic on $|t|<1$ since $|\rho_\ell(\theta)|\leq|\theta|^\ell\leq1$.
		Writing $\theta^j\overline\theta^{\,\ell-j}
		=|\theta|^2\theta^{j-1}\overline\theta^{\,\ell-1-j}$ in
		\eqref{eq:eulerian-cycle-average} gives $R'=|\theta|^2\mathcal A$, and
		therefore
		\[
		R''=|\theta|^2\bigl(1+2\theta_1\mathcal A+|\theta|^2\mathcal A^2\bigr)
		=|\theta|^2+2\theta_1R'+(R')^2.
		\]
		Put $D(t):=e^{-R(t)-\theta_1t}$, which is holomorphic and zero-free on
		$|t|<1$. Then $D'=-(R'+\theta_1)D$, and
		\[
		D''=\bigl((R'+\theta_1)^2-R''\bigr)D=(\theta_1^2-|\theta|^2)D
		=-\theta_2^2D,
		\qquad D(0)=1,\quad D'(0)=-\theta_1.
		\]
		The unique solution of this initial value problem is $\Delta_\theta$,
		in both cases of \eqref{eq:theta-trigonometric-factor}. Hence
		$\Delta_\theta=D$ has no zeros in $|t|<1$, and
		$-\log\Delta_\theta=R+\theta_1t$, which is
		\eqref{eq:cycle-average-generating-function}.
	\end{proof}
	The lemma gives a closed form for the series
	\[
	\Phi(u):=\sum_{\ell\geq1}\sigma_\ell^2u^\ell,\qquad |u|<1,
	\]
	which yields both the variance series $\mathcal C(t)=\Phi(t^2)$, and
	hence the correction factor in
	Theorem~\ref{thm:determinantal-approximation}, and the covariance
	kernel of Section~\ref{sec:joint-covariance}.
	
	\begin{corollary}[Variance series]\label{cor:variance-series}
		With the logarithm of $\Delta_\theta$ normalized at zero,
		\begin{equation}\label{eq:determinant-variance-closed-form}
			\Phi(u)
			=-\log(1-u)-\log\Delta_\theta(u)
			+(v-1-\operatorname{Re}\theta)u
			+\frac{\kappa-1-|\theta|^2}{2}u^2,
			\qquad |u|<1.
		\end{equation}
		For real $\theta$, $\Delta_\theta(u)=1-\theta u$.
	\end{corollary}
	
	\begin{proof}
		Split off the terms $\sigma_1^2=v$ and $\sigma_2^2=\kappa/2$. For
		$\ell\geq3$, use $\sigma_\ell^2=(1+\rho_\ell(\theta))/\ell$,
		$\sum_{\ell\geq3}u^\ell/\ell=-\log(1-u)-u-u^2/2$,
		\eqref{eq:cycle-average-generating-function} and
		$\rho_2(\theta)=|\theta|^2$ from
		Lemma~\ref{lem:eulerian-cycle-average}.
	\end{proof}
	
	The correction factor of Theorem~\ref{thm:determinantal-approximation}
	is $e^{-(1+\alpha)\mathcal C(\tau_\alpha(z))/(2\alpha^2)}$, with
	$\mathcal C(t)=\Phi(t^2)$. At the permanent $\alpha=1$ it is
	$e^{-\mathcal C(g_1(z))}=e^{-\Phi(g_1(z)^2)}$, and
	\begin{equation}\label{eq:permanent-determinant-correction}
		e^{-\mathcal C(t)}
		=(1-t^2)\Delta_\theta(t^2)
		\exp\!\left((1+\operatorname{Re}\theta-v)t^2
		+\frac{1+|\theta|^2-\kappa}{2}t^4\right).
	\end{equation}
	
	On the diagonal, for $\alpha\neq0$, the mode expansion
	\eqref{eq:mode-expansion} gives
	$|a_{\ell,\alpha}g_\alpha(z)^\ell|^2=|\tau_\alpha(z)|^{2\ell}/|\alpha|^2$, so
	the limiting second moment is
	\[
	\E|\Psi_\alpha(z)|^2=\exp\bigl(\mathcal C(|\tau_\alpha(z)|)/|\alpha|^2\bigr).
	\]
	By \eqref{eq:determinant-variance-closed-form} with $u=r^2$,
	\[
	\mathcal C(r)=-\log(1-r^2)-\log\Delta_\theta(r^2)+O(1),
	\qquad r\uparrow1,
	\]
	which refines the lower bound in Remark~\ref{rem:second-moment-domain}.
	
	\subsection{Joint covariance and correlations.}\label{sec:joint-covariance}
	The common Gaussian sequence gives an explicit bilinear covariance kernel for
	the joint limit. In this subsection, write
	\[
	Y_\alpha:=H_\alpha-\E H_\alpha=\sum_{\ell\geq1}a_{\ell,\alpha}G_\ell g_\alpha^\ell
	\]
	for the centered exponent in \eqref{eq:mode-expansion}.
	
	\begin{proposition}[Bilinear covariance kernel]\label{prop:covariance-kernel}
		For points $s=(\alpha,z)$ and $s'=(\gamma,w)$ in $\mathcal D$, put
		$t:=\tau_\alpha(z)\tau_\gamma(w)$. If $\alpha\gamma\neq0$, define
		\begin{equation}\label{eq:joint-covariance-kernel}
			K(s,s'):=\frac{\Phi(t)}{\alpha\gamma},
		\end{equation}
		with $\Phi$ as in \eqref{eq:determinant-variance-closed-form}.
		Here $|t|<1$. The continuous extension when $\alpha\gamma=0$ is
		$K(s,s')=v g_\alpha(z)g_\gamma(w)$. Then
		\begin{equation}\label{eq:joint-covariance-identities}
			\begin{aligned}
				\E[Y_\alpha(z)Y_\gamma(w)]&=K(s,s'),\\
				\E[Y_\alpha(z)\overline{Y_\gamma(w)}]
				&=K(s,\overline{s'}),
			\end{aligned}
		\end{equation}
		where $\overline{s'}=(\overline\gamma,\overline w)$.
		Moreover,
		\begin{equation}\label{eq:exponent-covariance-centering}
			H_\alpha(z)
			=Y_\alpha(z)-\frac12K((\alpha,z),(\alpha,z)).
		\end{equation}
	\end{proposition}
	
	\begin{proof}
		Independence of the real Gaussian modes gives
		\[
		\E[Y_\alpha(z)Y_\gamma(w)]
		=\sum_{\ell\geq1}\sigma_\ell^2
		a_{\ell,\alpha}a_{\ell,\gamma}\bigl(g_\alpha(z)g_\gamma(w)\bigr)^\ell.
		\]
		For $\alpha\gamma\neq0$, $a_{\ell,\alpha}g_\alpha^\ell=\tau_\alpha^\ell/\alpha$
		makes the summand $\sigma_\ell^2t^\ell/(\alpha\gamma)$. The definition
		of $\mathcal D$ gives $|\tau_\alpha(z)|<1$ and $|\tau_\gamma(w)|<1$, so
		$|t|<1$, and the sum is $\Phi(t)/(\alpha\gamma)$. At $\alpha\gamma=0$,
		only the term $\ell=1$, namely $vg_\alpha(z)g_\gamma(w)$, survives,
		which is the continuous extension since $\Phi(t)=vt+O(t^2)$.
		The identity
		$g_{\overline\gamma}(\overline w)=\overline{g_\gamma(w)}$
		gives the Hermitian covariance, and the Wick correction in
		\eqref{eq:mode-expansion} gives
		\eqref{eq:exponent-covariance-centering}.
	\end{proof}
	
	For real $\theta$, $\Delta_\theta(t)=1-\theta t$, recovering the
	logarithmic kernel with cycle weights $\theta^\ell$.
	For standard GOE and GUE entries, including their usual random
	diagonals, the parameters $(\theta,\kappa,v)$ are $(1,2,2)$ and
	$(0,1,1)$, respectively. The kernel then simplifies to
	$-2(\alpha\gamma)^{-1}\log(1-t)$ and
	$-(\alpha\gamma)^{-1}\log(1-t)$, with the same continuous
	interpretation when $\alpha\gamma=0$.
	
\begin{figure}[tbp]
		\centering
		\includegraphics[width=0.9\textwidth]{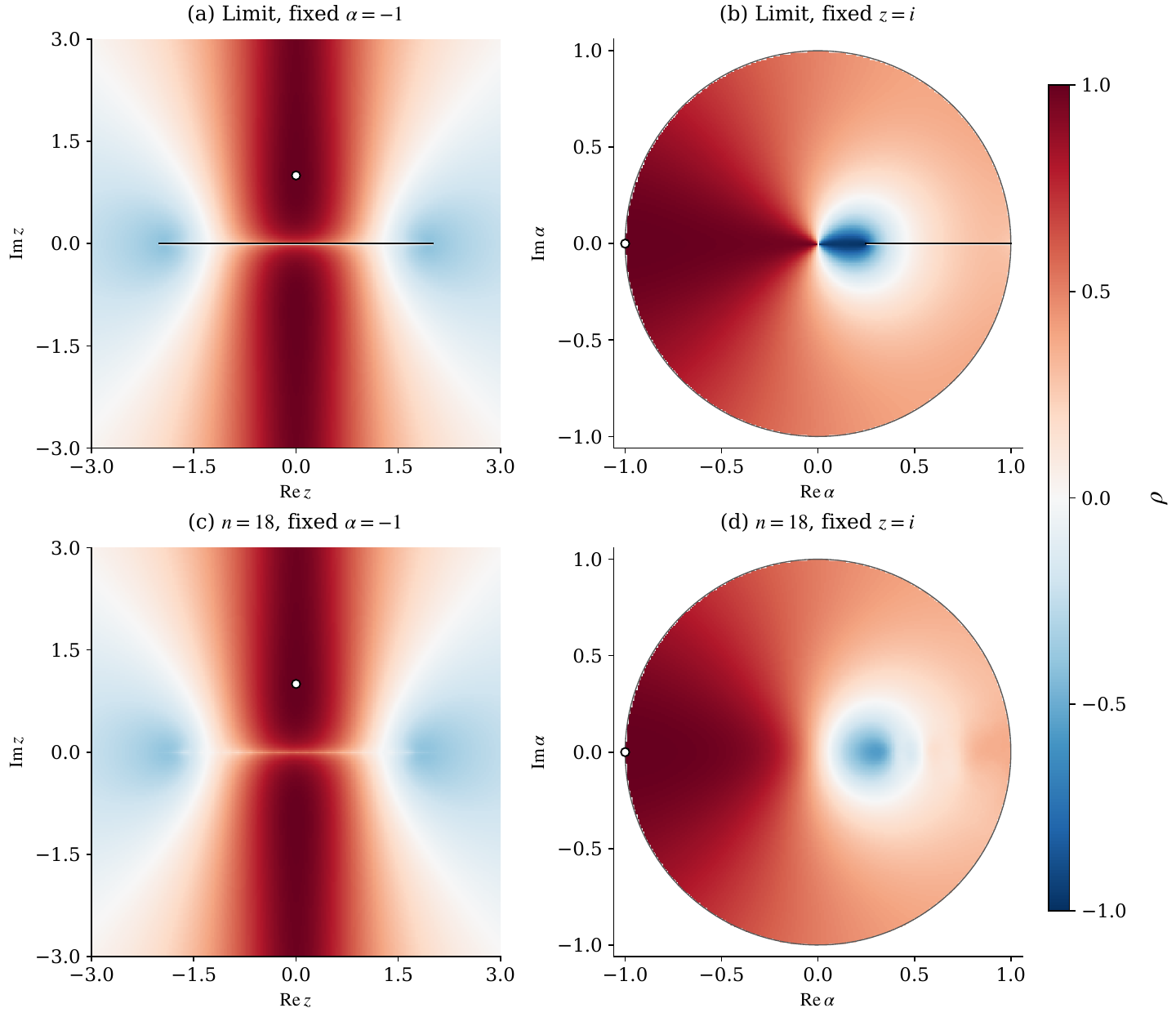}
		\caption{GUE correlations of $\log|\Psi_\alpha(z)|$ with
			$\log|\Psi_{-1}(i)|$ in (a,b), compared with empirical correlations
			of $\log|p_{18,\alpha}(z)|$ with $\log|p_{18,-1}(i)|$ in (c,d),
			using the same 4,096 GUE matrices as
			Figure~\ref{fig:limiting-correlation}. The left column fixes
			$\alpha=-1$; the right fixes $z=i$ and varies $\alpha$ over the unit
			disk. Open circles mark the reference point. All panels use the
			same color scale, with red and blue indicating positive and negative
			correlations. Black segments mark excluded limiting parameters:
			$z\in[-2,2]$ in (a) and $\alpha\in[1/4,1]$ in (b); the empirical
			panels include these segments. At $\alpha=0$ in (b), the limiting
			variance vanishes and the correlation is undefined; the empirical
			correlation in (d) is defined. GOE has the same limiting correlations.}
		\label{fig:determinant-correlation}
	\end{figure}

	The deterministic Wick correction does not affect the variance of
	$\log|\Psi_\alpha(z)|=\operatorname{Re}H_\alpha(z)$. With
	$s=(\alpha,z)$, the covariance identities give
	\begin{equation}\label{eq:log-modulus-variance}
		\Var\bigl(\log|\Psi_\alpha(z)|\bigr)
		=\frac12\operatorname{Re}\bigl[K(s,s)+K(s,\overline s)\bigr].
	\end{equation}
	The limiting GOE variance is
	twice the GUE variance, although their limiting correlations agree.
	Figure~\ref{fig:determinant-correlation} gives the corresponding
	correlations on slices through the determinant endpoint $\alpha=-1$,
	with reference point $(\alpha,z)=(-1,i)$.
	
	\section{The exterior regime for \texorpdfstring{$|\alpha|>1$}{|alpha| > 1}}\label{app:exterior-regime}
	For $|\alpha|>1$, our results give an asymptotic confinement region
	and identify holomorphic statistics of the zeros, but leave the
	limiting empirical distribution of the zeros open. The reference
	measure for these statistics is still the rotated semicircle law
	\[
	\mu_\alpha=(\beta_\alpha\id)_\#\mu_\scal,
	\qquad \beta_\alpha=\sqrt{-\alpha},
	\]
	supported on the segment $\Sigma_\alpha=\beta_\alpha[-2,2]$.
	The conclusions can be summarized as follows.
	\begin{itemize}
		\item \emph{Confinement.} With probability tending to one, all
		zeros lie in any prescribed neighborhood of the filled ellipse
		$\mathcal E_\alpha$ defined below.
		\item \emph{Holomorphic statistics.} Averages of functions
		holomorphic near $\mathcal E_\alpha$ converge to their integrals
		against $\mu_\alpha$. Equivalently, the empirical Cauchy transform
		converges to $g_\alpha$ outside the ellipse. We also identify
		the Gaussian fluctuations of these statistics.
		\item \emph{Determinantal approximation.}
		Theorem~\ref{thm:determinantal-approximation} is stated on all of
		$\mathcal D$, so it also holds outside $\mathcal E_\alpha$.
		\item \emph{Open question: zero law.} These conclusions do not show that
		the empirical zero measures converge, or determine
		whether the limiting support fills the ellipse. In particular,
		agreement with $\mu_\alpha$ for holomorphic statistics does not
		imply that the zeros follow the semicircle law on $\Sigma_\alpha$;
		the simulations in Figure~\ref{fig:zero-clouds} suggest that they do not.
	\end{itemize}
	We first describe the ellipse, then give the precise confinement and
	fluctuation statements obtained from Theorem~\ref{gaussian}, and
	finally explain why holomorphic statistics do not identify the zero
	law.
	
	\subsection{The ellipse geometry.}\label{sec:exterior-geometry}
	The domain of the functional limit at parameter $\alpha$ is
	\[
	\{z:(\alpha,z)\in\mathcal D\}
	=\{z\in\Omega_\alpha:|\tau_\alpha(z)|<1\}.
	\]
	For $|\alpha|>1$, define the filled ellipse
	\begin{equation}\label{eq:natural-ellipse}
		\mathcal E_\alpha:=
		\left\{\frac{\beta_\alpha}{\sqrt{|\alpha|}}
		\bigl((|\alpha|+1)x+i(|\alpha|-1)y\bigr):
		x,y\in\R,\ x^2+y^2\leq1\right\}.
	\end{equation}
	Its semiaxes are $|\alpha|+1$ and $|\alpha|-1$, with the major axis
	parallel to $\beta_\alpha=\sqrt{-\alpha}$. By
	Lemma~\ref{lem:conformal-disk-coordinate}, the rescaled Joukowski map
	$J_\alpha$ maps $0<|\tau|<|\alpha|^{1/2}$ conformally onto $\Omega_\alpha$.
	The condition $|\tau_\alpha(z)|<1$ therefore selects the image of the
	punctured unit disk. Writing $-\alpha=\beta_\alpha^2$, the unit circle
	maps to $J_\alpha(e^{i\vartheta})=e^{i\vartheta}+\beta_\alpha^2e^{-i\vartheta}$,
	which traces $\partial\mathcal E_\alpha$, so
	\[
	\{z:(\alpha,z)\in\mathcal D\}=\C\setminus\mathcal E_\alpha.
	\]
	The focal segment $\Sigma_\alpha$ lies strictly inside the ellipse.
	Put $\mathcal E_\alpha:=\Sigma_\alpha$ for $|\alpha|\leq1$, so that
	the displayed identity holds for every $\alpha$.
	At $|\alpha|=1$, the formula for $\mathcal E_\alpha$ reduces to
	$\Sigma_\alpha$. Since both descriptions depend continuously on
	$\alpha$ and agree on $|\alpha|=1$, the map
	$\alpha\mapsto\mathcal E_\alpha$ is Hausdorff continuous on all of $\C$.
	
	For $|\alpha|>1$, $1/g_\alpha$ is also the exterior conformal
	coordinate of $\mathcal E_\alpha$. Since $J_\alpha(t)=J_\alpha(-\alpha/t)$,
	the formula $J_\alpha(t)=t-\alpha/t$ maps $\{|t|>|\alpha|\}$
	conformally onto $\C\setminus\mathcal E_\alpha$, with inverse
	$1/g_\alpha$, the larger root of $t-\alpha/t=z$ in
	\eqref{eq:two-roots}. Since monic Faber polynomials depend only on
	this coordinate near infinity, the Dickson polynomials
	$\mathcal T_{\ell,\alpha}$ of Section~\ref{sec:dickson-statistics} are
	also the monic Faber polynomials of $\mathcal E_\alpha$, as for any
	ellipse confocal with $\Sigma_\alpha$.
	
	\subsection{Exterior confinement and statistics.}\label{sec:exterior-confinement}\label{sec:exterior-cauchy}\label{sec:exterior-statistics}
	For $|\alpha|>1$ and an open neighborhood $U$ of $\mathcal E_\alpha$,
	define $\mathcal L_{n,\alpha}(f)$ for $f\in\mathcal H(U)$ as in
	Section~\ref{sec:further-consequences}.
	
	\begin{corollary}[Exterior regime]\label{cor:exterior}
		\begin{enumerate}[label=(\alph*)]
			\item Let $B\subset\C$ be compact. For every $\varepsilon>0$, with
			probability tending to one, every zero $\zeta$ of $p_{n,\alpha}$
			satisfies $\dist(\zeta,\mathcal E_\alpha)<\varepsilon$,
			simultaneously for all $\alpha\in B$.
			\item For $|\alpha|>1$ and compact $K\Subset\C\setminus\mathcal E_\alpha$,
			the conclusions of Theorem~\ref{thm:cauchy-transform-fluctuations}
			hold on $\{\alpha\}\times K$.
			\item For $|\alpha|>1$ and $f_1,\ldots,f_r\in\mathcal H(U)$, the
			conclusions of Theorem~\ref{clt} and
			Corollary~\ref{cor:chebyshev-diagonalization} hold, provided the
			contour $\Gamma=\partial D$ encloses the whole filled ellipse,
			$\mathcal E_\alpha\subset D$, with $\overline D\Subset U$. In
			particular,
			$\mathcal L_\alpha(\mathcal T_{\ell,\alpha})
			=m_{\ell,\alpha}+\ell(-\alpha)^{\ell-1}G_\ell$.
		\end{enumerate}
	\end{corollary}
	
	\begin{proof}
		(a) Apply the proof of Theorem~\ref{zerofree} with this compact set
		$B$ and $\mathcal E_\alpha$ in place of $\Sigma_\alpha$. Hausdorff
		continuity of the ellipses makes
		$\{(\alpha,z):\alpha\in B,\
		z\in\widehat\C\setminus N_\varepsilon(\mathcal E_\alpha)\}$
		compact in $\widehat{\mathcal D}$, and
		Corollary~\ref{lem:extension-at-infinity} and nonvanishing of $\Psi$
		exclude zeros of $F_n$ there with probability tending to one. Since
		$\Sigma_\alpha\subset\mathcal E_\alpha$, the same holds for zeros of
		$p_{n,\alpha}$.
		
		(b) By Section~\ref{sec:exterior-geometry}, $\{\alpha\}\times K$ is a
		compact subset of $\mathcal D$, so
		Theorem~\ref{thm:cauchy-transform-fluctuations} applies there.
		
		(c) The proof of Theorem~\ref{clt} applies on $\{\alpha\}\times\Gamma$:
		part (a) places all zeros inside $\Gamma$ with probability tending to
		one, part (b) gives the fluctuations on $\Gamma$, and
		\eqref{eq:bias-contour-identity} identifies the bias integral as
		$B_\alpha(f)$. For Corollary~\ref{cor:chebyshev-diagonalization}, run
		its proof with $\Gamma=J_\alpha(\{|w|=R\})$ for some $R>|\alpha|$ close
		to $|\alpha|$, so that $\Gamma$ lies just outside $\mathcal E_\alpha$
		and its closed interior lies in $U$. Since $R>|\alpha|$, we have
		$g_\alpha(J_\alpha(w))=1/w$ on $|w|=R$ (Section~\ref{sec:exterior-geometry}), as used in the proof of
		Lemma~\ref{lem:dickson-biorthogonality}(b). Then
		\eqref{eq:chebyshev-contour-extraction} applies, and
		$\sup_\Gamma|\tau_\alpha|<1$ gives the same convergence estimates for
		the series.
	\end{proof}
	
	\subsection{Holomorphic statistics and the open zero law.}\label{sec:exterior-scope}
	Corollary~\ref{cor:exterior}(c) shows that, for every fixed
	$|\alpha|>1$ and function $f$ holomorphic on a neighborhood of
	$\mathcal E_\alpha$,
	\[
	\int f(z)\,d\nu_{n,\alpha}(z)-\int f(z)\,d\mu_\alpha(z)
	=O_{\mathbb P}(n^{-1}).
	\]
	In particular, all holomorphic moments approach those of
	$\mu_\alpha$. When $|\alpha|\leq1$, confinement to the segment
	$\Sigma_\alpha$ makes such moments sufficient to identify the zero
	law. For $|\alpha|>1$, the confinement region has interior, and
	holomorphic moments no longer determine a probability measure there.
	
	A concrete example is the uniform area measure on
	$\mathcal E_\alpha$, which has exactly the same holomorphic moments
	as the rotated semicircle law $\mu_\alpha$. To see this, put
	$r=|\alpha|$ and $u=\beta_\alpha/\sqrt r$, and let $W$ be uniform on
	the unit disk. Then $Z=u(rW+\overline W)$ is uniform on
	$\mathcal E_\alpha$. Rotational symmetry and
	$\E|W|^{2k}=1/(k+1)$ give
	\[
	\E Z^{2k}=\frac{1}{k+1}\binom{2k}{k}(-\alpha)^k
	=\int z^{2k}\,d\mu_\alpha(z),
	\qquad \E Z^{2k+1}=0.
	\]
	Consequently, these two distinct measures also have the same Cauchy
	transform outside $\mathcal E_\alpha$, so holomorphic limits cannot
	distinguish them.
	
	This is why our methods stop at the exterior. In the proof of
	Theorem~\ref{thm:almost-sure-semicircle}, uniqueness uses polynomial
	approximation on a segment, which fails on a filled ellipse.
	
\begin{figure}[tbp]
		\centering
		\includegraphics[width=0.9\textwidth]{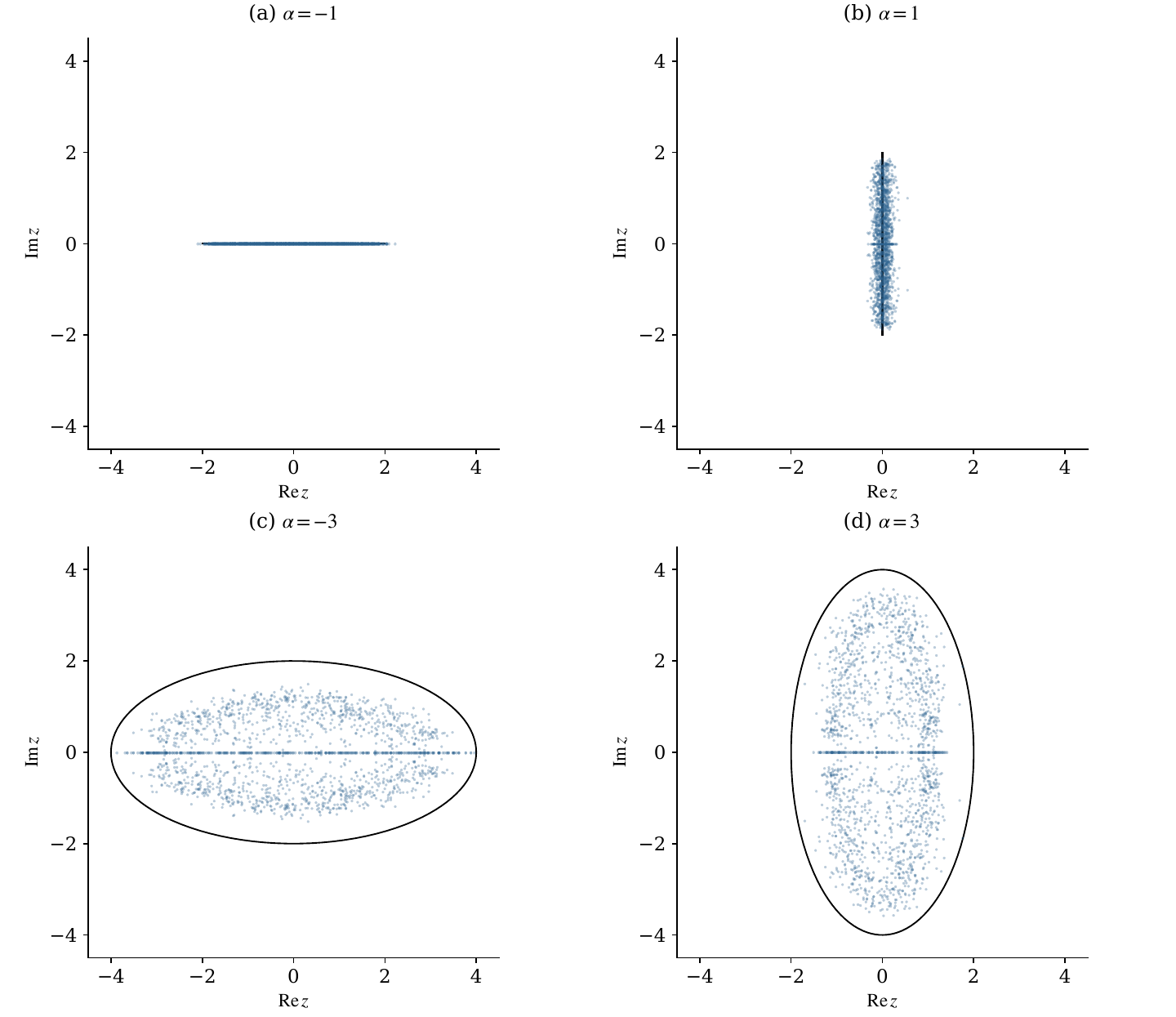}
		\caption{Zeros of $p_{18,\alpha}$ for $100$ independent GUE
			matrices, with the same matrices used in all four panels:
			(a) $\alpha=-1$, (b) $\alpha=1$, (c) $\alpha=-3$, and
			(d) $\alpha=3$. Each panel contains $1{,}800$ zeros and uses
			the same unscaled coordinate axes. Black segments in (a,b)
			mark $\Sigma_\alpha$; black curves in (c,d) mark
			$\partial\mathcal E_\alpha$, with semiaxes $4$ and $2$.
			The visible horizontal bands in (c) and (d) consist of real zeros.}
		\label{fig:zero-clouds}
	\end{figure}

	Figure~\ref{fig:zero-clouds} compares the zeros at $\alpha=\pm1$
	with those at $\alpha=\pm3$ for GUE matrices of size $n=18$.
	At this finite size, the latter exhibit a two-dimensional spread
	and a population of real zeros.

	\section{Cluster expansions and the trace method}\label{app:cluster-expansion}
	This appendix gives cluster-expansion and trace-method context for
	Steps~1--4 of Section~\ref{sec:heuristic}; it is not used in the
	proofs, whose estimates are in
	Sections~\ref{sec:expansion}--\ref{sec:fixed-order}.
	For simplicity, throughout this appendix we assume that the unscaled
	entry variables $\xi_{ij}$ and $\eta_i$ are bounded.
	
	The cycle expansion defines a \emph{raw gas}, which connects to the
	trace method, and its expectation, the \emph{mean gas}, carries the
	semicircle law (Sections~\ref{app:cluster-raw}--\ref{app:cluster-mean}).
	Normalizing by the mean gas and summing the attached mean dimers
	exactly gives the expansion of Proposition~\ref{prop:cycle-expansion};
	in the limit $A_{n,\alpha,m}\to g_\alpha^m$ it becomes the
	\emph{centered gas}
	(Sections~\ref{app:cluster-attachments}--\ref{app:cluster-resummed}),
	which is asymptotically ideal: one copy of a polymer gives a Gaussian
	mode and two copies give the Wick correction
	(Section~\ref{app:cluster-centered}).
	
	\subsection{The raw gas and the cluster expansion.}\label{app:cluster-raw}
	The cycle expansion assembles singletons, edges and simple cycles into
	vertex-disjoint families. Its pieces are the polymers of
	Section~\ref{sec:expansion}, but with uncentered edge activities, and
	compatible families again have disjoint vertex sets. An edge represents a transposition and is also called a
	\emph{dimer}. Let a formal variable $u$ mark occupied vertices. The raw
	activities of a singleton, an edge $e=\{i,j\}$ and a cycle $c$ of
	length $\ell\geq3$ are
	\begin{equation}\label{eq:raw-polymer-activities}
		w_{\{i\}}(u):=-uM_{ii},\qquad
		w_e(u):=\alpha u^2|M_{ij}|^2,\qquad
		w_c(u):=a_{\ell,\alpha}X_cu^\ell,
	\end{equation}
	suppressing the dependence on $\alpha$. Factoring out $z^n$ and grouping
	cycle orientations in the cycle decomposition
	(Lemma~\ref{lem:cycle-decomposition}) gives the exact
	polymer expansion
	\begin{equation}\label{eq:exact-polymer-gas}
		\Xi_{n,\alpha}(u)
		:=\sum_{\mathcal F\ \mathrm{compatible}}
		\prod_{p\in\mathcal F}w_p(u),
		\qquad
		p_{n,\alpha}(z)=z^n\,\Xi_{n,\alpha}(1/z).
	\end{equation}
	We call $\Xi_{n,\alpha}$, in which the empty family contributes $1$,
	the \emph{partition function} of the \emph{raw gas}.
	
	Its logarithm is organized by overlaps. Treat the activities as
	independent formal variables. Groups of polymers with disjoint vertex
	sets do not interact: the partition function factorizes over them, so
	its logarithm has no mixed terms. A single polymer, by contrast, can
	repeat: it contributes $\log(1+w_p)=w_p-w_p^2/2+\cdots$, although
	$\Xi_{n,\alpha}$ contains it at most once.
	
	To make this precise, let $\mathbf p=(p_1,\ldots,p_r)$ be an ordered
	tuple of polymers, possibly with repetitions. Its
	\emph{incompatibility graph} $I(\mathbf p)$ on $[r]$ joins $i\neq j$
	whenever $V(p_i)\cap V(p_j)\neq\varnothing$; repeated polymers are
	always joined. The tuple is a \emph{cluster} if $I(\mathbf p)$ is
	connected, and its \emph{Ursell coefficient} is
	\[
	\varphi(\mathbf p):=\sum_{A}(-1)^{|A|},
	\]
	summed over the edge sets $A\subseteq E(I(\mathbf p))$ for which
	$([r],A)$ is connected. In particular, $\varphi$ vanishes unless
	$\mathbf p$ is a cluster. Then, as formal power series,
	\begin{equation}\label{eq:formal-cluster-expansion}
		\log\Xi_{n,\alpha}(u)
		=\sum_{r\geq1}\frac1{r!}\sum_{\mathbf p}
		\varphi(\mathbf p)
		\prod_{j=1}^rw_{p_j}(u),
	\end{equation}
	where $1/r!$ compensates for ordering the tuple, and each
	$u$-coefficient is a finite sum \cite{KoteckyPreiss1986},
	\cite[Chapter~5]{FriedliVelenik2017}.
	
	\subsection{Tree clusters and the trace method.}\label{app:cluster-trees}
	At $\alpha=-1$, the raw partition function is a determinant, so
	\begin{equation}\label{eq:cluster-trace-identity}
		\log\Xi_{n,-1}(u)
		=\log\det(I-uM_n)
		=-\sum_{r\geq1}\frac{u^r}{r}\operatorname{Tr}M_n^r.
	\end{equation}
	The cluster expansion and the trace expansion therefore organize the
	same formal series in two ways: by overlapping polymers and
	by closed walks. In the trace-method proof of the semicircle law, the
	walks that survive in $\E\operatorname{Tr}M_n^{2k}$ as $n\to\infty$
	are those that traverse a tree, each edge once in each direction
	\cite[Section~12]{FleermannKirsch2022}. We show that the clusters made
	of the edges of a tree match these walks term by term;
	Section~\ref{app:cluster-mean} shows that, after expectation, they are
	also the leading clusters. Fix a labeled tree
	$T$ with $k\geq1$ edges, and let $d_v$ be its vertex degrees. The
	number $N(T)$ of closed walks traversing every edge once in each
	direction is
	\begin{equation}\label{eq:cluster-tree-walk-count}
		N(T)=2k\,w(T),\qquad
		w(T):=\prod_{v\in V(T)}(d_v-1)!.
	\end{equation}
	Indeed, such a walk is a depth-first traversal of $T$. For a starting
	vertex $v_0$, ordering the branches explored at each vertex gives
	$d_{v_0}!\prod_{v\neq v_0}(d_v-1)!$ walks. Summing over $v_0$ gives
	\eqref{eq:cluster-tree-walk-count}, since $\sum_vd_v=2k$.
	Every such walk has matrix weight $\prod_{e\in E(T)}|M_e|^2$. Thus the
	contribution of $T$ to \eqref{eq:cluster-trace-identity} is
	\begin{equation}\label{eq:cluster-tree-weight}
		-\frac{N(T)}{2k}u^{2k}
		\prod_{e\in E(T)}|M_e|^2
		=-w(T)u^{2k}
		\prod_{e\in E(T)}|M_e|^2.
	\end{equation}
	The cluster containing each edge of $T$ once as a dimer gives the same
	contribution, and its Ursell coefficient can be computed directly. Two
	dimers are incompatible exactly when they share a vertex, so the
	incompatibility graph of these $k$ dimers is the line graph of $T$: the
	dimers at a vertex $v$ form a clique $K_{d_v}$, and any two of these
	cliques share at most one dimer. They are therefore the blocks of the
	line graph, the sum defining $\varphi$ factorizes over them, and a clique $K_m$ contributes $(-1)^{m-1}(m-1)!$. Hence
	\begin{equation}\label{eq:cluster-tree-ursell}
		\varphi=\prod_{v\in V(T)}(-1)^{d_v-1}(d_v-1)!=(-1)^{k-1}w(T),
	\end{equation}
	since $\sum_v(d_v-1)=k-1$. The $k!$ orderings of the distinct dimers
	cancel the factor $1/k!$ in \eqref{eq:formal-cluster-expansion}, and
	the activities $w_e(u)=-u^2|M_e|^2$ at $\alpha=-1$ give
	\eqref{eq:cluster-tree-weight}. Thus the Ursell factor $(d_v-1)!$ at a
	vertex counts the orders in which a depth-first walk explores the
	branches there.
	
	\subsection{The mean gas.}\label{app:cluster-mean}
	To take expectations, split each edge weight into its mean and
	fluctuation as in \eqref{eq:cycle-split},
	$w_e(u)=\alpha u^2/n+a_{2,\alpha}X_eu^2$, and represent the two parts
	by a \emph{mean dimer} and a centered edge polymer. Both occupy the
	same two vertices, so they are incompatible and the partition function
	is unchanged. The other activities already have the centered form
	$a_{d(p),\alpha}X_pu^{d(p)}$. Hence
	\begin{equation}\label{eq:raw-gas-split}
		\Xi_{n,\alpha}(u)
		=\sum_{\mathcal F,\mathfrak m}
		\prod_{p\in\mathcal F}a_{d(p),\alpha}X_pu^{d(p)}
		\left(\frac{\alpha u^2}{n}\right)^{|\mathfrak m|},
	\end{equation}
	over compatible families $\mathcal F$ of centered polymers and matchings $\mathfrak m$
	of $[n]\setminus V(\mathcal F)$. Taking expectations removes every centered
	polymer, by independence and centering, and leaves a sum over
	matchings $\mathfrak m$ of $[n]$:
	\[
	\E\,\Xi_{n,\alpha}(u)
	=\sum_{\mathfrak m}\left(\frac{\alpha u^2}{n}\right)^{|\mathfrak m|}
	=\prod_{j=1}^n\bigl(1-u\beta_\alpha t_{n,j}\bigr).
	\]
	This deterministic gas is the \emph{mean gas}; the product form is
	$\E\,\Xi_{n,\alpha}(1/z)=z^{-n}P_{n,\alpha}(z)$ in terms of the
	zeros of $P_{n,\alpha}$. In its cluster expansion, a connected
	cluster of $k$ dimer occurrences on $m$ vertices contributes
	$O(n^{m-k})$ in total, since there are $O(n^m)$ choices of vertices
	and each dimer has weight of order $1/n$. Connectivity forces $m\leq k+1$, with
	equality exactly for trees of distinct dimers. These leading clusters
	are the expectations of the tree clusters of
	Section~\ref{app:cluster-trees}, since each tree has expected weight
	$n^{-k}$. At $\alpha=-1$, summing them gives the count
	$\E\operatorname{Tr}M_n^{2k}=nC_k+O(1)$ of the trace-method proof of
	the semicircle law \cite[Section~12]{FleermannKirsch2022}, where $C_k$
	is the $k$th Catalan number; for general $\alpha$, a $k$-dimer tree
	acquires the factor $(-\alpha)^k$. Coefficient by coefficient in $u$, therefore,
	\[
	\frac1n\log\E\,\Xi_{n,\alpha}(u)
	\longrightarrow\int\log(1-ux)\,d\mu_\alpha(x):
	\]
	the mean gas carries the semicircle law.
	
	\subsection{Normalization and mean-dimer attachments.}\label{app:cluster-attachments}
	The normalized polynomial
	$F_{n,\alpha}(z)=\Xi_{n,\alpha}(1/z)/\E\,\Xi_{n,\alpha}(1/z)$ divides
	out the mean gas, and with it the semicircle part. At the
	level of clusters, normalization subtracts the mean-gas
	expansion from the raw one:
	\[
	\log\frac{\Xi_{n,\alpha}(u)}{\E\,\Xi_{n,\alpha}(u)}
	=\log\Xi_{n,\alpha}(u)-\log\E\,\Xi_{n,\alpha}(u).
	\]
	This cancels every cluster made only of mean dimers. The remaining
	clusters contain centered polymers, but may still contain mean
	dimers attached to them.
	
	The simplest such cluster attaches a single mean dimer to a centered
	polymer, and already these clusters survive as $n\to\infty$. A
	centered cycle $c$ of length $\ell$ overlaps the $\ell(n-\ell)$ mean
	dimers with exactly one endpoint in $c$, each of weight $\alpha u^2/n$.
	The cluster formed by $c$ and one of them has Ursell coefficient $-1$,
	so the relative correction to $w_c(u)$ is
	\begin{equation}\label{eq:raw-dimer-attachment}
		-\frac{\ell(n-\ell)\alpha u^2}{n}
		\longrightarrow-\ell\alpha u^2
		\qquad\text{for fixed }\ell.
	\end{equation}
	Longer attachments contribute at the same order, so the raw expansion
	cannot be truncated at any fixed number of polymers: the attachments
	must be summed.
	
	\subsection{Matching resummation and its limit.}\label{app:cluster-resummed}
	The cluster expansion has identified the contributions that must be
	summed. We sum them in the finite partition function, before taking
	logarithms. For a compatible family of centered polymers on $m$
	vertices, the compatible mean dimers form a matching of the remaining
	vertices. At $u=1/z$, summing these matchings and dividing by the mean
	partition function replaces $u^m$ by $A_m:=A_{n,\alpha,m}(z)$. This is
	the exact expansion of Proposition~\ref{prop:cycle-expansion}. At
	finite $n$, the factor $A_{d(\mathcal F)}$ couples the centered
	polymers through their total degree.
	
	\emph{Asymptotic factorization.} For each fixed $m$,
	Section~\ref{sec:coefficients} proves $A_m\to g_\alpha(z)^m$; the
	cavity heuristic explains why. By
	\eqref{eq:matching-ratio-definition}, $A_m$ is the product of the
	one-step ratios $h_r:=\Her_{r-1}(z;-\alpha/n)/\Her_r(z;-\alpha/n)$ over $n-m<r\leq n$, each measuring
	the cost of deleting one vertex from the mean gas on $r$ vertices.
	That vertex is either unmatched, contributing $z$, or matched to one of
	the other $r-1$ vertices with weight $\alpha/n$. By the recurrence
	\eqref{eq:hermite-recurrence}, this gives
	$h_r=1/\bigl(z+\tfrac{(r-1)\alpha}{n}h_{r-1}\bigr)$. For $r$ close to
	$n$, deleting one more vertex should barely change the ratio, so
	$h_r\approx h_{r-1}\approx t$ with the cavity equation
	$t=1/(z+\alpha t)$. This is the quadratic \eqref{eq:g-quadratic}, and
	its relevant root is $g_\alpha(z)$; hence $A_m\approx g_\alpha(z)^m$.
	Lemma~\ref{lem:fixed-degree-limit} makes this rigorous: each one-step ratio
	is the Cauchy transform of the zeros of a smaller expected polynomial,
	which converges to $g_\alpha$.
	
	For a single centered polymer $p$ of degree $d$, set the other
	centered activities to zero and treat $x_p:=a_{d(p),\alpha}X_p$ as a
	formal variable. The normalized partition function is then
	$1+A_dx_p$, whose logarithm collects the raw clusters made
	of copies of $p$ and mean dimers:
	\begin{equation}\label{eq:cluster-single-polymer-resummation}
		\log(1+A_dx_p)
		\longrightarrow\log(1+g_\alpha(z)^dx_p)
	\end{equation}
	coefficient by coefficient in $x_p$. For a general compatible family, the
	limiting coefficient factorizes over its polymers:
	\[
	g_\alpha(z)^{d(\mathcal F)}
	=\prod_{p\in\mathcal F}g_\alpha(z)^{d(p)}.
	\]
	Thus, at fixed total degree, replacing $A_m$ by $g_\alpha(z)^m$
	removes the effective interaction generated by the mean dimers.
	The resulting \emph{centered gas} has individual activities
	$a_{d(p),\alpha}X_pg_\alpha(z)^{d(p)}$ and partition function
	$Q_{n,\alpha}(g_\alpha(z))$, while retaining the vertex-disjointness
	constraint.
	
	\begin{remark}[Rooted trees and matching resummation]\label{rem:cluster-rooted-trees}
		The limit $A_d\to g_\alpha(z)^d$
		has a trace-method meaning. Among the clusters of $\log(1+A_dx_p)$,
		the leading ones attach to $p$ a forest of mean-dimer trees, each
		rooted at a single vertex of $p$; repeated dimers, cycles, and trees
		joining two roots lose powers of $n$. The clique factorization of
		Section~\ref{app:cluster-trees} then shows that each root carries the
		generating function of rooted plane trees, with weight $-\alpha u^2$
		per edge:
		\[
		y:=\sum_{k\geq0}C_k(-\alpha u^2)^k,
		\qquad
		y\big|_{u=1/z}=zg_\alpha(z).
		\]
		The $d$ roots thus turn $u^d$ into $(uy)^d$, which is
		$g_\alpha(z)^d$ at $u=1/z$, and the one-edge term
		\eqref{eq:raw-dimer-attachment} is the first correction. In this
		sense the matching resummation absorbs the tree excursions of the
		trace method into the activity of each polymer. The rooted-tree
		recursion $y=1-\alpha u^2y^2$ is equivalent, at
		$u=1/z$, to the cavity equation $g=1/(z+\alpha g)$ above.
	\end{remark}
	
	\subsection{The centered gas and its leading clusters.}\label{app:cluster-centered}
	Write $\tilde w_p(u):=a_{d(p),\alpha}X_pu^{d(p)}$ for the activities of
	the centered gas. Its partition function is the polymer polynomial
	\eqref{eq:polymer-polynomial},
	\[
	Q_{n,\alpha}(u)
	=\sum_{\mathcal F\ \mathrm{compatible}}
	\prod_{p\in\mathcal F}\tilde w_p(u),
	\]
	and Proposition~\ref{prop:joint-polymer-approximation} gives
	$F_{n,\alpha}(z)-Q_{n,\alpha}(g_\alpha(z))\to0$ in $L^2$ of
	the supremum norm on compact subsets of $\mathcal D$.
	Allowing distinct polymers to share vertices gives the ideal-gas product
	\begin{equation}\label{eq:ideal-gas-product}
		\mathcal I_{n,\alpha}(u):=\prod_p(1+\tilde w_p(u)).
	\end{equation}
	Each polymer can still occur at most once. By the cluster identity
	\eqref{eq:formal-cluster-expansion}, the clusters consisting of copies
	of a single polymer $p$ sum to $\log(1+\tilde w_p)$. Thus
	\begin{equation}\label{eq:cluster-centered-split}
		\log Q_{n,\alpha}(u)
		=\log\mathcal I_{n,\alpha}(u)+\mathcal R^{\mathrm{mix}}_{n,\alpha}(u),
	\end{equation}
	where the remainder collects clusters containing distinct polymers.
	
	\emph{Mixed clusters.} A cluster of distinct overlapping polymers
	occupies fewer vertices than its total degree. In the squared
	degree-$m$ coefficient, Lemma~\ref{configuration} allows $m$ distinct
	vertices only for disjoint pairs of identical polymers, which a cluster
	containing two distinct overlapping polymers cannot provide. Thus at
	most $m-1$ distinct vertices occur, and $\mathcal R^{\mathrm{mix}}_{n,\alpha}$
	is $o_{\mathbb P}(1)$ at each fixed degree in $u$.
	
	\emph{Repeated polymers.} In
	$\log(1+\tilde w_p)=\tilde w_p-\tilde w_p^2/2+\cdots$, the linear terms
	sum to the polymer sums $S_{\ell,n}$. The quadratic terms involve
	$\sum_{p\in\mathcal P_{\ell,n}}X_p^2$, which has mean
	$(n)_\ell\sigma_\ell^2/n^\ell$ by \eqref{eq:cycle-variance-block} and
	concentrates. The higher powers are negligible, since each $|X_p|$ is
	of order $n^{-\ell/2}$. Hence, at each fixed degree in $u$,
	\begin{equation}\label{eq:cluster-centered-truncation}
		\log Q_{n,\alpha}(u)
		=\sum_p\tilde w_p(u)-\frac12\sum_p\tilde w_p(u)^2+o_{\mathbb P}(1),
	\end{equation}
	and Proposition~\ref{prop:polymer-sum-moments} identifies the limit as
	$\log Q_\alpha(u)$; Section~\ref{sec:exponential-limit} reaches the same
	conclusion by orthogonality, without taking logarithms. In the cluster
	expansion, the Wick correction thus comes from two copies of the same
	polymer. A related use of random cancellation in cluster expansions
	appears in \cite[Sections~2.3 and~3.2]{KoehlerLeung2026}.
	
	\emph{Gaussian modes.} The same vertex count explains why the polymer sums are asymptotically Gaussian: a joint cumulant of order at least three involves
	connected patterns of polymers sharing vertices, which lose
	powers of $n$, so only the pairings of identical polymers in
	Section~\ref{sec:power-sums} survive.
	
	\begin{remark}[The trace interpretation of the Gaussian modes]\label{rem:cluster-trace-modes}
		At $\alpha=-1$, the Gaussian modes have a familiar trace-method
		interpretation. The cycle contribution to the covariance comes from
		pairs of closed walks sharing a simple cycle, with tree excursions
		attached at its vertices; see
		\cite[Section~4.1, discussion following equation~(4.1)]{Kargin2015}.
		The shared cycle is our cycle polymer, while matching resummation
		absorbs the attached trees into the factor $g_{-1}(z)^\ell$, as
		explained in Remark~\ref{rem:cluster-rooted-trees}. The Chebyshev
		change of basis cancels contributions from shorter cycles with
		attached trees, isolating one Gaussian mode per degree; see
		\cite[Lemma~4.2 and Theorem~5.1]{SchenkerSchulzBaldes2007}.
		Our Dickson polynomials $\mathcal T_{\ell,\alpha}$ provide the
		corresponding mode isolation for general $\alpha$
		(Corollary~\ref{cor:chebyshev-diagonalization}).
	\end{remark}
	
	\makeatletter\let\sectionname\@empty\makeatother
	\section*{Acknowledgements and AI usage}
	At the time of writing \cite{KoehlerLeung2026}, the authors had also considered the setting in this paper. At that time they were aware that the terms in the raw cluster expansion (discussed in Appendix~\ref{app:cluster-expansion}) were related to the trees in the trace-method proof of the semicircle law, and had partial results proving non-sharp bounds on the zero-free region. The argument in this paper, which uses summation over matchings to get the centered polymer expansion and prove Theorem~\ref{gaussian}, was later discovered via significant contributions from GPT 5.5. Some editing and extensions were completed with help from GPT~5.6--6.1 and Claude Opus~5.5. The authors take responsibility for all errors. 
	
	We are grateful to Sidhanth Mohanty for a valuable conversation at an early stage in this project, and to Peter McCullagh for explaining \cite{McCullagh2012}.
	
	\bibliographystyle{plain}
	\bibliography{refs}
	
\end{document}